\documentclass[12pt]{amsart}
\usepackage{etex}

\usepackage[margin=1in]{geometry}

\usepackage{amscd,latexsym,amsthm,amsfonts,amssymb,amsmath,amsxtra,mathdots,tikz}
\usepackage[mathscr]{eucal}
\usepackage{xcolor}
\usepackage[normalem]{ulem}
\usepackage{hyperref}
\usepackage[all,cmtip]{xy}
\usepackage{enumitem}
\usepackage{tabularx}
\usepackage{textcomp}
\usepackage{caption}

\renewcommand\theequation{\thesection.\arabic{equation}}

\newcommand{\BA}{{\mathbb {A}}}

\newcommand{\BC}{{\mathbb {C}}}

\newcommand{\BZ}{{\mathbb {Z}}}

\newcommand{\CP}{{\mathcal {P}}}

\newcommand{\CS}{{\mathcal {S}}}

\newcommand{\Fg}{{\mathfrak {g}}}

\newcommand{\Fl}{{\mathfrak {l}}}

\newcommand{\GL}{{\mathrm{GL}}}

\newcommand{\GSp}{{\mathrm{GSp}}}

\newcommand{\GSO}{{\mathrm{GSO}}}
\newcommand{\GSpin}{{\mathrm{GSpin}}}
\newcommand{\Spin}{{\mathrm{Spin}}}
\newcommand{\PGSO}{{\mathrm{PGSO}}}

\newcommand{\PGL}{{\mathrm{PGL}}}

\newcommand{\SL}{{\mathrm{SL}}}

\newcommand{\GU}{{\mathrm{GU}}}
\newcommand{\HSpin}{{\mathrm{HSpin}}}

\newcommand{\SO}{{\mathrm{SO}}}

\newcommand{\Sym}{{\mathrm{Sym}}}

\newcommand{\Sp}{{\mathrm{Sp}}}

\newcommand{\pair}[1]{\langle {#1} \rangle}

\newcommand{\back}{\backslash}

\def\BA{{\mathbb A}}

\def\BC{{\mathbb C}}

\def\back{{\backslash}}

\newtheorem{thm}{Theorem}[section]

\newtheorem {conj}[thm]{Conjecture}

\newtheorem {ques/conj}[thm]{Question/Conjecture}

\newtheorem{defn}[thm]{Definition}
\newtheorem{rmk}[thm]{Remark}
\newtheorem{prpt}[thm]{Property}

\makeatletter

\newcommand{\Rmnum}[1]{\expandafter\@slowromancap\romannumeral #1@}
\makeatother

\begin{document}
\renewcommand{\theequation}{\arabic{equation}}
\numberwithin{equation}{section}

\title[Strongly tempered hyperspherical Hamiltonian spaces]{Strongly tempered hyperspherical Hamiltonian spaces}

\author{Zhengyu Mao}
\address{Department of Mathematics \& Computer Science\\
Rutgers University – Newark\\
Newark, NJ 07102, USA}
\email{zmao@rutgers.edu}

\author{Chen Wan}
\address{Department of Mathematics \& Computer Science\\
Rutgers University – Newark\\
Newark, NJ 07102, USA}
\email{chen.wan@rutgers.edu}

\author{Lei Zhang}
\address{Department of Mathematics\\
National University of Singapore, Singapore}
\email{matzhlei@nus.edu.sg}

\date{}



\begin{abstract}
In this paper, we give a complete list of strongly tempered hyperspherical Hamiltonian spaces. We show that the period integrals attached to the list 
contains many previously studied Rankin-Selberg integrals and period integrals, thus give a new conceptual understanding of these integrals. The list also proposes many new interesting period integrals to study.
\end{abstract}

\maketitle

\section{Introduction}

\subsection{BZSV Duality}

In \cite{BSV}, Ben-Zvi, Sakellaridis, and Venkatesh proposed a beautiful relative Langlands duality for hyperspherical Hamiltonian spaces (in this paper, we will call it BZSV duality). We briefly recall the datum in the duality. Throughout this paper, $k$ is a global field, $\BA=\BA_k$, $F$ is a local field, and $\psi$ is a non-trivial additive character of $\BA/k$ (resp. $F$) if we are in the global (resp. local) setting. Let $G$ be a split connected reductive group defined over $k$. In Section 3 of \cite{BSV}, Ben-Zvi, Sakellaridis, and Venkatesh defined a special category of G-Hamiltonian spaces called the hyperspherical G-Hamiltonian spaces. Moreover, they also showed that each hyperspherical G-Hamiltonian spaces is associated to a quadruple $\Delta=(G,H,\rho_H,\iota)$ where $H$ is a split reductive subgroup of $G$;  $\rho_H$  is a symplectic representation of $H$; and  $\iota$ is a homomorphism from $\SL_2$ into $G$ whose image commutes with $H$. For the rest of this paper, we will only discuss the quadruple instead of the Hamiltonian space associated to it (we will call such quadruple BZSV quadruple in this paper).

The BZSV duality concerns a pair of dual data $(\Delta,\hat{\Delta})$ where each side contains a BZSV quadruple: $\Delta=(G,H,\rho_H,\iota)$ and $\hat{\Delta}=(\hat{G},\hat{H}',\rho_{\hat{H}'},\hat{\iota}')$. The map $\iota$ induces an adjoint action of $H\times \SL_2$ on the Lie algebra $\Fg$  of $G$ and we can decompose it as
$$\oplus_{k\in I} \rho_k\otimes Sym^k$$
where $\rho_k$ is some representation of $H$ and $I$ is a finite subset of $\BZ_{\geq 0}$. We let $I_{odd}$ be the subset of $I$ containing all the odd numbers. In order for  $\Delta$ to be a BZSV quadruple, one of the (many) requirements is that the representation 
\begin{equation}\label{rho H iota}
\rho_{H,\iota}=\rho_H\oplus (\oplus_{i\in I_{odd}} \rho_i)
\end{equation}
 is a symplectic anomaly-free representation (see Section 5 of \cite{BSV}) of $H$. We refer the reader to \cite{BSV} for more details. Note that under BZSV duality, the group  $\hat{G}$ is the Langlands dual group of $G$ and $\hat{H}'=\hat{G}_{\Delta}$  can be viewed as the ``dual group" of the quadruple $\Delta$ (note that the groups $H$ and $\hat{H}'$ are not dual to each other in general, and the nilpotent orbits $\iota$ and $\hat{\iota}'$ are also not dual to each other in general). We recall the conjecture about period integrals in the BZSV duality.


  Let $\Delta=(G,H,\rho_H,\iota)$ and $\hat{\Delta}=(\hat{G},\hat{H}',\rho_{\hat{H}'},\hat{\iota}')$ be two quadruples that are dual to each other under the BZSV duality. We use $\rho_{H,\iota}$ and $\rho_{\hat{H}',\hat{\iota}'}$ to denote the symplectic anomaly-free representations associated to these quadruples. As we explained above, the maps $\iota$ and $\hat{\iota}'$ induce adjoint actions of $H\times \SL_2$ (resp. $\hat{H}'\times \SL_2$) on $\Fg$ (resp. $\hat{\Fg}$) and they can be decomposed as 
$$\Fg=\oplus_{k\in I} \rho_k\otimes Sym^k,\;\hat{\Fg}=\oplus_{k\in \hat{I}} \hat{\rho}_k\otimes Sym^k$$
where $\rho_k$ (resp. $\hat{\rho}_k$) are representations of $H$ (resp. $\hat{H}'$). It is clear that the adjoint representation of $H$ (resp. $\hat{H}'$) is a subrepresentation of $\rho_0$ (resp. $\hat{\rho}_0$). 

For an automorphic form $\phi$ of $G(\BA)$ (resp. $\hat{G}(\BA)$), we can define the period integral $\CP_{H,\iota,\rho_H}(\phi)$ (resp. $\CP_{\hat{H}',\hat{\iota}',\rho_{\hat{H}'}}(\phi)$) of it associated to the quadruple. Let's briefly recall the definition. We have a symplectic representation $\rho_{H,\iota}:H\rightarrow \Sp(V)$. Let $Y$ be a maximal isotropic subspace of $V$ and $\Omega_\psi$  be the Weil representation of $\widetilde{\Sp}(V)$ on the Schwartz space $\CS(Y(\BA))$. The anomaly free condition on $\rho_{H,\iota}$ ensures $\widetilde{\Sp}(V)$ splits over $Im(\rho_{H,\iota})$ and $\Omega_\psi$ restricts to a representation of $H(\BA)$ on $\CS(Y(\BA))$. We define the theta series
$$\Theta_{\psi}^{\varphi}(h)=\sum_{X\in Y(k)} \Omega_{\psi}(h)\varphi(X),\;h\in H(\BA),\varphi\in \CS(Y(\BA)),$$
and we can define the period integral to be
$$\CP_{H,\iota,\rho_H}(\phi,\varphi)=\int_{H(k)\back H(\BA)} \CP_\iota(\phi)(h)\Theta_{\psi}^{\varphi}(h)dh.$$
Here $\CP_\iota$ is the degenerate Whittaker period associated to $\iota$ (we refer the reader to Section 1.2 of \cite{MWZ} for its definition). To simplify the notation, we will omit the Schwartz function in the notion of the period and simply write it as $\CP_{H,\iota,\rho_H}(\phi)$ \footnote{when the nilpotent orbit associated to $\iota$ is not even, the degenerate Whittaker period $\CP_\iota$ is a Fourier-Jacobi coefficient and one also need to include an extra Schwartz function in its definition}. Similarly we can also define the period integral $\CP_{\hat{H}',\hat{\iota}',\rho_{\hat{H}'}}(\phi)$. The following conjecture is the main conjecture regarding global periods in BZSV duality.

\begin{conj}\label{BSV conj} (Ben-Zvi--Sakellaridis--Venkatesh, \cite{BSV})
\begin{enumerate}
\item Let $\pi$ be an irreducible discrete automorphic representation of $G(\BA)$ and let $\nu:\pi\rightarrow L^2(G(k)\back G(\BA))_{\pi}$ be an embedding. Then the period integral 
$$\CP_{H,\iota,\rho_H}(\phi),\;\phi\in Im(\nu)$$
is nonzero only if the Arthur parameter of $\pi$ factors through $\hat{\iota}':\hat{H}'(\BC)\times \SL_2(\BC)\rightarrow \hat{G}(\BC)$. If this is the case, $\pi$ is a lifting of a global tempered Arthur packet $\Pi$ of $H'(\BA)$ (the Langlands dual group of $\hat{H}'$). Then we can choose the embedding $\nu$ so that 
$$\frac{|\CP_{H,\iota,\rho_H}(\phi)|^2}{\pair{\phi,\phi}}``=" \frac{L(1/2,\Pi,\rho_{\hat{H}'})\cdot\Pi_{k\in \hat{I}}L(k/2+1,\Pi,\hat{\rho}_k)}{L(1,\Pi,Ad)^2},\;\phi\in Im(\nu).$$
Here $\pair{,}$ is the $L^2$-norm, and $``="$ means the equation holds up to some Dedekind zeta functions, some global constant determined by the component group of the global  L-packet associated to $\pi$, and some finite product over the ramified places  (including all the archimedean places).
\item Let $\pi$ be an irreducible discrete automorphic representation of $\hat{G}(\BA)$ and let $\nu:\pi\rightarrow L^2(\hat{G}(k)\back \hat{G}(\BA))_{\pi}$ be an embedding. Then the period integral 
$$\CP_{\hat{H}',\hat{\iota}',\rho_{\hat{H}'}}(\phi),\;\phi\in Im(\nu)$$ 
is nonzero only if the Arthur parameter of $\pi$ factors through $\iota:H(\BC)\times \SL_2(\BC)\rightarrow G(\BC)$. If this is the case,  $\pi$ is a lifting of a global tempered Arthur packet $\Pi$ of $\hat{H}(\BA)$ (the Langlands dual of $H$). Then we can choose the embedding $\nu$ so that 
$$\frac{|\CP_{\hat{H}',\hat{\iota}',\rho_{\hat{H}'}}(\phi)|^2}{\pair{\phi,\phi}}``=" \frac{L(1/2,\Pi,\rho_H)\cdot\Pi_{k\in I}L(k/2+1,\Pi,\rho_k)}{L(1,\Pi,Ad)^2},\;\phi\in Im(\nu).$$
\end{enumerate}
\end{conj}

\begin{rmk}\label{rmk BSV conj}
The above conjecture is usually called the Ichino-Ikeda type conjecture. To state an explicit identity instead of $``="$, one needs to make two adjustments on the right-hand side of the equation.
\begin{itemize}
\item In the ramified places, instead of using the local L-function, one needs to use the so-called local relative character defined by the (conjectural) Plancherel decomposition (see Section 17 of \cite{SV} and Section 9 of \cite{BSV}).
\item One also needs to add some Dedekind zeta functions on the right-hand side determined by the groups $G$ and $H$ (in all the known examples, those zeta functions are the $L$-function of the dual $M^{\vee}$ to the motive $M$ associated to $G,H$ introduced by Gross in \cite{Gross}), as well as some global constant determined by component group of the global  L-packet associated to $\pi$ (see Section 14.6.4 of \cite{BSV}) for these two quadruples.
\end{itemize}

\end{rmk}

\begin{rmk}\label{rmk dualtiy}
In \cite{BSV}, they also formulated many other conjectures for the duality (i.e., local/global geometric conjecture, local conjecture for Plancherel decomposition). The expectation is that those conjectures would uniquely determine the duality. In this paper we will only focus on their conjecture for period integrals. We also want to point out that given a general BZSV quadruple $\Delta=(G,H,\rho_H,\iota)$, at this moment there is no algorithm to compute the dual quadruple $\hat{\Delta}$. The only exception is for the so-called polarized case (i.e., when $\rho_H=0$) where the algorithm is given in Section 4 of \cite{BSV} (most quadruples considered in this paper are not polarized). As a result, given two BZSV quadruples $\Delta$ and $\hat{\Delta}$, at this moment one can only provide evidence for the duality between them by studying the various conjectures  (i.e., local/global geometric conjecture, local conjecture for Plancherel decomposition, global conjecture for period integrals) in \cite{BSV}.
\end{rmk}

\subsection{Strongly tempered BZSV quadruples}

\begin{defn}
We say the quadruple $\Delta=(G,H,\rho_H,\iota)$ is strongly tempered if $\hat{G}=\hat{H}'Z_{\hat{G}}$, i.e. the ``dual group" of $\Delta$ is equal to the dual group of $G$ up to center. We say the quadruple is reductive if $\iota$ is trivial.
\end{defn}

If the quadruple $\Delta=(G,H,\rho_H,\iota)$ is strongly tempered, then Conjecture 1.1(1) states that for all global tempered L-packet $\Pi$ of $G(\BA)$ \footnote{when $\hat{G}\neq \hat{H}'$, we need to make some assumptions on the central character of $\Pi$ so that its Langlands parameter factors through $\hat{H}'$}, there exists $\pi\in \Pi$ and $\nu:\pi\rightarrow L^2(G(k)\back G(\BA))_{\pi}$ such that 
\begin{equation}\label{periodst}
	\frac{|\CP_{H,\iota,\rho_H}(\phi)|^2}{\pair{\phi,\phi}}``=" \frac{L(1/2,\Pi,\rho_{\hat{H}'})}{L(1,\Pi,Ad)},\;\phi\in Im(\nu).
	\end{equation}
In other words, it means that the norm square of the period integral $\CP_{H,\iota,\rho_H}(\phi)$  is essentially equal to the central value of an automorphic L-function on every tempered global L-packet.

The most well-known example of strongly tempered quadruple is the Gross-Prasad model $(G,H,\rho_H,\iota)=(\SO_{2n+1}\times \SO_{2n},\SO_{2n},0,1)$. In this case the dual quadruple is given by 
$$(\hat{G},\hat{G},\hat{\rho},1)=(\Sp_{2n}\times \SO_{2n},\Sp_{2n}\times \SO_{2n},std_{\Sp_{2n}}\otimes std_{\SO_{2n}},1).$$
In this case, Conjecture \ref{BSV conj}(1) is just the Ichino-Ikeda conjecture in \cite{II} and Conjecture \ref{BSV conj}(2) is just the Rallis inner product formula for the theta correspondence between $\Sp_{2n}$ and $\SO_{2n}$.

\begin{rmk}
Conjecturally the quadruple is strongly tempered if and only if the integral
\begin{equation}\label{relchar}
	\int_{H(F)} \CP_{\iota}(\phi)(h)\varphi(h) dh
	\end{equation}
is absolutely convergent for all tempered matrix coefficient $\phi$ of $G(F)$.
Here $F=k_v$ is a local field for some $v\in |k|$,  $\CP_{\iota}$ is the local analogue of the global degenerate Whittaker period, and $\varphi(h)$ is a matrix coefficient of the local Weil representation of $H(F)$ associated to the symplectic representation $\rho_H$ (although the unipotent integral $\CP_\iota$ is not necessarily convergent and it needs to be regularized, see examples in \cite{Beu2, LM, Wal2, Wan2, WZ}). The local relative character in Remark \ref{rmk BSV conj} is given by the integral \eqref{relchar} where $\phi$ is the matrix coefficient of $\pi_v$; and $\pi_v$ is the local component of $\pi$ at $v$ which is a tempered representation of $G(F)$.
\end{rmk} 

In \cite{MWZ}, we proposed a relative trace formula comparison that relates the periods	$\CP_{H,\iota,\rho_H}(\phi)$ associated to any BZSV quadruple 
$(G,H,\rho_H,\iota)$ to the periods $\CP_{H_0,\iota_0,\rho_{H_0}}(\phi_0)$ associated to a strongly tempered BZSV quadruple $(G_0,H_0,\rho_{H_0},\iota_0)$.  Thus it is natural to consider Conjecture~\ref{BSV conj} first for the strongly tempered BZSV quadruples. In this paper we provide and study a complete list of strongly tempered BZSV quadruples (and hence a complete list of strongly tempered hyperspherical Hamiltonian spaces).

 By duality, in order to classify the strongly tempered quadruple $\Delta$, it is enough to classify its dual quadruple 
$$\hat{\Delta}=(\hat{G},\hat{H}',\hat{\rho},1).$$
Since $\hat{H}'Z_{\hat{G}}=\hat{G}$, it is enough to classify all the BZSV quadruples of the form
$$(\hat{G},\hat{G},\hat{\rho},1).$$

By \cite{BSV}, a quadruple $\hat{\Delta}=(\hat{G},\hat{G},\hat{\rho},1)$ is a BZSV quadruple if it satisfies the following three conditions.
\begin{enumerate}
\item The symplectic representation $\hat{\rho}$ is anomaly-free (see \cite[Section~5]{BSV}).
\item The symplectic representation $\hat{\rho}$ is multiplicity free.
\item The generic stabilizer of the representation $\hat{\rho}$ of $\hat{G}$ is connected.
\end{enumerate}
The set of multiplicity-free symplectic representations were classified by Knop \cite{K} and Losev \cite{Lo} independently. In this paper we will use the list in \cite{K}. By \cite[Theorem~2.3]{K}, the classification is reduced to that of symplectic representations that are {\it saturated} and multiplicity free, which are listed in Table 1, 2, 11, 12, 22, S of \cite{K}. In this paper we write down the strongly tempered quadruples that are (up to isogeny) the duals of $(\hat{G},\hat{G},\hat{\rho},1)$ when $\hat{\rho}$ is the symplectic representations listed in Knop's tables. In order to find the dual quadruple, we will provide a systematic way to write down $H$ and $\iota$ (see Property~\ref{main property}). On the other hand the choice of $\rho_H$ has been done in an ad hoc way at this moment.

\begin{rmk}
Condition (3) above is related to the Type N spherical root. Whenever this condition fails, we should expect some covering group to appear in the dual quadruple $\Delta=(G,H,\rho_H,\iota)$. This is not covered in BZSV's framework at this moment. Nonetheless, for some of the cases in \cite{K} that do not satisfy (3), we are still able to write down a candidate for the dual of the quadruple $\hat\Delta$ from some existing automorphic integrals in previous literatures. \footnote{In this paper, we will not check the connectedness condition for representations in \cite{K}, we will leave it as an exercise for the reader.}.
\end{rmk}

 
 \subsection{Statement of main results}
 
We first consider representations not in Table S of \cite{K} (because Table S of \cite{K} is an infinite table), i.e. consider all quadruples $\hat{\Delta}=(\hat{G},\hat{G},\hat{\rho},1)$ satisfy the following two conditions:

\begin{enumerate}
\item The symplectic representation $\hat{\rho}$ is anomaly-free.
\item The symplectic representation $\hat{\rho}$ appears in Table 1, 2, 11, 12, 22 of \cite{K}.
\end{enumerate}

For each of them, we will write down a quadruple $\Delta=(G,H,\rho_H,\iota)$ and claim it is dual to $\hat\Delta$ up to isogeny, or more precisely it is dual to $(\hat{G},\widehat{G/Z_{\Delta}}, \hat{\rho},1)$ where $Z_{\Delta}=Z_G\cap ker(\rho_H)$ and $Z_G$ is the center of $G$. To support the claim we provide evidence through the three main theorems below. Our results are summarized in the 6 tables at the end of this paper (Table \ref{red list}, \ref{red list extra}, \ref{non-red list 1}, \ref{non-red list 1 extra}, \ref{non-red list 2} and \ref{non-red list 2 extra}, the first two tables are for reductive  cases while the last four tables are for non-reductive  cases). 



\begin{thm}\label{main thm 1}
For all the reductive cases (Table \ref{red list} and \ref{red list extra}) except the quadruple $(\GL_6\times \GL_2,\GL_2\times S(\GL_4\times \GL_2),\wedge^2\otimes std_{\GL_2})$, and for all quadruples in Table \ref{non-red list 1} and \ref{non-red list 1 extra}, the local relative character of the period integral $\CP_{H,\rho_H,\iota}$ is equal to the $L-$value in Conjecture \ref{BSV conj}(1) at unramified places, namely equals $\frac{L(1/2,\Pi,\hat\rho)}{L(1,\Pi,Ad)}$ for the unramified representation $\Pi$.
\end{thm}
Recall that the local relative character at unramified places is defined in \eqref{relchar} with $\phi$ and $\varphi$ being unramified matrix coefficients normalized to be $1$ at identity, and with suitably chosen Haar measures. It is easy to check for all cases in Table \ref{red list} -- \ref{non-red list 2 extra}, the integral \eqref{relchar} is absolutely convergent.
\begin{rmk}
For the quadruple $(\GL_6\times \GL_2,\GL_2\times S(\GL_4\times \GL_2),\wedge^2\otimes std_{\GL_2})$ and for all  quadruples in Table \ref{non-red list 2} and \ref{non-red list 2 extra}, as far as we know, their local relative characters have not been computed at unramified places. Although we believe they can be computed by the same method as in \cite{II} and \cite{WZ}.
\end{rmk}

\begin{thm}\label{main thm 2}
	For the quadruples in Table \ref{red list}, \ref{non-red list 1} and \ref{non-red list 2},  Conjecture \ref{BSV conj}(2) holds, if we assume (when applicable) the global period integral conjectures in  \cite{GGP, GGP2,II} for Gan-Gross-Prasad models. 
\end{thm}

\begin{rmk}
In most cases for Theorem~\ref{main thm 2} and some cases for Theorem~\ref{main thm 1} we utilize the theta correspondence. We summarize the results needed for theta correspondence in Section~\ref{sec theta}. 
\end{rmk}

\begin{rmk}
In \cite{GGP2}, the authors only formulated a global conjecture regarding the non-vanishing of the period integrals for non-tempered Arthur L-packets (Conjecture 9.11 of \cite{GGP2}). An Ichino-Ikeda type conjecture for the period  is not available in \cite{GGP2} because of the difficulty in the definition of local relative character in the non-tempered case (see the last paragraph of Section 9 of \cite{GGP2}). Thus strictly speaking, for some cases in Theorem~\ref{main thm 2} we can only claim the nonvanishing part of Conjecture~\ref{BSV conj}(2). However the identity in Conjecture~\ref{BSV conj}(2) disregards the local factors at bad places, thus to prove it we only need an Ichino-Ikeda type conjecture without specifying the local factors at bad places. The formulation of such a conjecture is well known and we assume this version of the conjecture in Theorem~\ref{main thm 2}.
\end{rmk}


Beside the above two theorems, we provide one further evidence  for the duality for all the non-reductive quadruples. To state the evidence, we need to say a little more about the Hamiltonian space associated to the quadruple. Let $\Delta=(G,H,\rho_H,\iota)$ be a BZSV quadruple. Let $M$ be the centralizer of $\{\iota(diag(t,t^{-1}))|\;t\in \GL_1\}$ in $G$. It is easy to see that $M$ is a Levi of $G$ and $H\subset M$. We define
$$\Delta_{red}=(M,H,1,\rho_{H,\iota})$$
where the representation $\rho_{H,\iota}$ has been defined in \eqref{rho H iota}. The Hamiltonian $G$-space associated to $\Delta$ is defined by certain induction of the Hamiltonian $M$-space associated to $\Delta_{red}$ (see Section 3 of \cite{BSV} for details). In Section 4.2.2 of \cite{BSV}, they proposal a conjecture about the relation between the dual quadruples of $\Delta$ and $\Delta_{red}$. We will recall this conjecture in Conjecture \ref{Whittaker induction}. Now we are ready to state the third evidence.


\begin{thm}\label{main thm 3}
For any quadruple $\Delta=(G,H,\rho_H,\iota)$ in Table \ref{non-red list 1}, \ref{non-red list 1 extra}, \ref{non-red list 2} and \ref{non-red list 2 extra}, the corresponding quadruple $\Delta_{red}=(M,H,1,\rho_{H,\iota})$ is a quadruple in Table \ref{red list} and \ref{red list extra}. Moreover, the duality for the quadruples $\Delta$ and $\Delta_{red}$ \footnote{here by saying the duality for $\Delta$ (resp. $\Delta_{red}$) we mean the duality between $\Delta$ (resp. $\Delta_{red}$) and the quadruple $(\hat{G},\hat{G},1,\hat{\rho}$) (resp. $(\hat{M},\hat{M},1,\hat{\rho})$) where $\hat{\rho}$ is the corresponding symplectic representation in Table \ref{red list}-\ref{non-red list 2 extra} } is compatible with Conjecture \ref{Whittaker induction}.
\end{thm}

\begin{rmk}
Most of the quadruples in Table \ref{red list} and \ref{red list extra} come from Tables 1, 11, 2, 12, 22 of \cite{K}. There are some exceptions; the quadruples given in \eqref{GL(2) 5}, \eqref{GL(2) 4}, \eqref{GL(2) 3},  \eqref{GL(2) 31} and \eqref{GL(2) 1} are strongly tempered and dual to $\hat\rho$ from Table S in \cite{K}.
\end{rmk}
\begin{rmk}
For quadruples in Table \ref{non-red list 1}, \ref{non-red list 1 extra} and \ref{non-red list 2}, Theorem \ref{main thm 1} and \ref{main thm 2} already provide strong evidence for the duality of $(G,H,\rho_H,\iota)$. Combining with Theorem \ref{main thm 3}, we get strong evidence of Conjecture \ref{Whittaker induction} for quadruples in these three tables.
\end{rmk}

Lastly we consider Table S of \cite{K}. The representations coming out of this table are glued together from various representations of this table that already appeared in Table 1, 2, 11, 12, 22 of \cite{K}. Since the length can be arbitrary (i.e. we can glue any number of certain representations together), so this table produces infinitely many representations. In Section 9, for all the representations $\hat{\rho}$ coming from Table S that are anomaly-free and with connected generic stabilizer, we will describe a way to glue the dual quadruples which gives the dual of the quadruple $(\hat{G},\hat{G},\hat{\rho},1)$.

More precisely, given representations $(\hat{G}_i, \hat{\rho}_i)$ in Table S of \cite{K}, and let $(\hat{G},\hat{\rho})$ be the gluing of those representations. Assume that $\hat{\rho}$ is anomaly-free and its generic stabilizer is connected. We will describe the dual quadruple $\Delta$ of $\hat{\Delta}=(\hat{G},\hat{G},\hat{\rho},1)$ in terms of the dual quadruple $\Delta_i$ of $(\hat{G}_i, \hat{G}_i,\hat{\rho}_i,1)$. Roughly speaking, $\Delta$ is glued from $\Delta_i$ where the gluing process will be described in Section 9. To justify our construction, we will prove the following theorem.

 \begin{thm}\label{main thm S}
	With the notation above,	Conjecture \ref{BSV conj} for $(\Delta,\hat{\Delta})$ follows from Conjecture \ref{BSV conj} for $(\Delta_i,\hat{\Delta}_i)$.
	\end{thm}

In this paper, we provide the evidence of duality mainly through the period integral aspect, i.e., Conjecture \ref{BSV conj}. As we mentioned in Remark \ref{rmk dualtiy}, there are other ways to justify the duality, for example from the geometric conjectures (e.g. \cite{D, FW, BFGT,BFT,TY,FTY}) and local Plancherel conjectures (e.g. \cite{D}, \cite{FW}). We will not consider those conjectures in this paper. We just want to remark that Theorem~\ref{main thm 1} provides numerical evidence for the local Plancherel conjecture in Proposition 9.2.1 of \cite{BSV}, but we will not digress in these directions here.

\subsection{Rankin-Selberg integrals and special values of period integrals}

To end this introduction, we would like to point out that the list of strongly tempered quadruples we found in this paper recovers many existing integrals such as the Rankin-Selberg integrals in  \cite{BF90}, \cite{BG92a}, \cite{BG00}, \cite{BG92-ann}, \cite{G2}, \cite{Gin93}, \cite{G95}, \cite{GH04}, \cite{JPSS}, \cite{JS90}, \cite{PPS89}, \cite{PS} and the period integrals in \cite{GGP}, \cite{GJR01}, \cite{WZ}. It also produces many new interesting period integrals for studying.

A simple example that leads to a Rankin-Selberg integral is the quadruple \eqref{2.1 m=n}:
$$(\GL_n\times\GL_n,\GL_n,T(std_{\GL_n}),1)$$ which is dual to
$$(\GL_n\times\GL_n,\GL_n\times\GL_n, T(std_{\GL_n}\otimes std_{\GL_n}),1).$$
The attached period integral is
$$\int_{\GL_n(k)\back \GL_n(\BA)}\phi_1(g)\phi_2(g)\Theta^\Phi(g)\ dg$$
where $\phi_1\in\pi_1,\phi_2\in \pi_2$ are cusp forms in irreducible unitary cuspidal automorphic representations $\pi_1$ and $\pi_2$ on $\GL_n$ and $\Theta^\Phi(g)$  is a theta series on $\GL_n$ explicitly given by 
$$\Theta^\Phi(g)=|\det g|^{-\frac12}\sum_{\xi\in k^n}\Phi(\xi g).$$
Let $\xi_0=(0,0,\ldots,0,1)$, then we can identify $\Phi(g)$ with the sum of $|\det g|^{-\frac12}\Phi(0)$ and a mirabolic Eisenstein series 
$$E^\Phi(g)=|\det g|^{-\frac12}\sum_{\gamma\in P_0(k)\back \GL_n(k)}\Phi(\gamma g)$$
where $P_0$ is the mirabolic subgroup that fixes $\xi_0$. 
This period integral is just the specialization of the well-known Rankin-Selberg integral for tensor product $L-$function \cite{JPSS} evaluated at a specified value. 

The theory of Rankin-Selberg integrals is a very successful theory, producing many integral representations to study $L$-functions. A noted drawback of this theory is that the integrals are mostly developed in an ad hoc way. The list provided in this paper can actually fit many of the Rankin-Selberg integrals into the framework of BZSV duality. To be precise, those Rankin-Selberg integrals (evaluated at certain value) are simply the period integrals attached to some strongly tempered BZSV quadruples whose dual is closely related to the L-functions associated to the Rankin-Selberg integrals. The following is a list of such Rankin-Selberg integrals.

\begin{itemize}
	\item Integrals for exterior square $L-$functions by Bump-Friedberg \cite{BF90}.
	\item Integrals for Spin $L-$function by Bump-Ginzburg \cite{BG92a}, \cite{BG00} and \cite{G95}.
	\item  Integrals for standard $L-$functions of exceptional groups $E_6$  by Ginzburg \cite{G2}.
	\item Multivariable Rankin-Selberg integrals by Ginzburg-Hundley \cite{GH04} and Pollack-Shah  \cite{PS}.
	\item Rankin-Selberg convolution by Jacquet-Piatetski-Shapiro-Shalika \cite{JPSS}.
	\item Integrals for exterior square $L-$functions by Jacquet-Shalika \cite{JS90}.
\end{itemize}
The above list exhausts all currently known Rankin-Selberg integrals utilizing the mirobolic Eisenstein series. There are also examples above that use the Eisenstein series of other types (e.g., the ones in \cite{GH04} and \cite{PS}).
	
	 Our list provides more candidates for Rankin-Selberg integrals. For example, Model 12 of Table \ref{non-red list 2 extra} suggests considering the following Rankin-Selberg integral of $G=\GSO_8$, which should produce the standard L-function and the Half-Spin L-function. Let $\pi$ be a generic cuspidal automorphic representation of $\GSO_8(\BA)$, $\phi\in \pi$ and $P=MN$ be a maximal parabolic subgroup $\GSO_8$ with its Levi subgroup $M=\GL_2\times \GSO_4$. Let $H=S(\GL_2\times \GSO_4)$ be a subgroup of $M$ and let $E(h,s_1,s_2)$ be an automorphic function on $H$ induced from the trivial function on $\GL_2$ and the Borel Eisenstein series of $\GSO_4$ ($s_1,s_2$ are the parameter of the Eisenstein series). It is easy to see that one can take a Fourier-Jacobi coefficient of $\phi$ along the unipotent subgroup $N$ that produces an automorphic function on $H$. We will denote it by $\CP_{N}(\phi)$. Then, the integral associated to Model 12 of Table \ref{non-red list 2 extra} is just
$$\int_{H(k)\back H(\BA)/Z_G(\BA)}\CP_N(\phi)(h)E(h,s_1,s_2) dh.$$
In the spirit of Conjecture \ref{BSV conj}, we expect this to be the integral representation of the L-function $L(s_1,\pi,\rho_1)L(s_2,\pi,\rho_2)$ where $\rho_1$ (resp. $\rho_2$) is the standard representation (resp. Half-Spin representation) of $\Spin_8(\BC)$.

Meanwhile the majority of the quadruples in our list have period integrals that cannot be considered as specializations of Rankin-Selberg integrals. In some cases, the identities between the periods and the $L-$values in Conjecture~\ref{BSV conj} are consequences of Gan-Gross-Prasad conjectures \cite{GGP, GGP2,II}) and the Conjectures in \cite{WZ}. There is also one case where the integral is predicted by the work of Ginzburg-Jiang-Rallis \cite{GJR01} on the central value of symmetric cube $L-$functions.  Of more interest are the many cases where the conjectured identity in Conjecture~\ref{BSV conj} is new and unrelated to the conjectures mentioned above. For example each of the quadruple in tables \ref{non-red list 2} and \ref{non-red list 2 extra} gives such a new conjecture. 

We now list one example from Table \ref{red list extra} that not only provides a new Ichino-Ikeda type conjecture for a strongly tempered quadruple but also can be used to explain the Rankin-Selberg in \cite{GH04}. The example is Model 3 of Table \ref{red list extra}. The quadruple is reductive and is given by
$$\Delta=(G,H,\rho_H)=(\GSp_4\times \GSpin_8\times \GL_2,S(\GSpin_8\times G(\Sp_4\times \SL_2)), std_{\Sp_4}\otimes std_{\Spin_8}\oplus \HSpin_8\otimes std_{\SL_2}).$$
Let $\pi$ be a cuspidal generic automorphic representation of $G(\BA)$, $\phi\in \pi$ and $\Theta_{\rho_H}$ be the theta series associated to the symplectic representation $\rho_H$. Then the period integral is given by
$$\CP_{\Delta}(\phi)=\int_{H(k)\back H(\BA)/Z_{\Delta}(\BA)}\phi(h)\Theta_{\rho_H}(h)dh.$$
In the spirit of Conjecture \ref{BSV conj}, we expect the square of this period integral to be equal to 
$$\frac{L(1/2,\Pi,\hat{\rho})}{L(1,\Pi,Ad)}$$
where $\hat{\rho}$ is the representation $std_{\Sp_4}\otimes std_{\Spin_8}\oplus \HSpin_8\otimes std_{\SL_2}$ of $\widehat{G/Z_{\Delta}}(\BC)$. This is a new period integral that has not been considered before. If we replace the cusp form on $\GSp_4$ and $\GL_2$ by Borel Eisenstein series, then the period integral $\CP_{\Delta}$ becomes the Rankin-Selberg integral in \cite{GH04}.

	\begin{rmk} In this paper we will also encounter some representations $\hat\rho$ whose generic stabilizer is not connected. While these representations do not fit
		in the current framework of BZSV quadruple, one can still consider the associated period integrals and they are  related to the previously studied   integrals on covering groups, in \cite{BG92-ann, Gin93, GJR01, GRS, PPS89,  Tk14}.
	\end{rmk}



\subsection{Organization of the paper}
In Section 2, we will explain our strategy for writing down the dual quadruple. In Sections 3-7, we will consider Tables 1, 2, 11, 12, and 22 of \cite{K}. In Section 8 we summarize our findings in six tables. In Section 9 we will discuss Table S of \cite{K}.

\subsection{Acknowledgement} We thank Yiannis Sakellaridis, Akshay Venkatesh and Hiraku Nakajima for many helpful discussions. We thank Friedrich Knop for answering our question for some cases in \cite{K}. The work of the first author is partially supported by the Simons Collaboration Grant. The second author's work is partially supported by the NSF grant DMS-2103720, DMS-2349836 and a Simons Travel Grant. 
The work of the third author is partially supported by AcRF Tier 1 grants 	A-0004274-00-00 and A-0004279-00-00 of the National University of Singapore.

\section{Our strategy}

\subsection{Notation and convention}
In this paper, for a group $G$ of Type $A_n$ (resp. $B_n$, $C_n$, $D_n$, $G_2$, $E_6$, $E_7$), we use $std_G$ to denote the $n$-dimensional (resp. $2n+1$-dimensional, $2n$-dimensional, $2n$-dimensional, 7-dimensional, 27-dimensional, 56-dimensional) standard representation of $G$. We use $\Spin_{2n}$ (resp. $\Spin_{2n+1}$) to denote the Spin representation of the reductive group of Type $D_n$ (resp. $B_n$) and we use $\HSpin_{2n}$ to denote the Half-Spin representation of reductive group with Type $D_n$. We use $Sym^n$ (resp. $\wedge^n$) to denote the n-th symmetric power (resp. exterior power) of a reductive group of Type $A$. We use $\wedge_{0}^{3}$ to denote the third fundamental representation of a reductive group of Type $C_3$. Lastly, for a representation $\rho$ of $G$, we use $\rho^\vee$ to denote the dual representation and $T(\rho)$ to denote $\rho\oplus \rho^\vee$.

In this paper, we always use $l$ to denote the similitude character of a similitude group. If we have two similitude group $GH_1$ and $GH_2$, we let
$$G(H_1\times H_2)=\{(h_1,h_2)\in GH_1\times GH_2|\;l(h_1)=l(h_2)\},$$
$$S(GH_1\times GH_2)=\{(h_1,h_2)\in GH_1\times GH_2|\;l(h_1)l(h_2)=1\}.$$
Similarly we can also define $G(H_1\times \cdots \times H_n)$ and $S(GH_1\times \cdots \times GH_n)$. For example, 
$$S(\GL_{2}^{3})=S(\GL_2\times \GL_2\times \GL_2)=\{(h_1,h_2,h_3)\in \GL_{2}^{3}|\;\det(h_1h_2h_3)=1\}.$$

All the nilpotent orbits considered in this paper are principal in a Levi subgroup (this is also the case in \cite{BSV}). As a result, we will use the Levi subgroup or just the root type of the Levi subgroup to denote the nilpotent orbit (the zero nilpotent orbit is denoted by $1$). For a split reductive group $G$, we will use $T_G$ to denote a maximal split torus of $G$ (a minimal Levi subgroup).

For a BZSV quadruple $\hat{\Delta}=(\hat{G},\hat{G},\hat{\rho},1)$, there are many other quadruples that is essentially equal to $\hat{\Delta}$ up to some central isogeny. To be specific, one can take any group $\hat{H}$ of the same root Type as $\hat{G}$ such that the representation $\hat{\rho}$ can also be defined on $\hat{H}$. Then one can choose any group $\hat{G}'$ containing $\hat{H}$ such that $\hat{G}'=\hat{H}Z_{\hat{G}'}$. The quadruple $(\hat{G}',\hat{H},\hat{\rho}, 1)$ is essentially equal to $\hat{\Delta}$ up to some central isogeny. For example, both $(\PGL_{2}^3,\PGL_2,0,1)$ and $(\GL_{2}^{3},\GL_2,0,1)$ can be viewed as trilinear $\GL_2$-model. The dual quadruple of them are $(\SL_{2}^3,\SL_{2}^3,\hat{\rho},1)$ and $(\GL_{2}^3,S(\GL_{2}^3),\hat{\rho},1)$ where $\hat{\rho}$ is the tensor product of $\SL_{2}^3$ and $S(\GL_{2}^3)$ respectively, and they are equal to each other up to some central isogeny. While there are various choices of dual quadruples pairs $(\Delta,\hat\Delta)$ associated to $\hat\rho$ due to the isogeny issue, in this paper, for each representation $\hat{\rho}$ in \cite{K}, we will only write down one quadruple $\Delta=(G,H,\rho_H,\iota)$ whose dual quadruple $\hat\Delta$ is $(\hat{G},\widehat{G/Z_\Delta},\hat{\rho},1)$ where $Z_{\Delta}=Z_G\cap ker(\rho_H)$.

\begin{rmk}
	In our proof of Theorem \ref{main thm 1}, we frequently quote the unramified computation in \cite{II} and \cite{WZ}. The settings in \cite{II} and \cite{WZ} may actually differ from ours through finite isogeny or central isogeny. It is clear that the computation can be adapted and the results there still apply. For example, in \cite{II}, they computed the local relative character for the Gross-Prasad model $(\SO_{n+1}\times \SO_n,\SO_n)$ at unramified places. Their results can be also applied to models like $(\GL_4\times \GSp_4,\GSp_4)$ (which is essentially the Gross-Prasad model $(\SO_6\times \SO_5,\SO_5)$ up to some central isogeny). 
\end{rmk}

\subsection{Theta correspondence for classical groups}\label{sec theta}
In this paper we will frequently use theta correspondence for classical groups. We will briefly review it in this subsection. We start with the theta correspondence for the general linear group. Let $n\geq m\geq 1$ and $G=H_1\times H_2=\GL_n\times \GL_m$. We use $V$ to denote the underlying vector space of the representation $\rho=std_{\GL_n}\otimes std_{\GL_m}$ of $G$. For $\varphi\in \CS(V(\BA))$, we define the theta function
$$\Theta_{\psi}^{\varphi}(g)=\sum_{X\in V(k)}\rho(g)\varphi(X),\;g\in G(\BA)$$
which is an automorphic function on $G(\BA)=H_1(\BA)\times H_2(\BA)$. Let $\pi$ be a cuspidal automorphic representation of $H_2(\BA)$. For $\phi\in L^{2}(H_2(k)\back H_2(\BA))_{\pi}$, the integral 
$$\int_{H_2(k)\back H_2(\BA)} \Theta_{\psi}^{\varphi}(h_1,h_2)\phi(h_2)dh_2$$
gives an automorphic function on $H_1(\BA)$ which will be denoted by $\Theta(\phi)$.

\begin{thm}\label{theta GL}(\cite{L})
We have
	$$\{\Theta(\phi)|\;\phi\in L^{2}(H_2(k)\back H_2(\BA))_{\pi}\}=\{E(\phi',1)|\; \phi'\in L^{2}(H_2(k)\back H_2(\BA))_{\pi}\}$$ 
	where $E(\phi',1)$ is the Eisenstein series on $H_1(\BA)=\GL_n(\BA)$ induced from $\phi'$ and the identity function on $\GL_{n-m}(\BA)$. Moreover, for $\phi_1,\phi_2\in L^{2}(H_2(k)\back H_2(\BA))_{\pi}$, we have the Rallis inner product formula
$$\int_{H_2(k)\back H_2(\BA)/Z_{H_2}(\BA)} \int_{H_1(k)\back H_1(\BA)} \int_{H_1(k)\back H_1(\BA)}  \Theta_{\psi}^{\varphi}(h_1,h_2)  \Theta_{\psi}^{\varphi}(h_1',h_2) E(\phi_1,1)(h_1) E(\phi_2,1)(h_1')dh_1dh_1' dh_2  $$
$$``=" {\rm Res}_{s=\frac{n-m}{2}}L(s+\frac{1}{2},\pi)\cdot \int_{H_2(k)\back H_2(\BA)/Z_{H_2}(\BA)} \phi_1(h_2)\phi_2(h_2)dh_2.$$
\end{thm}

\begin{rmk}
When $m=1$, the above theorem implies that if we integrate the theta series on $\GL_n$ associated to the symplectic representation $T(std_{n})$ over the center of $\GL_n$ we will get the mirabolic Eisenstein series of $\GL_n$. We will frequently use this fact in later discussions.
\end{rmk}

For the unramified computation, we also need the local theta correspondence for unramified representation. Let $F$ be a p-adic local field that is a local place of $k$. We use $\phi_{\rho}(h_1,h_2)$ to denote the local spherical matrix coefficient of the Weil representation with $\phi_{\rho}(1,1)=1$. Let $\pi$ be a tempered unramified representation of $H_2(F)$, $\phi_\pi$ (resp. $\phi_{\pi,1}$) be the unramified matrix coefficient of $\pi$ (resp. $Ind_{\GL_m\times \GL_{n-m}}^{\GL_n}(\pi\otimes 1)$) with $\phi_\pi(1)=\phi_{\pi,1}(1)=1$.  

\begin{thm}\label{theta GL unramified}(\cite{L})
With the notation above, we have
$$\int_{H_2(F)} \phi_{\rho}(h_1,h_2) \phi_\pi(h_2) dh_2=L(\frac{n-m+1}{2},\pi)\cdot \phi_{\pi,1}(h_1).$$
\end{thm}

Next we study the theta correspondence between $\SO_{2n}$ and $\Sp_{2m}$ with $n\geq m\geq 1$. Let $G=H_1\times H_2=\SO_{2n}\times \Sp_{2m}$ and we use $V$ to denote the underlying vector space of the representation $\rho=std_{\SO_{2n}}\otimes std_{\Sp_{2m}}$ of $G$. Let $Y$ be a maximal isotropic subspace of $V$, we can define $\Theta_{\psi}^{\varphi}(g)$ an automorphic function on $G(\BA)$ as in the introduction, for any Schwartz function $\varphi$ on $Y$. 

Let $\Pi$ be a cuspidal tempered global Arthur packet of $H_2(\BA)=\Sp_{2m}(\BA)$ and let $\Pi'$ be its lifting to $H_1(\BA)=\SO_{2n}(\BA)$ under the map $\SO_{2m+1}(\BC)\times \SL_2(\BC)\rightarrow \SO_{2n}(\BA)$ whose restrict to $\SL_2$ is the principal embedding from $\SL_2$ to $\SO_{2n-2m-1}$ (if $n>m$ then $\Pi'$ is a non-tempered Arthur L-packet) \footnote{in fact here $\Pi'$ should be an Arthur packet of $O_{2n}(\BA)$ which is the union of two Arthur packets of $\SO_{2n}(\BA)$ differed by the outer automorphism}. For $\phi\in L^{2}(H_2(k)\back H_2(\BA))_{\pi}$, the integral
$$\int_{H_2(k)\back H_2(\BA)} \Theta_{\psi}^{\varphi}(h_1,h_2)\phi(h_2) dh_2$$
gives an automorphic function on $H_1(\BA)=\SO_{2n}(\BA)$ which will be denoted by $\Theta(\phi)$. Then the following theorem holds.

\begin{thm}\label{theta}(\cite{KR, Y, GQT})
	With the notation above, the representation
	$$\{\Theta(\phi)|\;\phi\in L^{2}(\Sp_{2m}(k)\back \Sp_{2m}(\BA))_{\Pi}\}$$ 
	of $\SO_{2n}(\BA)$ is a direct sum of some distinct irreducible representations belonging to the Arthur L-packet $\Pi'$ of $H_1(\BA)=\SO_{2n}(\BA)$. Moreover, for $\phi_1,\phi_2\in \Pi'$, we have the Rallis inner product formula
$$\int_{H_2(k)\back H_2(\BA)} \int_{H_1(k)\back H_1(\BA)} \int_{H_1(k)\back H_1(\BA)}  \Theta_{\psi}^{\varphi}(h_1,h_2)  \Theta_{\psi}^{\varphi}(h_1',h_2) \phi_1(h_1) \phi_2(h_1')dh_1dh_1' dh_2  $$
$$``="{\rm Res}_{s=\frac{2n-2m-1}{2}}L(s+\frac{1}{2},\Pi')\cdot \int_{H_1(k)\back H_1(\BA)} \phi_1(h_1)\phi_2(h_1)dh_1.$$
\end{thm}

For the unramified computation, we also need the local theta correspondence for unramified representation. Let $F$ be a p-adic local field that is a local place of $k$. We use $\phi_{\rho}(h_1,h_2)$ to denote the local spherical matrix coefficient of the Weil representation with $\phi_{\rho}(1,1)=1$. Let $\pi$ be a tempered unramified representation of $H_2(F)$ and $\pi'$ be its lifting to $H_1(F)$ (which is also unramified). Let $\phi_\pi$ (resp. $\phi_{\pi'}$) be the unramified matrix coefficient of $\pi$ (resp. $\pi'$) with $\phi_\pi(1)=\phi_{\pi'}(1)=1$. 

\begin{thm}\label{theta unramified}(\cite{L})
With the notation above, we have
$$\int_{H_2(F)} \phi_{\rho}(h_1,h_2) \phi_\pi(h_2) dh_2=L(n-m,\pi')\cdot \phi_{\pi'}(h_1).$$
\end{thm}

The theta correspondence between $\SO_{2m}$ and $ \Sp_{2n}$ (resp. $\GSO_{2n}$ and $\GSp_{2m}$, $\GSO_{2m}$ and $\GSp_{2n}$) is similar and we will skip it here.

\subsection{A conjecture of the duality under certain induction}
We recall the notion from the introduction. Let $\Delta=(G,H,\rho_H,\iota)$ be a BZSV quadruple. Let $M$ be the centralizer of $\{\iota(diag(t,t^{-1}))|\;t\in \GL_1\}$ in $G$. It is easy to see that $M$ is a Levi of $G$ and $H\subset M$. We define
$$\Delta_{red}=(M,H,1,\rho_{H,\iota})$$
where the representation $\rho_{H,\iota}$ has been defined in \ref{rho H iota}. It is clear that $\Delta$ is reductive if and only if $\Delta=\Delta_{red}$. 

In Section 4.2.2 of \cite{BSV}, Ben-Zvi--Sakellaridis--Venkatesh made a conjecture about the relation between the dual of $\Delta$ and $\Delta_{red}$. To state their conjecture, we first need a definition.

\begin{defn}
Let $M$ be a Levi subgroup of $G$ and $\rho$ be an irreducible representation of $M$ with the highest weight $\varpi_M$. There exists a Weyl element $w$ of $G$ such that $w\varpi_M$ is a dominant weight of $G$ \footnote{the choice of $w$ is not unique but $w\varpi_M$ is uniquely determined by $\varpi_M$}. We define $(\rho)_{M}^{G}$ to be the irreducible representation of $G$ whose highest weight is $w\varpi_M$. In general, if $\rho=\oplus_i\rho_i$ is a finite-dimensional representation of $M$ with $\rho_i$ irreducible, we define
$$(\rho)_{M}^{G}=\oplus_i(\rho_i)_{M}^{G}.$$
\end{defn}

Now we are ready to state the conjecture.

\begin{conj}\label{Whittaker induction}
With the notation above. If the dual of $\Delta_{red}$ is given by $\hat{\Delta}_{red}=(\hat{M}, \hat{H}_M', \rho',\hat{\iota}')$, then the dual of $\Delta$ is given by
$$(\hat{G},\hat{H}', (\rho')_{\hat{H}_M'}^{\hat{H}'},\hat{\iota}' )$$
where $\hat{H}'$ is generated by $\hat{H}_M'$ and $\{Im(\iota_\alpha)|\;\alpha\in \Delta_{\hat{G}}-\Delta_{\hat{M}}\}$. Here $\Delta_{\hat{G}}$ (resp. $\Delta_{\hat{M}}$) is the set of simple roots of $\hat{G}$ (resp. $\hat{M}$) and $\iota_\alpha:\SL_2\rightarrow \hat{G}$ is the embedding associated to $\alpha$.
\end{conj}

\subsection{General strategy}
Let $\hat{\Delta}=(\hat{G},\hat{G},\hat{\rho},1)$ be a quadruple such that $\hat{\rho}$ is an anomaly-free symplectic representation of $\hat{G}$, and it appears in Table 1, 2, 11, 12, 22 of \cite{K}. Our goal is to write down a dual quadruple (up to isogeny) $\Delta=(G,H,\rho_H,\iota)$.

The data in Knop's tables of \cite{K}, besides $(\hat{G},\hat{\rho})$, also contains the following two items: a Levi subgroup $\hat{L}$ of $\hat{G}$ and a Weyl group $\hat{W}_V$ written in the form of $W_{\hat{H}}$ where $\hat{H}$ is the root type (e.g. $A_n,B_n,C_n$, etc). (In \cite{K} the notations are $L,G,W_V$ in place of $\hat{L},\hat{G},\hat{W}_V$  respectively.) Our key observation is that two data $(H,\iota)$ of the dual quadruple $\Delta=(G,H,\rho_H,\iota)$ are given by the following properties.

\begin{prpt}\label{main property}
	\begin{enumerate}
		\item The root type of $H$ is dual to the root type of $\hat{W}_V$ in the tables of \cite{K}.
		\item The nilpotent orbit ${\mathcal O}_\iota$ associated to $\iota$ is the principal nilpotent orbit of $L$ where $L$ is the dual Levi of $\hat{L}$.
	\end{enumerate}
\end{prpt}

\begin{rmk}
	Basically, the Weyl group $\hat{W}_V$ can be viewed as the ``little Weyl group" of the quadruple $\hat{\Delta}=(\hat{G},\hat{G},\hat{\rho},1)$, and $\hat{\Fl}$ in tables of \cite{K}  is an analogue of $\hat{\Fl}_X$ in Table 3 of \cite{Knop}.
\end{rmk}

As a result, it remains to find out what is $\rho_H$. We do not have a systematic way to write down $\rho_H$. Instead we propose a $\rho_H$ in an ad hoc way and then provide evidence for the duality between $\Delta=(G,H,\rho_H,\iota)$ and  $(\hat{G},\widehat{G/Z_\Delta},\hat{\rho},1)$. 

We provide two strong evidences for the duality. The first one is evidence for Conjecture \ref{BSV conj}, i.e., Theorem \ref{main thm 1} and \ref{main thm 2}. The second evidence is for non-reductive models. For those models, we will show that the duality is compatible with Conjecture \ref{Whittaker induction}.

In the sections that follow, we will go through Knop's list of representations $\hat \rho$. For each $\rho$ we write down a quadruple $(G,H,\rho_H,\iota)$. When the quadruple is not reductive, we will also write down $\Delta_{red}$ which is dual to another representation $(\hat M,\hat\rho_M)$ in Knop's list and verify that Theorem \ref{main thm 3} holds.
For cases in Table \ref{red list}, \ref{red list extra}, \ref{non-red list 1} and \ref{non-red list 1 extra}, we give references where the local relative character is calculated in the unramified places, thus verifying Theorem \ref{main thm 1}. We also verify Theorem~\ref{main thm 2} for the global periods associated to the dual side $\hat\Delta$  for cases in Table \ref{red list}, \ref{non-red list 1} and \ref{non-red list 2}.



\section{Models in Table 1 of \cite{K}}\label{sec Table 1}
In this section we will consider Table 1 of \cite{K}, this is for the case when $\hat{\rho}$ is an irreducible representation of $\hat{G}$. It is easy to check that the representations in (1.2), (1.8), (1.9) and (1.10) of \cite{K} are not anomaly free and the representation in (1.1) of \cite{K} is only anomaly free when $p=2n$ is even. Hence it remains to consider the following cases. Note that we only write the root type of $\hat{\Fl}$ and we write 0 if it is abelian. Also we separate the cases when $\hat{\Fl}$ is abelian and when $\hat{\Fl}$ is not abelian. These are precisely the cases where the dual quadruple is reductive/non-reductive (see Property \ref{main property}).

\begin{figure}[h!]
\begin{tabular}{| c | c | c |c|}
\hline
Number in \cite{K} & $(\hat{G},\hat{\rho})$ & $\hat{W}_V$ &  $\hat{\Fl}$ \\
\hline
(1.1), p=2m & $(\Sp_{2m}\times \SO_{2m},std_{\Sp_{2m}}\otimes std_{\SO_{2m}})$ & $D_m$ & $0$ \\
\hline

(1.1), p=2m+2 & $(\Sp_{2m}\times \SO_{2m+2},std_{\Sp_{2m}}\otimes std_{\SO_{2m+2}})$ & $C_m$ & $0$ \\
\hline

(1.3), m=2  & $(\Spin_5\otimes \Spin_7,\Spin_5\otimes \Spin_7)$ & $C_2\times A_1$ & 0\\

\hline

(1.3), m=3  & $(\Sp_6\otimes \Spin_7,std_{\Sp_6}\otimes \Spin_7)$ & $C_3\times B_3$ & 0\\

\hline

(1.3), m=4  & $(\Sp_8\otimes \Spin_7,std_{\Sp_8}\otimes \Spin_7)$ & $D_4\times B_3$ & 0\\

\hline

(1.6) & $(SL_2,Sym^3)$ & $A_1$ &  $0$\\
\hline

\hline
\end{tabular}
\captionof{table}{Reductive models in Table 1 of \cite{K}}
\label{table 1 reductive}
\end{figure}

\begin{figure}[h!]
\begin{tabular}{| c | c | c |c|}
\hline
Number in \cite{K} & $(\hat{G},\hat{\rho})$ & $\hat{W}_V$ &  $\hat{\Fl}$ \\
\hline

(1.1), $p=2n<2m$ & $(\Sp_{2m}\times \SO_{2n},std_{\Sp_{2m}}\otimes std_{\SO_{2n}})$ & $D_n$ & $C_{m-n}$ \\
\hline

(1.1), $p=2n>2m+2$ & $(\Sp_{2m}\times \SO_{2n},std_{\Sp_{2m}}\otimes std_{\SO_{2n}})$ & $C_m$ & $D_{n-m}$ \\
\hline

(1.3), m=1  & $(\SL_2\times \Spin_{7},std_{\SL_2}\otimes \Spin_7)$ & $A_1$ & $A_2$\\

\hline

(1.3), $m>4$  & $(\Sp_{2m}\otimes \Spin_7,std_{\Sp_{2m}}\otimes \Spin_7)$ & $D_4\times B_3$ & $C_{m-4}$\\

\hline

(1.4) & $(\SL_2\times \Spin_{9},std_{\SL_2}\otimes \Spin_9)$ & $A_1\times A_1$ &  $A_2$\\
\hline

(1.5), n=11 & $(\Spin_{11},\Spin_{11})$ & $A_1$ &  $A_4$\\
\hline

(1.5), n=12 & $(\Spin_{12},\HSpin_{12})$ & $A_1$ &  $A_5$\\
\hline

(1.5), n=13 & $(\Spin_{13},\Spin_{13})$ & $B_2$ &  $A_2\times A_2$\\
\hline

(1.7) & $(SL_6,\wedge^3)$ & $A_1$ &  $A_2\times A_2$\\
\hline

(1.11) & $(E_{7},std_{E_7})$ & $A_1$ &  $E_6$\\
\hline
\end{tabular}
\captionof{table}{Non-reductive models in Table 1 of \cite{K}}
\label{Table 1 non-red}
\end{figure}

\subsection{The reductive case}\label{sec Table 1 red}
In this subsection we consider the reductive cases, i.e., the ones in Table \ref{table 1 reductive}. The nilpotent orbit $\iota$ is trivial for all these cases so we will ignore it. 

For (1.1) with $p=2m$ (resp. $p=2m+2$), the associated quadruple $\Delta$ is 
\begin{equation}\label{1.1 p=2m}
(G,H,\rho_H)=(\SO_{2m+1}\times \SO_{2m},\SO_{2m}, 0)
\end{equation}

\begin{equation}\label{1.1 p=2m+2}
(\text{resp.} (G,H,\rho_H)=(\SO_{2m+1}\times \SO_{2m+2},\SO_{2m+1}, 0))
\end{equation}
which is just the reductive Gross-Prasad model. The unramified computations in \cite{II} prove Theorem \ref{main thm 1} in these two cases. For the dual side, Theorem~\ref{theta} applied to the theta correspondence  between $\SO_{2m}\times \Sp_{2m}$ (resp. $\SO_{2m+2}\times \Sp_{2m}$) implies Conjecture \ref{BSV conj}(2) and this proves Theorem \ref{main thm 2}.

For (1.3) with $m=2$, the associated quadruple $\Delta$ is $$(G,H,\rho_H)=(\GSp_6\times \GSp_4,G(\Sp_4\times \Sp_2),0)$$ 
which is the model $(\GSp_6\times \GSp_4,G(\Sp_4\times \Sp_2))$ studied in \cite{WZ}. The unramified computations in \cite{WZ} prove Theorem \ref{main thm 1} in this case.

For (1.3) with $m=3$, the associated quadruple $\Delta$ is 
$$(G,H,\rho_H)=(\GSp_6\times \GSpin_{7},S(\GSp_6\times \GSpin_{7}),std_{\Sp_6}\otimes \Spin_7).$$
For (1.3) with $m=4$, the associated quadruple $\Delta$ is
\begin{equation}\label{1.3 m=4}
(G,H,\rho_H)=(\GSp_6\times \GSpin_{9},S(\GSp_6\times \GSpin_{8}),std_{\Sp_6}\otimes \HSpin_8).
\end{equation}
Theorem~\ref{main thm 1} and \ref{main thm 2} for two cases can be established by the same argument as Model (11.11) of \cite{K} (see \eqref{11.11 p=2m-1} and \eqref{11.11 p=2m+1} of Section \ref{sec Table 11 red}) together with the triality of $D_4$. 

For (1.6), it is clear that the generic stabilizer of $\hat{\rho}$ in $\hat{G}$ is not connected, hence it does not belong to the current framework of BZSV duality. However, for this specific case, by the work of \cite{GJR01}, we expect there is an associated quadruple of the form $(\GL_2,\GL_2,\rho_H,1)$  where $\rho_H$ is no longer an anomaly free symplectic representation, but rather we understand that $\rho_H$ corresponds to the theta series on $H=\GL_2$ defined via the cubic covering of $\GL_2$ as in \cite{GJR01}. There is a covering group involved in the integral since the generic stabilizer is not connected. In \cite{GJR01} it is established that the nonvanishing of $\CP_{H,\iota,\rho_H}(\phi)$ is equivalent to the nonvanishing of $L(1/2,\Pi,\hat\rho)$. We expect further that Conjecture~\ref{BSV conj}(1) holds in this case.
	

By the discussion above, the strongly tempered quadruple associated to Table \ref{table 1 reductive} (without the row corresponding to (1.6)) is given as follows. Note that $\iota$ is trivial for all these cases.

\begin{figure}[h!]
\begin{tabular}{| c   |c| c|}
\hline
 (G, H, $\rho_H$)&  $\hat{\rho}$ \\
\hline
 ($\SO_{2m+1}\times \SO_{2m}$, $\SO_{2m}$, 0)& $std_{\Sp_{2m}}\otimes std_{\SO_{2m}}$ \\
\hline

 ($\SO_{2m+2}\times \SO_{2m+1}$, $\SO_{2m+1}$, 0)& $std_{\Sp_{2m}}\otimes std_{\SO_{2m+2}}$ \\
\hline

($\GSp_6\times \GSp_4$, $G(\Sp_4\times \Sp_2)$,0) & $\Spin_5\otimes \Spin_7$\\

\hline

 ($\GSp_6\times \GSpin_{7}$, $S(\GSp_6\times \GSpin_{7})$, $std_{\Sp_6}\otimes \Spin_7$) & $std_{\Sp_6}\otimes \Spin_7$ \\

\hline

 $(\GSp_6\times \GSpin_{9},S(\GSp_6\times \GSpin_{8}),std_{\Sp_6}\otimes \HSpin_8)$  & $std_{\Sp_8}\otimes \Spin_7$\\

\hline

\end{tabular}
\captionof{table}{Dual quadruples of Table \ref{table 1 reductive}}
\label{table 1 red list}
\end{figure}

\subsection{The non-reductive case}\label{sec Table 1 non-red}
In this subsection we consider the non-reductive cases, i.e., the ones in Table \ref{Table 1 non-red}. 

For (1.1) with $p=2n<2m$, the associated quadruple $\Delta$ is 
$$(\SO_{2m+1}\times \SO_{2n},\SO_{2n},0,(\GL_1)^{n}\times \SO_{2m-2n+1}\times T_{\SO_{2n}})$$
and it is the Gross-Prasad period for $\SO_{2m+1}\times \SO_{2n}$. For (1.1) with $p=2n>2m+2$, the associated quadruple $\Delta$ is 
$$(\SO_{2m+1}\times \SO_{2n},\SO_{2m+1},0,T_{\SO_{2m+1}}\times  (\GL_1)^{m}\times \SO_{2n-2m})$$
and it is still the Gross-Prasad period for $\SO_{2m+1}\times \SO_{2n}$. In these two cases $\Delta_{red}$ are given by \eqref{1.1 p=2m}, \eqref{1.1 p=2m+2}. It is clear that Theorem \ref{main thm 3} holds  in these two cases. The unramified computation in \cite{II} proves Theorem \ref{main thm 1} for these two cases. Theorem~\ref{theta} applied to the theta correspondence between  $\SO_{2n}\times \Sp_{2m}$ implies Conjecture \ref{BSV conj}(2) and proves Theorem \ref{main thm 2} for these two cases.

For (1.3) when $m=1$, the associated quadruple $\Delta$ is

\begin{equation}\label{GSp(6)xGL(2)}
(\GSp_6\times \GL_2,\GL_2,0,(\GL_3\times \GL_1)\times T_{\GL_2})
\end{equation}
and it is the model $(\GSp_6\times \GL_2,\GL_2\ltimes U)$ studied in \cite{WZ}. In this case $\Delta_{red}=((\GL_2)^3,\GL_2,0,1)$ (which a special case of \eqref{1.1 p=2m+2} with $m=1$). It is clear that Theorem \ref{main thm 3} holds in this case and the unramified computation in \cite{WZ} proves Theorem \ref{main thm 1} in this case.

For (1.3) when $m>4$, the associated quadruple $\Delta$ is
$$(\GSpin_{2m+1}\times \GSp_6,S(\GSpin_8\times \GSp_6),std_{\Sp_6}\otimes \HSpin_8, L)$$
where $L$ is the Levi subgroup whose projection to $\GSpin_{2m+1}$ (resp. $\GSp_6$) is of the form $(\GL_1)^4\times \GSpin_{2m-7}$ (resp. the maximal torus). The nilpotent orbit induces a Bessel period for the unipotent radical of the parabolic subgroup $P=MU$ with $M=(\GL_1)^{m-4}\times \GSpin_9\times \GSp_6$ whose stabilizer is $\GSpin_8\times \GSp_6$ and we can naturally embed $H$ into the stabilizer. In this case $\Delta_{red}$ is given by \eqref{1.3 m=4} and it is clear that Theorem \ref{main thm 3} holds. Theorem~\ref{main thm 1} and \ref{main thm 2} for this model can be established by the same argument as \eqref{11.11 non-red} in Section \ref{sec Table 11 non-red} together with the triality of $D_4$.

For (1.4), the associated quadruple $\Delta$ is 
\begin{equation}\label{GSp(8)xGL(2)}
(\GSp_8\times \GL_2,G(\SL_2\times \SL_2),0,\GL_3\times \GL_1\times \GL_1\times T_{\GL_2}).
\end{equation}
The nilpotent orbit induces a Bessel period for the unipotent radical of the parabolic subgroup $P=MU$ with $M=\GL_2\times \GSp_4\times \GL_2$ whose stabilizer is $G(\SL_2\times \SL_2)\times \GL_2$. We embeds $H$ into the stabilizer so that the induced embedding from $H$ into $M$ is given by the natural embeddings of $H$ into $\GSp_4$ and into $\GL_2\times \GL_2$. In this case $\Delta_{red}=(\GSp_4\times \GL_2\times \GL_2,G(\SL_2\times \SL_2),0,1)$ which is essentially the Gross-Prasad model for $\SO_5\times \SO_4$. If we replace the cusp form on $\GL_2$ by an Eisenstein series, we recover the Rankin-Selberg integrals in \cite{BG00}. It is clear that Theorem \ref{main thm 3} holds in this case and the unramfied computation in \cite{BG00} proves Theorem \ref{main thm 1} in this case.

For (1.5) when $n=11$, the associated quadruple $\Delta$ is
\begin{equation}\label{GSp(10) case}
(\GSp_{10},\GL_2,0,\GL_5\times \GL_1)
\end{equation}
and it is the model $(\GSp_{10},\GL_2\ltimes U)$ studied in \cite{WZ}. In this case $\Delta_{red}=((\GL_2)^3,\GL_2,0,1)$. It is clear that Theorem \ref{main thm 3} holds in this case and the unramified computation in \cite{WZ} proves Theorem \ref{main thm 1} in this case.

For (1.5) when $n=12$, the associated quadruple $\Delta$ is
$$(\GSO_{12},\GL_2,0,\GL_6\times \GL_1)$$
and it is the model $(\GSO_{12},\GL_2\ltimes U)$ studied in \cite{WZ}. In this case $\Delta_{red}=((\GL_2)^3,\GL_2,0,1)$.  It is clear that Theorem \ref{main thm 3} holds in this case and the unramified computation in \cite{WZ} proves Theorem \ref{main thm 1} in this case.

For (1.5) when $n=13$, the associated quadruple $\Delta$ is
$$(\GSp_{12},\GSp_4,0,\GL_3\times \GL_3\times \GL_1).$$
The nilpotent orbit induces a Bessel period for the unipotent radical of the parabolic subgroup $P=MU$ with $M=\GL_4\times \GSp_4$ whose stabilizer is $H=\GSp_4$. In this case $\Delta_{red}=(\GSp_4\times \GL_4,\GSp_4,0,1)$ which is essentially the Gross-Prasad model for $\SO_6\times \SO_5$. It is clear that Theorem \ref{main thm 3} holds in this case. In this case the unramified computation can be done in a similar way as \cite{WZ}, which will give Theorem~\ref{main thm 1}.

For (1.7), the associated quadruple $\Delta$ is
$$(\GL_6,\GL_2,0,\GL_3\times \GL_3)$$
and it is the Ginzburg-Rallis model $(\GL_6,\GL_2\ltimes U)$ studied in \cite{WZ}. In this case $\Delta_{red}=((\GL_2)^3,\GL_2,0,1)$. It is clear that Theorem \ref{main thm 3} holds in this case and the unramified computation in \cite{WZ} proves Theorem \ref{main thm 1} in this case.

For (1.11), the associated quadruple $\Delta$ is
$$(E_{7},\PGL_2,0,GE_6)$$
and it is the model $(E_7,\PGL_2\ltimes U)$ studied in \cite{WZ}. In this case $\Delta_{red}=((\PGL_2)^3,\PGL_2,0,1)$. It is clear that Theorem \ref{main thm 3} holds in this case and the unramified computation in \cite{WZ} proves Theorem \ref{main thm 1} in this case.

By the discussion above, the strongly tempered quadruple associated to Table \ref{Table 1 non-red} is given as follows. Here for $\iota$, we only list the root type of the Levi subgroup $L$ of $G$ such that $\iota$ is principal in $L$. 

\begin{figure}[h!]
\begin{tabular}{| c | c | c |}
\hline
$(G,H,\rho_H)$ & $\iota$ &$\hat{\rho}$\\
\hline

 $(\SO_{2m+1}\times \SO_{2n},\SO_{2n},0)$ &$B_{m-n}$ & $std_{\Sp_{2m}}\otimes std_{\SO_{2n}}$  \\
\hline

 $(\SO_{2m+1}\times \SO_{2n},\SO_{2m+1},0)$ &$D_{n-m}$ & $std_{\Sp_{2m}}\otimes std_{\SO_{2n}}$ \\
\hline

$(\GSp_6\times \GL_2,\GL_2,0)$  & $A_2$& $std_{\GL_2}\otimes \Spin_7$ \\

\hline

 $(\GSpin_{2m+1}\times \GSp_6,S(\GSpin_8\times \GSp_6),std_{\Sp_6}\otimes \HSpin_8)$ &$B_{m-4}$ & $std_{\Sp_{2m}}\otimes \Spin_7$  \\

\hline

$(\GSp_8\times \GL_2,G(\SL_2\times \SL_2),0)$ &$A_2$ & $std_{\GL_2}\otimes \Spin_9$  \\
\hline

$(\GSp_{10},\GL_2,0)$ &$A_4$& $\Spin_{11}$  \\
\hline

$(\GSO_{12},\GL_2,0)$ &$A_5$ & $\HSpin_{12}$ \\
\hline

$(\GSp_{12},\GSp_4,0)$ &$A_2\times A_2$ & $\Spin_{13}$   \\
\hline

$(\GL_6,\GL_2,0)$ & $A_2\times A_2$ & $\wedge^3$  \\
\hline

$(E_{7},\PGL_2,0)$ & $E_6$ & $std_{E_7}$   \\
\hline
\end{tabular}
\captionof{table}{Dual quadruples of Table \ref{Table 1 non-red}}
\label{Table 1 non-red list}
\end{figure}

\section{Models in Table 2}\label{sec Table 2}

In this section we will consider Table 2 of \cite{K}, this is for the case when $\hat{\rho}=T(\hat{\tau})$ is the direct sum of two irreducible representations of $\hat{G}$ that are dual to each other. All the representations in Table 2 of \cite{K} are anomaly free, so we need to consider all of them. We still separate the cases based on whether
$\hat{\Fl}$ is abelian or not.

\begin{figure}[h!]
\begin{tabular}{| c | c | c |c|}
\hline
Number in \cite{K} & $(\hat{G},\hat{\rho})$ & $\hat{W}_V$ &  $\hat{\Fl}$ \\
\hline

(2.1), m=n & $(\GL_n\times \GL_n,T(std_{\GL_n}\otimes std_{\GL_n}))$ & $A_{n-1}$ &0 \\
\hline

(2.1), m=n+1 and (2.4), n=2 & $(\GL_{n+1}\times \GL_n,T(std_{\GL_{n+1}}\otimes std_{\GL_n}))$ & $A_{n-1}$ &0 \\
\hline

(2.3) & $(\GL_n,T(Sym^2))$ & $A_{n-1}$ &0 \\

\hline

(2.6), m=n=2 & $(\Sp_4\times \GL_2,T(Std_{\Sp_4}\otimes Std_{\GL_2}))$ & $A_{1}\times A_1$ &0 \\

\hline

(2.6), m=2, n=3 & $(\Sp_4\times \GL_3,T(Std_{\Sp_4}\otimes Std_{\GL_3}))$ & $C_2\times A_2$ &0 \\

\hline

(2.6), m=2, n=4 & $(\Sp_4\times \GL_4,T(Std_{\Sp_4}\otimes Std_{\GL_4}))$ & $C_2\times A_3$ &0 \\

\hline

(2.6), m=2, n=5 & $(\Sp_4\times \GL_5,T(Std_{\Sp_4}\otimes Std_{\SL_5}))$ & $C_2\times A_3$ &0 \\

\hline

(2.6), m=n=3 & $(\Sp_6\times \GL_3),T(Std_{\Sp_6}\otimes Std_{\GL_3}))$ & $A_{3}\times A_2$ &0 \\

\hline

\end{tabular}
\captionof{table}{Reductive models in Table 2 of \cite{K}}
\label{Table 2 red}
\end{figure}

\begin{figure}[h!]
\begin{tabular}{| c | c | c |c|}
\hline
Number in \cite{K} & $(\hat{G},\hat{\rho})$ & $\hat{W}_V$ &  $\hat{\Fl}$ \\
\hline

(2.1), $m>n+1$, and (2.4), $n>2$ & $(\GL_m\times \GL_n,T(std_{\GL_m}\otimes std_{\GL_n}))$ & $A_{n-1}$ & $A_{m-n-1}$ \\

\hline

(2.2), n=2m & $(\GL_{2m},T(\wedge^2))$ & $A_{m-1}$ & $(A_1)^m$ \\
\hline

(2.2),n=2m+1 & $(\GL_{2m+1},T(\wedge^2))$ & $A_{m-1}$ & $(A_1)^m$ \\
\hline

(2.5) & $(\Sp_{2n},T(std_{\Sp_{2n}})$ & $0$ & $C_{m-1}$ \\
\hline

(2.6), $m>2$, n=2 & $(\Sp_{2m}\times \SL_2,T(Std_{\Sp_{2m}}\otimes Std_{\SL_2}))$ & $A_{1}\times A_1$ &$C_{m-2}$ \\

\hline

(2.6), m=2, $n>5$ & $(\Sp_4\times \SL_m,T(Std_{\Sp_4}\otimes Std_{\SL_m}))$ & $C_2\times A_3$ &$A_{m-5}$ \\

\hline

(2.6), $m>3$, n=3 & $(\Sp_{2m}\times \SL_3,T(Std_{\Sp_{2m}}\otimes Std_{\SL_3}))$ & $A_{3}\times A_2$ &$C_{m-3}$ \\

\hline

(2.7), m=2k & $(\SO_{2k}, T(std_{\SO_{2k}}))$ & $A_1$ & $D_{k-1}$ \\

\hline

(2.7), m=2k+1 & $(\SO_{2k+1}, T(std_{\SO_{2k+1}}))$ & $A_1$ & $B_{k-1}$ \\

\hline

(2.8), n=7 & $(\Spin_7, T(\Spin_7))$ & $A_1$ & $A_2$ \\

\hline

(2.8), n=9 & $(\Spin_9, T(\Spin_9))$ & $A_1\times A_1$ & $A_2$ \\

\hline

(2.8), n=10 & $(\Spin_{10}, T(\HSpin_{10}))$ & $A_1$ & $A_3$ \\

\hline

(2.9) & $(G_2,T(std_{G_2}))$ & $A_1$ & $A_1$ \\
\hline

(2.10) & $(E_6,T(std_{E_6}))$ & $A_2$ & $D_4$ \\
\hline

\end{tabular}
\captionof{table}{Non-reductive models in Table 2 of \cite{K}}
\label{Table 2 non-red}
\end{figure}

\newpage

\subsection{The reductive case}\label{sec Table 2 red}
In this subsection we consider the reductive cases, i.e., the ones in Table \ref{Table 2 red}. 

For (2.1) with $m=n$, the associated quadruple $\Delta$ is given by 
\begin{equation}\label{2.1 m=n}
(G,H,\rho_H,\iota)=(\GL_n\times \GL_n,\GL_n,T(std_{\GL_n}),1).
\end{equation}
For (2.1) with $m=n+1$ and (2.4) with $n=2$, the associated quadruple $\Delta$ is given by 
\begin{equation}\label{2.1 m=n+1}
(G,H,\rho_H,\iota)=(\GL_{n+1}\times \GL_n,\GL_n,0,1).
\end{equation}
The period integrals in these two cases are exactly the Rankin-Selberg integral for $\GL_n\times \GL_n$ and $\GL_{n+1}\times \GL_n$ in \cite{JPSS}. The result in loc. cit. proves Conjecture \ref{BSV conj}(1) and Theorem \ref{main thm 1}. For the dual side, Theorem~\ref{theta GL} applied to the theta correspondence for $\GL_n\times \GL_{n+1}$ and $\GL_n\times \GL_n$  imply Conjecture \ref{BSV conj}(2) and this proves Theorem \ref{main thm 2}.

For (2.3),  the generic stabilizer of $\hat{\rho}$ in $\hat{G}$ is not connected,  hence it does not belong to the current framework of the BZSV duality. However, for this specific case, by the Rankin-Selberg integral in \cite{BG92-ann, PPS89,  Tk14}, we know that the dual integral should be the one in \cite{BG92-ann, PPS89,  Tk14}. As the generic stabilizer is not connected, there are covering groups involved in the integral.

For (2.6) with $m=n=2$, the associated quadruple $\Delta$ is given by 
\begin{equation}\label{2.6 m=n=2}
(G,H,\rho_H,\iota)=(\GSp_4\times \GL_2, G(\SL_2\times \SL_2), T(std_{\GL_2,2},),1)
\end{equation}
where the embedding of $H$ into $G$ is given by the canonical embedding from $\GSpin_4=G(\SL_2\times \SL_2)$ into $\GSpin_5=\GSp_4$ and the projection of $G(\SL_2\times \SL_2)$ into $\GL_2$ via the first $\GL_2$-copy. The representation $\rho_H$ is the standard representation of the second $\GL_2$-copy of $H$. This integral is essentially the Gross-Prasad model for $\SO_5\times \SO_4$ except we replace the cusp form on one $\GL_2$-copy by the theta series. The unramified computation in \cite{II} proves Theorem \ref{main thm 1} in this case. For the dual side, Conjecture \ref{BSV conj}(2) follows from Theorem~\ref{theta GL} applied to  the theta correspondence of $\GL_2\times \GL_4$ and Gan-Gross-Prasad conjecture (Conjecture 9.11 of \cite{GGP2}) for non-tempered Arthur packet for the pair $(\GL_4\times \GSp_4,\GSp_4)$ which is essentially the Gross-Prasad period for $\SO_6\times \SO_5$. This proves Theorem \ref{main thm 2}.

For (2.6) with $m=2, n=3$, the associated quadruple $\Delta$ is given by 
$$(G,H,\rho_H,\iota)=(\GSp_4\times \GL_3, \GSp_4\times \GL_3, T(std_{\GSp_4}\otimes std_{\GL_3}),1).$$
By the theta correspondence for $\GL_3\times \GL_4$ (note that the theta function constructed from $T(std_{\GSp_4}\otimes std_{\GL_3})$ is the  restriction of the theta function from $T(std_{\GL_4}\otimes std_{\GL_3})$), the integral over $\GL_3$ of a cusp form on $\GL_3$ with the theta series associated to $\rho_H$ produces an Eisenstein series of $\GL_4$ induced from the cusp form on $\GL_3$ and the trivial character of $\GL_1$. Then the integral over $\GSp_4$ is just the period integral for the pair $(\GL_4\times \GSp_4,\GSp_4)$ which is essentially the Gross-Prasad period for $\SO_6\times \SO_5$. The unramified computation in \cite{II} and Theorem~\ref{theta GL unramified} applied to theta correspondence for $\GL_3\times \GL_4$ proves Theorem \ref{main thm 1} in this case. For the dual side, Conjecture \ref{BSV conj}(2) follows from Theorem~\ref{theta GL} applied to the theta correspondence of $\GL_4\times \GL_3$ and the global period integral conjecture for the pair $(\GL_4\times \GSp_4,\GSp_4)$ (which is essentially the Gross-Prasad period for $\SO_6\times \SO_5$) in \cite{GGP}. This proves Theorem \ref{main thm 2}.

For (2.6) with $m=2, n=4$, the associated quadruple $\Delta$ is 
\begin{equation}\label{2.6, m=2 n=4}
(\GSp_4\times \GL_4, S(\GSp_4\times \GL_4),std_{\Sp_4}\otimes \wedge^2\oplus T(std_{\GL_4})).
\end{equation}
By the theta correspondence for $\GSp_4\times \GSO_6$, the integral over $\Sp_4$ of a cusp form on $\GSp_4$ with the theta series associated to $\rho_H$ produces an automorphic form of $\GL_4$. Then the integral over $\GL_4$ is just  the Rankin-Selberg integral of $\GL_4\times \GL_4$ as in \cite{JPSS} \footnote{In this paper we will frequently use the fact that the theta series associated to $\rho\oplus \rho'$ is the product of the theta series associated to $\rho$ and $\rho'$.}. The Rankin-Selberg integral in  \cite{JPSS} and Theorems~\ref{theta GL} and \ref{theta GL unramified} applied to theta correspondence for $\GSp_4\times \GSO_6$ proves Conjecture \ref{BSV conj}(1) and Theorem \ref{main thm 1} in this case. For the dual side, Conjecture \ref{BSV conj}(2) follows from Theorem~\ref{theta GL} applied to the theta correspondence of $\GL_4\times \GL_4$ and the global period integral conjecture for the pair $(\GL_4\times \GSp_4,\GSp_4)$ (which is essentially the Gross-Prasad period for $\SO_6\times \SO_5$) in \cite{GGP}. This proves Theorem \ref{main thm 2}. This is a very interesting case because both $\Delta$ and $\hat{\Delta}$ are strongly tempered and they are not equal to each other.

For (2.6) with $m=2, n=5$, the associated quadruple $\Delta$ is 
\begin{equation}\label{(2.6) m=2 n=5}
(\GSp_4\times \GL_5, S(\GSp_4\times \GL_4),std_{\Sp_4}\otimes \wedge^2).
\end{equation}
By the theta correspondence for $\GSp_4\times \GSO_6$, the integral over $\Sp_4$ of a cusp form on $\GSp_4$ with the theta series associated to $\rho_H$ produces an automorphic form of $\GL_4$. Then the integral over $\GL_4$ is just the Rankin-Selberg integral of $\GL_5\times \GL_4$.  The Rankin-Selberg integral in  \cite{JPSS} and Theorems~\ref{theta} and \ref{theta unramified} applied to theta correspondence $\GSp_4\times \GSO_6$ proves Conjecture \ref{BSV conj}(1) and Theorem \ref{main thm 1} in this case. For the dual side, Conjecture \ref{BSV conj}(2) follows from Theorem~\ref{theta GL} applied to the theta correspondence of $\GL_4\times \GL_5$ and the global period integral conjecture for the pair $(\GL_4\times \GSp_4,\GSp_4)$ (which is essentially the Gross-Prasad period for $\SO_6\times \SO_5$) in \cite{GGP}. This proves Theorem \ref{main thm 2}.

For (2.6) with $m=n=3$, the associated quadruple $\Delta$ is given by
\begin{equation}\label{2.6 m=n=3}
(\GSpin_7\times \GL_3,\GSpin_6\times \GL_3,T(\HSpin_6\otimes std_{\GL_3})).
\end{equation}
By the theta correspondence for $\GL_3\times \GL_4$ (note that $\GSpin_6$ is essentially $\GL_4$ up to some central isogeny which won't affect the unramified computation) the integral over $\GL_3$ of a cusp form on $\GL_3$ with the theta series associated to $\rho_H$ produces an Eisenstein series of $\GSpin_6$ induced from the cusp form on $\GL_3$ and the trivial character of $\GL_1$. Then the integral over $\GSpin_6$ is just the period integral for the Gross-Prasad model of $\GSpin_7\times \GSpin_6$. The unramified computation in \cite{II} and Theorem~\ref{theta GL unramified} applied to theta correspondence for $\GL_3\times \GL_4$ proves Theorem \ref{main thm 1} in this case. For the dual side, Conjecture \ref{BSV conj}(2) follows from Theorem~\ref{theta} applied to the theta correspondence of $\GSp_6\times \GSO_6$ and the Rankin-Selberg integral of $\GL_4\times \GL_3$. This proves Theorem \ref{main thm 2}.

By the discussion above, the strongly tempered quadruple associated to Table \ref{Table 2 red} is given as follows. Note that $\iota$ is trivial for all these cases.

\begin{figure}[h!]
\begin{tabular}{| c | c |}
\hline
$(G,H,\rho_H)$ & $\hat{\rho}$  \\
\hline

 $(\GL_n\times \GL_n,\GL_n,T(std_{\GL_n}))$ & $T(std_{\GL_n}\otimes std_{\GL_n})$  \\
\hline

 $(\GL_{n+1}\times \GL_n,\GL_n,0)$ & $T(std_{\GL_{n+1}}\otimes std_{\GL_n})$  \\
\hline

 $(\GSp_4\times \GL_2,G(\SL_2\times \SL_2),T(std_{\GL_2,2}))$ & $T(Std_{\GSp_4}\otimes Std_{\GL_2})$ \\

\hline

 $(\GSp_4\times \GL_3,H=G,T(std_{\GSp_4}\otimes std_{\GL_3}))$ & $T(Std_{\GSp_4}\otimes Std_{\GL_3})$\\

\hline

 $(\GSp_4\times \GL_4, S(\GSp_4\times \GL_4),std_{\Sp_4}\otimes \wedge^2\oplus T(std_{\GL_4}))$ & $T(Std_{\GSp_4}\otimes Std_{\GL_4})$  \\

\hline

 $(\GSp_4\times \GL_5, S(\GSp_4\times \GL_4),std_{\Sp_4}\otimes \wedge^2)$ & $T(Std_{\GSp_4}\otimes Std_{\GL_5})$ \\

\hline

 $(\GSpin_7\times \GL_3,\GSpin_6\times \GL_3,T(\HSpin_6\otimes std_{\GL_3}))$ & $T(Std_{\GSp_6}\otimes Std_{\GL_3})$ \\

\hline

\end{tabular}
\captionof{table}{Dual quadruples of Table \ref{Table 2 red}}
\label{Table 2 red list}
\end{figure}

\subsection{The non-reductive case}\label{sec Table 2 non-red}
For (2.1) with $m>n+1$ and (2.4) with $n>2$, the associated quadruple $\Delta$ is given by
$$(G,H,\rho_H,\iota)=(\GL_m\times \GL_n,\GL_n,0, (\GL_1^{n}\times \GL_{m-n}\times T_{\GL_n} ).$$
When $m-n$ is odd (resp. even), the nilpotent orbit induces a Bessel period (resp. Fourier-Jacobi period) for the unipotent radical of the parabolic subgroup $P=MU$ with $M=(\GL_1)^{m-n-1}\times \GL_{n+1}\times \GL_n$ (resp. $M=(\GL_1)^{m-n}\times \GL_{n}\times \GL_n$) whose stabilizer in $M$ is $\GL_n\times \GL_n$. We can diagonally embed $H$ into the stabilizer. In this case $\Delta_{red}$ is given by the quadruple \eqref{2.1 m=n+1} (resp. \eqref{2.1 m=n}). It is clear that Theorem \ref{main thm 3} holds in this case. The period integral in this case is closely related to the Rankin-Selberg integral in \cite{JPSS}. However the difference is not negligible and we do not claim Theorem \ref{main thm 1} for this case.  For the dual side, Conjecture \ref{BSV conj}(2) follows from Theorem~\ref{theta GL} applied to the theta correspondence for $\GL_n\times \GL_m$. This proves Theorem \ref{main thm 2}.

For (2.2) with $n=2m$, the associated quadruple $\Delta$ is given by
$$(\GL_{2m},\GL_m,T(std_{\GL_m}),(\GL_2)^m).$$
The nilpotent orbit induces a Bessel period for the unipotent radical of the parabolic subgroup $P=MU$ with $M=\GL_m\times \GL_m$ whose stabilizer in $M$ is $H=\GL_m$. In this case $\Delta_{red}$ is given by \eqref{2.1 m=n}. It is clear that Theorem \ref{main thm 3} holds in this case. The period integral in this case is exactly the Rankin-Selberg integral in \cite{JS90}. The result in loc. cit. proves Conjecture \ref{BSV conj}(1) and Theorem \ref{main thm 1}. 

For (2.2) with $n=2m+1$, the associated quadruple $\Delta$ is given by
$$(\GL_{2m+1},\GL_m,0,(\GL_2)^m\times \GL_1).$$
The nilpotent orbit induces a Fourier-Jacobi period for the unipotent radical of the parabolic subgroup $P=MU$ with $M=\GL_m\times \GL_1\times \GL_m$ whose stabilizer in $M$ is $\GL_n\times \GL_1$. We can naturally embed $H$ into the stabilizer. In this case $\Delta_{red}$ is given by \eqref{2.1 m=n}. It is clear that Theorem \ref{main thm 3} holds in this case. The period integral in this case is exactly the Rankin-Selberg integral in \cite{JS90}. The result in loc. cit. proves Conjecture \ref{BSV conj}(1) and Theorem \ref{main thm 1}.

For (2.5), the associated quadruple $\Delta$ is given by 
$$(\SO_{2m+1},\SO_2,0,\SO_{2m-1}\times \GL_1).$$
It is the Gross-Prasad model of $\SO_{2m+1}\times \SO_2$ and $\Delta_{red}$ is given by \eqref{1.1 p=2m} when $m=1$. It is clear that Theorem \ref{main thm 3} holds in this case. The unramified computation in \cite{II} proves Theorem \ref{main thm 1}. For the dual side, Conjecture \ref{BSV conj}(2) follows from Theorem~\ref{theta} applied to the theta correspondence  for $\Sp_{2m}\times \SO_2$ and this proves Theorem \ref{main thm 2}.

For (2.6) with $m>2, n=2$, the associated quadruple $\Delta$ is given by 
$$(G,H,\rho_H,\iota)=(\GSpin_{2m+1}\times \GL_2, G(\SL_2\times \SL_2), T(std_{\GL_2}),(\GL_1)^{2}\times \GSpin_{2m-3}\times T_{\GL_2,2}).$$
The nilpotent orbit $\iota$ induces a Bessel period on the unipotent radical of the parabolic subgroup $P=MU$ with $M=\GSpin_{5}\times (\GL_1)^{m-2}\times \GL_2$ whose stabilizer in $M$ is $\GSpin_4\times \GL_2$. We then embeds $H=G(\SL_2\times \SL_2)$ into $\GSpin_4\times \GL_2$ via the identity map on $\GSpin_4$ and the projection of $G(\SL_2\times \SL_2)$ into $\GL_2$ via the first $\GL_2$-copy. The representation $\rho_H$ is the standard representation of the second $\GL_2$-copy of $H$. This integral is essentially the Gross-Prasad model for $\GSpin_{2m+1}\times \GSpin_4$ except we replace the cusp form on one $\GL_2$-copy by theta series. 
In this case $\Delta_{red}$ is given by \eqref{2.6 m=n=2}. It is clear that Theorem \ref{main thm 3} holds in this case. The unramified computation in \cite{II} proves Theorem \ref{main thm 1}. For the dual side, Conjecture \ref{BSV conj}(2) follows from Theorem~\ref{theta} applied to the theta correspondence  for $\GSp_{2n}\times \GSO_4$ and the Rankin-Selberg integral of $\GL_2\times \GL_1$. This proves Theorem \ref{main thm 2}.

For (2.6) with $m=2, n>5$, the associated quadruple $\Delta$ is 
$$(\GSp_4\times \GL_n, S(\GSp_4\times \GL_4),std_{\Sp_4}\otimes \wedge^2,T_{\GSp_4}\times (\GL_1)^{4}\times \GL_{n-4}).$$
When $n$ is odd (resp. even),  the nilpotent orbit induces a Bessel period (resp. Fourier-Jacobi period) for the unipotent radical of the parabolic subgroup $P=MU$ with $M=\GSp_4\times \GL_5\times (\GL_1)^5$ (resp. $M=\GSp_4\times \GL_4\times (\GL_1)^4$) whose stabilizer in $M$ is $\GSp_4\times \GL_4$. We can naturally embed $H$ into the stabilizer. In this case $\Delta_{red}$ is given by \eqref{(2.6) m=2 n=5} (resp. \eqref{2.6, m=2 n=4}). It is clear that Theorem \ref{main thm 3} holds in this case. For the dual side, Conjecture \ref{BSV conj}(2) follows from Theorem~\ref{theta GL} applied to the theta correspondence of $\GL_n\times \GL_4$ and the global period integral conjecture for the pair $(\GL_4\times \GSp_4,\GSp_4)$ (which is essentially the Gross-Prasad period for $\SO_6\times \SO_5$) in \cite{GGP}. This proves Theorem \ref{main thm 2}.

For (2.6) with $m>3,n=3$, the associated quadruple $\Delta$ is given by
$$(\GSpin_{2m+1}\times \GL_3,\GSpin_6\times \GL_3,T(\HSpin_6\otimes std_{\GL_3}),(\GL_1)^{3}\times \GSpin_{2m-5}\times T_{\GL_3}).$$
The nilpotent orbit $\iota$ induces a Bessel period on the unipotent radical of the parabolic subgroup $P=MU$ with $M=\GSpin_{7}\times (\GL_1)^{m-3}\times \GL_3$ whose stabilizer in $M$ is $H=\GSpin_6\times \GL_3$. In this case $\Delta_{red}$ is given by \eqref{2.6 m=n=3}. It is clear that Theorem \ref{main thm 3} holds in this case. The unramified computation in \cite{II} and Theorem~\ref{theta GL unramified} applied to theta correspondence for $\GL_4\times \GL_3$ proves Theorem \ref{main thm 1}. For the dual side, Conjecture \ref{BSV conj}(2) follows from Theorem~\ref{theta} applied to the theta correspondence of $\GSp_{2n}\times \GSO_6$ and the Rankin-Selberg period for $\GL_4\times \GL_3$. This proves Theorem \ref{main thm 2}.

For (2.7) with $m=2k$, the associated quadruple $\Delta$ is
$$(\GSpin_{2k},\GSpin_3,T(\Spin_3),\GL_1\times \GSpin_{2k-2}).$$
This is essentially the Gross-Prasad model for $\GSpin_{2k}\times \GSpin_3$ except we replace the cusp form on $\GSpin_3$ by a theta series. In this case $\Delta_{red}$ is given by \eqref{2.1 m=n} when $n=2$. It is clear that Theorem \ref{main thm 3} holds in this case. The unramified computation in \cite{II} proves Theorem \ref{main thm 1}.

For (2.7) with $m=2k+1$, the generic stabilizer of $\hat{\rho}$ in $\hat{G}$ is not connected,  hence it does not belong to the current framework of the BZSV duality. However, for this specific case, by the Rankin-Selberg integral in \cite{GRS}, we know that the dual integral should be the one in \cite{GRS}. As the generic stabilizer is not connected, there are covering groups involved in the integral.

For (2.8) with $n=7$, the associated quadruple $\Delta$ is given by 
$$(\GSp_6,\GL_2,T(std_{\GL_2}),\GL_3\times \GL_1).$$
This is essentially the same as the quadruple \eqref{GSp(6)xGL(2)} except we replace the cusp form on $\GL_2$ by theta series. The period integral in this case is exactly the Rankin-Selberg integral in \cite{BG92a} and $\Delta_{red}$ is given by \eqref{2.1 m=n} when $m=2$. It is clear that Theorem \ref{main thm 3} holds in this case. The unramfied computation in \cite{BG92a} and \cite{WZ} proves Theorem \ref{main thm 1}.

For (2.8) with $n=9$, the associated quadruple $\Delta$ is 
\begin{equation}
(\GSp_8,G(\SL_2\times \SL_2),T(std_{\GL_2,2}),\GL_3\times \GL_1\times \GL_1).
\end{equation}
where $std_{\GL_2,2}$ is the standard representation of the second $\GL_2$-copy. This is essentially the same as the quadruple \eqref{GSp(8)xGL(2)} except we replace the cusp form on $\GL_2$ by theta series  and the period integral in this case is exactly the Rankin-Selberg integral in \cite{BG00}. In this case $\Delta_{red}$ is given by \eqref{2.6 m=n=2}. It is clear that Theorem \ref{main thm 3} holds in this case. The unramfied computation in \cite{BG00} proves Theorem \ref{main thm 1}.

For (2.8) with $n=10$, the associated quadruple $\Delta$ is 
$$(\PGSO_{10},\GL_2,0,\GL_4\times \GL_1).$$
The nilpotent orbit $\iota$ induces a Fourier-Jacobi period  on the unipotent radical of the parabolic subgroup $P=MU$ with $M=\GL_2\times \GL_2\times \SO_2$ whose stabilizer in $M$ is $H=\GL_2$ (here the embedding is given by $h\mapsto (h,h,diag(\det(h),1))$). In this case $\Delta_{red}$ is given by \eqref{2.1 m=n} when $n=2$. It is clear that Theorem \ref{main thm 3} holds in this case.   This integral is very close to the Rankin-Selberg integral in \cite{G95}, though we again do not claim Theorem~\ref{main thm 1} in this case.

For (2.9), the generic stabilizer of $\hat{\rho}$ in $\hat{G}$ is not connected,  hence it does not belong to the current framework of the BZSV duality. However, for this specific case, by the Rankin-Selberg integral in \cite{Gin93}, we know that the dual integral should be the one in \cite{Gin93}. As the generic stabilizer is not connected, there are covering groups involved in the integral.

For (2.10), the associated quadruple $\Delta$ is 
$$(GE_6, \GL_3,T(std_{\GL_3}),D_4).$$
In this case $\Delta_{red}$ is given by \eqref{2.1 m=n} when $n=3$. The period integral associated to it is exactly the Rankin-Selberg integral in \cite{G2}. It is clear that Theorem \ref{main thm 3} holds in this case.  The unramified compuation in \cite{G2}  proves Theorem \ref{main thm 1}.

By the discussion above, the strongly tempered quadruple associated to Table \ref{Table 2 non-red} is given as follows. Here for $\iota$, we only list the root type of the Levi subgroup $L$ of $G$ such that $\iota$ is principal in $L$.

\begin{figure}[h!]
\begin{tabular}{| c | c | c |}
\hline
$(G,H,\rho_H)$ & $\iota$ & $\hat{\rho}$ \\
\hline

 $(\GL_m\times \GL_n,\GL_n,0)$ & $A_{m-n-1}$ & $T(std_{\GL_m}\otimes std_{\GL_n})$ \\

\hline

$(\GL_{2m},\GL_m,T(std_{\GL_m}))$ & $(A_1)^m$ & $T(\wedge^2)$ \\
\hline

$(\GL_{2m+1},\GL_m,0)$ & $(A_1)^m$ & $T(\wedge^2)$ \\
\hline

 $(\SO_{2m+1},\SO_2,0)$ & $B_{m-1}$ & $T(std_{\Sp_{2n}})$ \\
\hline

 $(\GSpin_{2m+1}\times \GL_2, G(\SL_2\times \SL_2), T(std_{\GL_2}))$ & $B_{m-2}$ & $T(Std_{\GSp_{2m}}\otimes Std_{\GL_2})$ \\

\hline

 $(\GSp_4\times \GL_n, S(\GSp_4\times \GL_4),std_{\Sp_4}\otimes \wedge^2,(\GL_1)^5)$ & $A_{n-5}$ & $T(Std_{\Sp_4}\otimes Std_{\SL_n})$ \\

\hline

 $(\GSpin_{2m+1}\times \GL_3,\GSpin_6\times \GL_3,T(\HSpin_6\otimes std_{\GL_3}))$ & $B_{m-3}$ & $T(Std_{\Sp_{2m}}\otimes Std_{\SL_3})$ \\

\hline

$(\GSpin_{2k},\GSpin_3,T(\Spin_3))$ & $D_{k-1}$ & $T(std_{\SO_{2k}})$ \\

\hline

$(\GSp_6,\GL_2,T(std_{\GL_2}))$ & $A_2$ & $T(\Spin_7)$ \\

\hline

$(\GSp_8,G(\SL_2\times \SL_2),T(std_{\GL_2}))$ & $A_2$ & $T(\Spin_9)$ \\

\hline

$(\PGSO_{10},\GL_2,0)$ & $A_3$ & $T(\HSpin_{10})$ \\

\hline

$(GE_6, \GL_3,T(std_{\GL_3}))$ & $D_4$ & $T(std_{E_6})$ \\
\hline

\end{tabular}
\captionof{table}{Dual quadruples of Table \ref{Table 2 non-red}}
\label{Table 2 non-red list}
\end{figure}

\section{Models in Table 11}\label{sec Table 11}

In this section we will consider Table 11 of \cite{K}, this is for the case when $\hat{\rho}$ is the direct sum of two distinct irreducible symplectic representations of $\hat{G}$. It is easy to check that the representations in (11.5), (11.8), (11.13), (11.14), (11.15) of \cite{K} are not anomaly free and the representation in (11.1) (resp. (11.11)) of \cite{K} is only anomaly free when $n$ is even (resp. $p$ odd). Hence it remains to consider the following cases. We still separate the cases based on whether
$\hat{\Fl}$ is abelian or not.

\begin{figure}[h!]
\begin{tabular}{| c | c | c |c|}
\hline
Number in \cite{K} & $(\hat{G},\hat{\rho})$ & $\hat{W}_V$ &  $\hat{\Fl}$ \\
\hline

(11.7)& $(\Sp_4\times \Spin_8\times \SL_2, std_{\Sp_4}\otimes std_{\Spin_8}\oplus \HSpin_8\otimes std_{\SL_2})$ & $C_2\times D_4\times A_1$ & 0   \\

\hline

(11.9) & $(\SL_2\times \Spin_7\times SL_2, std_{\SL_2}\otimes \Spin_7\oplus \Spin_7 \otimes std_{\SL_2})$ & $(A_1)^3\times B_2$ & 0       \\

\hline

(11.10) & $(\SL_2\times \SO_6\times \SL_2, std_{\SL_2}\otimes std_{\SO_6}\oplus std_{\SO_6}\otimes std_{\SL_2})$ & $A_1\times A_1\times B_2$ & 0 \\

\hline

(11.11), p=2m+1 & $(\SO_{2m+1}\times \Sp_{2m}, std_{\SO_{2m+1}}\otimes std_{\Sp_{2m}}\oplus std_{\Sp_{2m}})$ & $B_m\times C_m$ & 0 \\

\hline

(11.11), p=2m-1 & $(\SO_{2m-1}\times \Sp_{2m}, std_{\SO_{2m-1}}\otimes std_{\Sp_{2m}}\oplus std_{\Sp_{2m}})$ & $B_{m-1}\times D_{m}$ & 0 \\

\hline

\end{tabular}
\captionof{table}{Reductive models in Table 11 of \cite{K}}
\label{Table 11 red}
\end{figure}

\begin{figure}[h!]\leftskip-2cm
\begin{tabular}{| c | c | c |c|}
\hline
Number in \cite{K} & $(\hat{G},\hat{\rho})$ & $\hat{W}_V$ &  $\hat{\Fl}$ \\
\hline

(11.1), n=2k & $(\SL_2\times \SO_{2k}\times \SL_2, std_{\SL_2}\otimes std_{\SO_{2k}}\oplus std_{\SO_{2k}}\otimes std_{\SL_2})$ & $A_1\times A_1\times B_2$ & $D_{k-2}$ \\

\hline

(11.2) & $(\Spin_{12},\HSpin_{12}^{+}\oplus \HSpin_{12}^{-})$ & $(A_1)^2\times B_2$ & $A_1\times A_1$ \\

\hline

(11.3) &  $(\SL_2\times \Spin_{12}, std_{\SL_2}\otimes std_{\Spin_{12}}\oplus \HSpin_{12})$ & $(A_1)^3$ & $A_3$  \\

\hline

(11.4) & $(\Sp_4\times \Spin_{12}, std_{\Sp_4}\otimes std_{\Spin_{12}}\oplus \HSpin_{12})$ & $C_2\times A_1\times D_4$ & $A_1$   \\

\hline

(11.6) & $(\SL_2\times \Spin_8\times SL_2, std_{\SL_2}\otimes std_{\Spin_8}\oplus \HSpin_8 \otimes std_{\SL_2})$ & $(A_1)^3$ & $A_1$   \\

\hline

(11.11), $p=2k+1>2m+1$ & $(\SO_{2k+1}\times \Sp_{2m}, std_{\SO_{2k+1}}\otimes std_{\Sp_{2m}}\oplus std_{\Sp_{2m}})$ & $B_m\times C_m$ & $B_{k-m}$ \\

\hline

(11.11), $p=2n-1<2m-1$ & $(\SO_{2n-1}\times \Sp_{2m}, std_{\SO_{2n-1}}\otimes std_{\Sp_{2m}}\oplus std_{\Sp_{2m}})$ & $B_{n-1}\times D_{n}$ & $C_{m-n}$ \\

\hline

(11.12) & $(\Sp_6, \wedge_{0}^3\oplus std_{\Sp_6})$ & $A_1\times A_1$ & $A_1$ \\

\hline

\end{tabular}
\captionof{table}{Non-reductive models in Table 11 of \cite{K}}
\label{Table 11 non-red}
\end{figure}

\subsection{The reductive case}\label{sec Table 11 red}
For (11.7), the associated quadruple $\Delta$ is 
\begin{equation}\label{11.7}
(\GSp_4\times \GSpin_8\times \GL_2,S(\GSpin_8\times G(\Sp_4\times \SL_2)), std_{\Sp_4}\otimes std_{\Spin_8}\oplus \HSpin_8\otimes std_{\SL_2}).
\end{equation}
 Note that when we take principal series on $\GSp_4$ and $\GL_2$, this period integral recovers the Rankin-Selberg integral in \cite{GH04}. The unramified computation in loc. cit. proves Theorem \ref{main thm 1} in this case. This quadruple is self-dual.

For (11.9), the associated quadruple $\Delta$ is given by
$$(\GSp_6\times \GSO_4,S(\GSO_4\times G(\Sp_4\times \SL_2)), std_{\SO_4}\times std_{\Sp_4}).$$
By the theta correspondence for $\GSO_4\times \GSp_4$, the integral over $\SO_4$ of a cusp form on $\GSO_4$ with the theta series associated to $\rho_H$ produces an automorphic form on $\GSp_4$. Then the integral over $G(\Sp_4\times \SL_2)$ is just the period integral for the pair $(\GSp_6\times \GSp_4,G(\Sp_4\times \Sp_2))$ in \cite{WZ}. The unramified computation in \cite{WZ} and Theorem~\ref{theta unramified} applied to theta correspondence for $\GSO_4\times \GSp_4$ proves Theorem \ref{main thm 1} in this case.

For (11.10), the associated quadruple $\Delta$ is given by
\begin{equation}\label{11.10}
(\GL_4\times \GSO_4,S(\GSp_4\times \GSO_4),std_{\SO_4}\times std_{\Sp_4}).
\end{equation}
By the theta correspondence for $\GSO_4\times \GSp_4$, the integral over $\SO_4$ of a cusp form on $\GSO_4$ with the theta series associated to $\rho_H$ produces an automorphic form on $\GSp_4$. Then the integral over $\GSp_4$ is just the period integral for the pair $(\GL_4\times \GSp_4,\GSp_4)$ which is essentially the Gross-Prasad model for $\SO_6\times \SO_5$. The unramified computation in \cite{II} and Theorem~\ref{theta unramified} applied to theta correspondence for $\GSO_4\times \GSp_4$ proves Theorem \ref{main thm 1} in this case. For the dual side, Conjecture \ref{BSV conj} follows from the theta correspondence for $\SO_{6}\times \Sp_{4}$ (here we view $\SL_2\times \SL_2$ as a subgroup of $\Sp_4$) and the global period integral conjecture for the Gross-Prasad model $\SO_{5}\times \SO_{4}$ in \cite{GGP}. This proves Theorem \ref{main thm 2}. 

For (11.11) when $p=2m+1$, the associated quadruple $\Delta$ is given by 
\begin{equation}\label{11.11 p=2m+1}
(\SO_{2m+1}\times \Sp_{2m},H=G, std_{\SO_{2m+1}}\otimes std_{\Sp_{2m}}\oplus std_{\Sp_{2m}}).
\end{equation}
By the theta correspondence for $\SO_{2m+2}\times \Sp_{2m}$,  the integral over  $\Sp_{2m}$ of a cusp form on $\Sp_{2m}$ with the theta series associated to $\rho_H$ produces an automorphic form on $\SO_{2m+2}$. Then the integral over $\SO_{2m+1}$ is just the period integral for the Gross-Prasad period for $\SO_{2m+2}\times \SO_{2m+1}$. The unramified computation in \cite{II} and Theorem~\ref{theta unramified} applied to theta correspondence for $\SO_{2m+2}\times \Sp_{2m}$ proves Theorem \ref{main thm 1} in this case. This quadruple is self-dual and it is clear that Conjecture \ref{BSV conj} follows from the theta correspondence for $\SO_{2m+2}\times \Sp_{2m}$ and the global period integral conjecture for the Gross-Prasad model of $\SO_{2m+2}\times \SO_{2m+1}$ in \cite{GGP}. This proves Theorem \ref{main thm 2}.

For (11.11) when $p=2m-1$, the associated quadruple $\Delta$ is given by 
\begin{equation}\label{11.11 p=2m-1}
(\SO_{2m+1}\times \Sp_{2m-2},\SO_{2m}\times \Sp_{2m-2}, std_{\SO_{2m}}\otimes std_{\Sp_{2m-2}}).
\end{equation}
By the theta correspondence for $\SO_{2m}\times \Sp_{2m-2}$,  the integral over  $\Sp_{2m-2}$ of a cusp form on $\Sp_{2m}$ with the theta series associated to $\rho_H$ produces an automorphic form on $\SO_{2m}$. Then the integral over $\SO_{2m}$ is just the Gross-Prasad period for $\SO_{2m+1}\times \SO_{2m}$. The unramified computation in \cite{II} and Theorem~\ref{theta unramified} applied to theta correspondence for $\SO_{2m}\times \Sp_{2m-2}$ proves Theorem \ref{main thm 1} in this case. For the dual side, Conjecture \ref{BSV conj} follows from the theta correspondence for $\SO_{2m}\times \Sp_{2m-2}$ and the global period integral conjecture for the Gross-Prasad model $\SO_{2m}\times \SO_{2m+1}$ in \cite{GGP}. This proves Theorem \ref{main thm 2}.

By the discussion above, the strongly tempered quadruple associated to Table \ref{Table 11 red} is given as follows (note that $\iota$ is trivial for all these cases) where
$$\ast=(\GSp_4\times \GSpin_8\times \GL_2,S(\GSpin_8\times G(\Sp_4\times \SL_2)), std_{\Sp_4}\otimes std_{\Spin_8}\oplus \HSpin_8\otimes std_{\SL_2})$$

\begin{figure}[h!]
\begin{tabular}{| c |c|}
\hline
$(G,H,\rho_H)$ & $\hat{\rho}$ \\
\hline

 $\ast$ & $std_{\Sp_4}\otimes std_{\Spin_8}\oplus \HSpin_8\otimes std_{\SL_2}$   \\

\hline

$(\GSp_6\times \GSO_4,S(\GSO_4\times G(\Sp_4\times \SL_2)), std_{\SO_4}\times std_{\Sp_4})$ & $std_{\SL_2}\otimes \Spin_7\oplus \Spin_7 \otimes std_{\SL_2}$    \\

\hline

$(\GL_4\times \GSO_4,S(\GSp_4\times \GSO_4),std_{\SO_4}\times std_{\Sp_4})$ & $std_{\SL_2}\otimes std_{\SO_6}\oplus std_{\SO_6}\otimes std_{\SL_2}$  \\

\hline

 $(\SO_{2m+1}\times \Sp_{2m},H=G, std_{\SO_{2m+1}}\otimes std_{\Sp_{2m}}\oplus std_{\Sp_{2m}})$ & $std_{\SO_{2m+1}}\otimes std_{\Sp_{2m}}\oplus std_{\Sp_{2m}}$  \\

\hline

 $(\SO_{2m+1}\times \Sp_{2m-2},\SO_{2m}\times \Sp_{2m-2}, std_{\SO_{2m}}\otimes std_{\Sp_{2m-2}})$ & $std_{\SO_{2m-1}}\otimes std_{\Sp_{2m}}\oplus std_{\Sp_{2m}}$  \\

\hline

\end{tabular}
\captionof{table}{Dual quadruples of Table \ref{Table 11 red}}
\label{Table 11 red list}
\end{figure}

\subsection{The non-reductive case}\label{sec Table 11 non-red}
For (11.1) when $n=2k$, the associated quadruple $\Delta$ is 
$$(\GSpin_{2k}\times \GSO_4,S(\GSp_4\times \GSO_4),std_{\SO_4}\times std_{\Sp_4}, \GSpin_{2k-4}\times (\GL_1)^2\times T_{\GSO_4}).$$
The nilpotent orbit $\iota$ induces a Bessel period on the unipotent radical of the parabolic subgroup $P=MU$ with $M=\GSpin_{6}\times (\GL_1)^{k-3}\times \GSO_4$ whose stabilizer in $M$ is $\GSpin_5\times \GSO_4$. We can embed $H$ into the stabilizer as in \eqref{11.10} and $\Delta_{red}$ is given by \eqref{11.10}. It is clear that Theorem \ref{main thm 3} holds in this case. The unramified computation in \cite{II} and Theorem~\ref{theta unramified} applied to theta correspondence for $\GSO_4\times \GSp_4$ proves Theorem \ref{main thm 1} in this case. For the dual side, Conjecture \ref{BSV conj} follows from the theta correspondence for $\SO_{2k}\times \Sp_{4}$ (here we view $\SL_2\times \SL_2$ as a subgroup of $\Sp_4$) and the global period integral conjecture for the Gross-Prasad model $\SO_{5}\times \SO_{4}$ in \cite{GGP}. This proves Theorem \ref{main thm 2}.

For (11.2), the associated quadruple $\Delta$ is

$$(\GSO_{12},S(\GSp_4\times \GSO_4),0, \GL_2\times \GL_2\times (\GL_1)^3).$$
The nilpotent orbit $\iota$ induces a Fourier-Jacobi period on the unipotent radical of the parabolic subgroup $P=MU$ with $M=\GL_4\times \GSO_4$ whose stabilizer in $M$ is $H$. In this case $\Delta_{red}$ is given by \eqref{11.10}. It is clear that Theorem \ref{main thm 3} holds in this case.

For (11.3), we first introduce a reductive quadruple which belongs to Table S of \cite{K}. Let $G=(\GL_2)^5$ and $H=S(\GL_2\times \GL_2\times \GL_2)$ where the embedding $H\rightarrow G$ is given by mapping the first $\GL_2$-copy into the first $\GL_2$-copy, and mapping the second (resp. third) $\GL_2$-copy diagonally into the second and third (resp. fourth and fifth) $\GL_2$-copy. Let $\rho_H=std_{\GL_2}\otimes std_{\GL_2}\otimes std_{\GL_2}$ be the triple product representation and $\iota$ be trivial. The quadruple
\begin{equation}\label{GL(2) 5}
\Delta_0=(G,H,\rho_H,\iota)=((\GL_2)^5,S(\GL_2\times \GL_2\times \GL_2),std_{\GL_2}\otimes std_{\GL_2}\otimes std_{\GL_2},1)
\end{equation}
will be used to explain several models in this paper. This quadruple comes from Table S of \cite{K}, it is obtained by combining two copies of Model (S.3) with $n=4$. We claim the dual quadruple is given by
$$\hat{\Delta}_0=(\hat{G},\widehat{G/Z_{\Delta}},\hat{\rho},1),\;\hat{\rho}=std_{\GL_2,1}\otimes std_{\GL_2,2}\otimes std_{\GL_2,3}\oplus std_{\GL_2,1}\otimes std_{\GL_2,4}\otimes std_{\GL_2,5}$$
where $std_{\GL_2,i}$ represents the standard representation of the $i$-th $\GL_2$-copy. To justify the duality, we will prove Theorem \ref{main thm 1} and Theorem \ref{main thm 2} for this case.

We start with Theorem \ref{main thm 1}. By the theta correspondence for $\GSp_2\times \GSO_4$, the integral of a cusp form on the first $\GL_2$-copy with the theta series produces cusp forms on the other two $\GL_2$-copies of $H$. Then the period integral over the remaining two copies of $\GL_2$ are just the period for two trilinear $\GL_2$-models (i.e., the first, second, third $\GL_2$-copies and the first, fourth, fifth $\GL_2$-copies ). Then Theorem \ref{main thm 1} follows from the unramified computation in \cite{II}. In fact, in this case, Conjecture \ref{BSV conj}(1) follows from the result in \cite{HK} and Theorem~\ref{theta unramified} applied to theta correspondence for $\GSp_2\times \GSO_4$. For the dual side, Conjecture \ref{BSV conj}(2) in this case is also a direct consequence of the result in \cite{HK} and Theorem~\ref{theta} applied to theta correspondence for $\GSp_2\times \GSO_4$. This proves Theorem \ref{main thm 2}. Later in Section 9, we will use a similar argument to prove Theorem \ref{main thm S} for most of the cases.

For (11.3) the associated quadruple $\Delta$ is
\begin{equation}\label{11.3}
	(\GSO_{12}\times \PGL_2,S(\GL_2\times \GSO_4),0,\GL_4\times (\GL_1)^3\times T_{\PGL_2}).
	\end{equation}
The nilpotent orbit $\iota$ induces a Fourier-Jacobi period on the unipotent radical of the parabolic subgroup $P=MU$ with $M=\GL_2\times \GL_2\times \GSO_4\times \PGL_2$ whose stabilizer in $M$ is $S(\GL_2\times \GSO_4)\times \GL_2$. We can embed $H$ into the stabilizer by mapping the $\GL_2$-copy of $H$ into the $\GL_2$-copy of the stabilizer and by mapping $\GSO_4=\GL_2\times \GL_2/\GL_{1}^{diag}$ into $\GSO_4\times \PGL_2$ via the idenity map on $\GSO_4$ and the projection map $\GSO_4=\GL_2\times \GL_2/\GL_{1}^{diag}\rightarrow \PGL_2$ via the firts $\GL_2$-copy of $\GSO_4$. It is clear that the induced embedding from $H$ into $M$ is the same as \eqref{GL(2) 5}. In this case $\Delta_{red}$ is given by \eqref{GL(2) 5}. It is clear that Theorem \ref{main thm 3} holds in this case.

For (11.4), the associated quadruple $\Delta$ is
$$(\GSp_4\times \GSpin_{12}, S(\GSpin_8\times G(\Sp_4\times \SL_2)), std_{\Sp_4}\otimes std_{\Spin_8}, T_{\GSp_4}\times \GL_2\times (\GL_1)^5).$$
The nilpotent orbit $\iota$ induces a Fourier-Jacobi period on the unipotent radical of the parabolic subgroup $P=MU$ with $M=\GSp_4\times \GL_2\times \GSpin_8$ whose stabilizer in $M$ is $\GSpin_4\times S(\GL_2\times \GSpin_8)$ and we can naturally embed $H$ into the stabilizer. In this case $\Delta_{red}$ is given by \eqref{11.7}. It is clear that Theorem \ref{main thm 3} holds in this case.

For (11.6), the associated quadruple $\Delta$ is
\begin{equation}\label{11.6}
	(\GSO_{8}\times \GSO_4,S(\GL_2\times \GSO_4),0,\GL_2\times (\GL_1)^3\times T_{\GSO_4}).
	\end{equation}
The nilpotent orbit $\iota$ induces a Fourier-Jacobi period on the unipotent radical of the parabolic subgroup $P=MU$ with $M=\GSO_4\times \GL_2\times \GSO_4$ whose stabilizer in $M$ is $S(\GSO_4\times \GL_2)\times \GSO_4$. We can embed $H$ into the stabilizer by making the $\GL_2$-copy of $H$ into the $\GL_2$-copy of the stabilizer and by mapping the $\GSO_4$-copy of $H$ diagonally into the two $\GSO_4$-copies of the stabilizer. It is clear that the induced embedding from $H$ into $M$ is the same as \eqref{GL(2) 5}. In this case $\Delta_{red}$ is given by \eqref{GL(2) 5}. It is clear that Theorem \ref{main thm 3} holds in this case.

For (11.11) when $p=2k+1>2m+1$, the associated quadruple $\Delta$ is
$$(\SO_{2m+1}\times \Sp_{2k},\SO_{2m+1}\times \Sp_{2m},std_{\SO_{2m+1}}\otimes std_{\Sp_{2m}}, T_{\SO_{2m+1}}\times \Sp_{2k-2m}\times (\GL_1)^m ).$$
The nilpotent orbit $\iota$ induces a Fourier-Jacobi period on the unipotent radical of the parabolic subgroup $P=MU$ with $M=\Sp_{2m}\times (\GL_1)^{k-m}\times \SO_{2m+1}$ whose stabilizer in $M$ is $H$. In this case $\Delta_{red}$ is given by \eqref{11.11 p=2m+1}. It is clear that Theorem \ref{main thm 3} holds in this case. For the dual side, Conjecture \ref{BSV conj}(2) follows from Theorem~\ref{theta} applied to the theta correspondence for $\Sp_{2m}\times \SO_{2k+2}$ and the Gan-Gross-Prasad conjecture (Conjecture 9.11 of \cite{GGP2}) for non-tempered Arthur packet of the Gross-Prasad model of $\SO_{2k+2}\times \SO_{2k+1}$. This proves Theorem \ref{main thm 2}.

For (11.11) when $p=2n-1<2m-1$, the associated quadruple $\Delta$ is
\begin{equation}\label{11.11 non-red}
(\SO_{2m+1}\times \Sp_{2n-2},\SO_{2n}\times \Sp_{2n-2},  std_{\SO_{2n}}\otimes std_{\Sp_{2n-2}},\SO_{2m-2n+1}\times (\GL_1)^n\times T_{\Sp_{2n-2}}).
\end{equation}
In this case $\Delta_{red}$ is given by \eqref{11.11 p=2m-1}. It is clear that Theorem \ref{main thm 3} holds in this case. By the theta correspondence for $\SO_{2n}\times \Sp_{2n-2}$,  the integral over  $\Sp_{2n-2}$ of a cusp form on $\Sp_{2n}$ with the theta series associated to $\rho_H$ produces an automorphic form on $\SO_{2n}$. Then the integral over $\SO_{2n}$ is just the Gross-Prasad period for $\SO_{2m+1}\times \SO_{2n}$.   The unramified computation in \cite{II} and Theorem~\ref{theta unramified} applied to theta correspondence for $\SO_{2n}\times \Sp_{2n-2}$ proves Theorem \ref{main thm 1} in this case. For the dual side, Conjecture \ref{BSV conj}(2) follows from the theta correspondence for $\Sp_{2m}\times \SO_{2n}$ and the global period integral conjecture for the Gross-Prasad period of $\SO_{2n}\times \SO_{2n-1}$ in \cite{GGP}. This proves Theorem \ref{main thm 2}.

For (11.12), the associated quadruple $\Delta$ is
$$(\GSpin_7, \GL_2, S(\GL_2\times \GL_2), std_{\GL_2}, \GL_2\times (\GL_1)^2).$$
The nilpotent orbit $\iota$ induces a Fourier-Jacobi period on the unipotent radical of the parabolic subgroup $P=MU$ with $M=\GSpin_3\times \GL_2$ whose stabilizer in $M$ is $H$. The representation $\rho_H$ is the standard representation on the first $\GL_2$-copy. In this case $\Delta_{red}$ is given by \eqref{11.11 p=2m+1} when $m=1$. It is clear that Theorem \ref{main thm 3} holds in this case.

By the discussion above, the strongly tempered quadruple associated to Table \ref{Table 11 non-red} is given as follows. Here for $\iota$, we only list the root type of the Levi subgroup $L$ of $G$ such that $\iota$ is principal in $L$ and
$$\ast=(\GSpin_4\times \GSpin_{12}, S(\GSpin_8\times G(\Sp_4\times \SL_2)), std_{\Sp_4}\otimes std_{\Spin_8}).$$

\begin{figure}[h!]
\begin{tabular}{| c | c | c |}
\hline
$(G,H,\rho_H)$ & $\iota$ & $\hat{\rho}$  \\
\hline

$(\GSpin_{2k}\times \GSO_4,S(\GSp_4\times \GSO_4),std_{\SO_4}\times std_{\Sp_4})$ & $D_{k-2}$ & $std_{\SL_2}\otimes std_{\SO_{2k}}\oplus std_{\SO_{2k}}\otimes std_{\SL_2}$ \\

\hline

$(\GSO_{12},S(\GSp_4\times \GSO_4),0)$ & $A_1\times A_1$ & $\HSpin_{12}^{+}\oplus \HSpin_{12}^{-}$  \\

\hline

$(\GSO_{12}\times \PGL_2,S(\GL_2\times \GSO_4),0)$ & $A_3$ &  $std_{\SL_2}\otimes std_{\Spin_{12}}\oplus \HSpin_{12}$ \\

\hline

$\ast$ & $A_1$ & $std_{\Sp_4}\otimes std_{\Spin_{12}}\oplus \HSpin_{12}$ \\

\hline

$(\GSO_{8}\times \GSO_4,S(\GL_2\times \GSO_4),0)$ & $A_1$ & $std_{\SL_2}\otimes std_{\Spin_8}\oplus \HSpin_8 \otimes std_{\SL_2}$ \\

\hline

$(\SO_{2m+1}\times \Sp_{2k},\SO_{2m+1}\times \Sp_{2m},std_{\SO_{2m+1}}\otimes std_{\Sp_{2m}})$ & $C_{k-m}$ & $std_{\SO_{2k+1}}\otimes std_{\Sp_{2m}}\oplus std_{\Sp_{2m}}$\\

\hline

$(\SO_{2m+1}\times \Sp_{2n-2},\SO_{2n}\times \Sp_{2n-2},  std_{\SO_{2n}}\otimes std_{\Sp_{2n-2}})$ & $B_{m-n}$ & $std_{\SO_{2n-1}}\otimes std_{\Sp_{2m}}\oplus std_{\Sp_{2m}}$\\

\hline

$(\GSpin_7, S(\GL_2\times \GL_2), std_{\GL_2})$ & $A_1$ & $\wedge^3\oplus std_{\Sp_6}$ \\

\hline

\end{tabular}
\captionof{table}{Dual quadruples of Table \ref{Table 11 non-red}}
\label{Table 11 non-red list}
\end{figure}

\section{Models in Table 12}\label{sec Table 12}

In this section we will consider Table 12 of \cite{K}, this is for the case when $\hat{\rho}$ is the direct sum of three irreducible representations of $\hat{G}$ with two of them dual to each other (i.e. $\hat{\rho}=\hat{\rho}_0\oplus T(\hat{\tau})$). It is easy to check that the representations in (12.4), (12.9), (12.10), (12.11), (12.11) of \cite{K} are not anomaly free. Hence it remains to consider the following cases. We still separate the cases based on whether
$\hat{\Fl}$ is abelian or not.

\begin{figure}[h!]
\begin{tabular}{| c | c | c |c|}
\hline
Number in \cite{K} & $(\hat{G},\hat{\rho})$ & $\hat{W}_V$ &  $\hat{\Fl}$ \\
\hline

(12.5) & $(\SL_6\times \SL_2, \wedge^3\oplus T(std_{\SL_6}\otimes std_{\SL_2} ))$ &  $A_1\times A_1\times A_3$ &0 \\

\hline

(12.7), m=1 & $(\SL_2\times \SL_4, std_{\SL_2}\otimes \wedge^2\oplus T(std_{\SL_4}))$ & $A_1\times A_1$ & 0   \\

\hline

(12.7), m=2 & $(\Sp_4\times \SL_4, std_{\Sp_4}\otimes \wedge^2\oplus T(std_{\SL_4}))$ & $C_2\times A_3$ & 0   \\

\hline

(12.7), m=3 & $(\Sp_6\times \Spin_6, std_{\Sp_6}\otimes std_{\Spin_6}\oplus T(\HSpin_6))$ & $A_3\times A_3$ & 0  \\

\hline

(12.8) & $(\SL_2\times \SL_4\times \SL_2, std_{\SL_2}\otimes \wedge^2\oplus T(std_{\SL_4}\otimes std_{\SL_2}))$ & $A_1\times A_1\times A_3$ & 0 \\

\hline

\end{tabular}
\captionof{table}{Reductive models in Table 12 of \cite{K}}
\label{Table 12 red}
\end{figure}

\begin{figure}[h!]
\begin{tabular}{| c | c | c |c|}
\hline
Number in \cite{K} & $(\hat{G},\hat{\rho})$ & $\hat{W}_V$ &  $\hat{\Fl}$ \\
\hline

(12.1) & $(\Spin_{12}, \HSpin_{12}\oplus T(std_{\Spin_{12}}))$ & $A_1\times A_1\times A_1$ & $A_3$ \\
\hline

(12.2) & $(\SL_2\times \Spin_{10}, std_{\SL_2}\otimes std_{\Spin_{10}}\oplus T(std_{\Spin_{10}}))$ & $A_1\times A_1\times A_3$ & $A_1$ \\

\hline

(12.3) & $(\SL_2\times \Spin_{8}, std_{\SL_2}\otimes std_{\Spin_{8}}\oplus T(std_{\Spin_{8}}))$ & $A_1\times A_1\times A_1$ & $A_1$ \\

\hline

(12.6) & $(\SL_6, \wedge^3\oplus T(std_{\SL_6} ))$ &  $A_1\times A_1$ & $A_1\times A_1$ \\

\hline

(12.7), $m>3$ & $(\Sp_{2m}\times \SO_6, std_{\Sp_{2m}}\otimes std_{\SO_6}\oplus T(\HSpin_6))$ & $A_3\times A_3$ & $C_{m-3}$  \\

\hline

\end{tabular}
\captionof{table}{Non-reductive models in Table 12 of \cite{K}}
\label{Table 12 non-red}
\end{figure}

\subsection{The reductive case}\label{sec Table 12 red}

For (12.5), the associated quadruple $\Delta$ is
$$(\GL_6\times \GL_2,\GL_2\times S(\GL_4\times \GL_2),\wedge^2\otimes std_{\GL_2}).$$
At this moment we do not have much evidence that the above is the dual quadruple other than the fact that $\wedge^2\otimes std_{\GL_2}$ is the only feasible choice of symplectic representation. We believe an unramified computation similar to \cite{II} and \cite{WZ} can confirm the duality in this case. 

For (12.7) with $m=1$, the associated quadruple $\Delta$ is
$$(\GL_4\times\GL_2,\GL_2\times \GL_2,0).$$
This is the model $(\GL_4\times \GL_2,\GL_2\times \GL_2)$ studied in \cite{WZ} and the unramified computation in \cite{WZ} proves Theorem \ref{main thm 1} in this case.  For the dual side, Conjecture \ref{BSV conj}(2) follows from Theorem~\ref{theta} applied to the theta correspondence of $\GSp_{2}\times \GSO_{6}$ and Gan-Gross-Prasad conjecture (Conjecture 9.11 of \cite{GGP2}) for non-tempered Arthur packet of the Rankin-Selberg integral of $\GL_4\times \GL_4$.  This proves Theorem \ref{main thm 2}.

For (12.7) with $m=2$, the associated quadruple $\Delta$ is 
$$(\GL_4\times \GSp_4, \GL_4\times \GSp_4, T(std_{\GL_4}\otimes std_{\GSp_4})).$$
Observe that this is the dual to the quadruple in \eqref{2.6, m=2 n=4}, thus both Theorems~\ref{main thm 1} and \ref{main thm 2} have been proved there. 


For (12.7) with $m=3$, the associated quadruple $\Delta$ is
\begin{equation}\label{12.7}
(\GSpin_7\times \GSpin_6,\GSpin_6\times \GSpin_6,T(\HSpin_6\otimes \HSpin_6)).
\end{equation}
By the theta correspondence for $\GL_4\times \GL_4$, the integral over the second $\GSpin_6$-copy of a cusp form on $\GSpin_6$ with the theta series associated to $\rho_H$ produces the same cusp form with an extra central value of the Spin L-function. Then the integral over the other copy of $\GSpin_6$ is just the period integral for the Gross-Prasad model $\GSpin_7\times \GSpin_6$. The unramified computation in \cite{II} and Theorem~\ref{theta GL unramified} applied to theta correspondence for $\GL_4\times \GL_4$ proves Theorem \ref{main thm 1} in this case. For the dual side, Conjecture \ref{BSV conj}(2) follows from the theta correspondence for $\GSp_6\times \GSO_6$ and the Rankin-Selberg integral of $\GL_4\times \GL_4$. This proves Theorem \ref{main thm 2}.

For (12.8), the associated quadruple $\Delta$ is
\begin{equation}\label{12.8}
(\GL_2\times \GL_4\times \GL_2,S(\GL_2\times \GL_4)\times \GL_2, std_{\GL_2}\otimes \wedge^2\oplus T(std_{\GL_4}\times std_{\GL_2})).
\end{equation}
 Note that when we put principal series on both $\GL_2$ copies, this period integral recovers the Rankin-Selberg integral in \cite{PS}. The unramified computation in \cite{PS} proves Theorem \ref{main thm 1} in this case. This quadruple is self-dual.

By the discussion above, the strongly tempered quadruple associated to Table \ref{Table 12 red} is given as follows ($\iota$ is trivial for all these cases) where
$$\ast=(\GL_2\times \GL_4\times \GL_2,S(\GL_2\times \GL_4)\times \GL_2, std_{\GL_2}\otimes \wedge^2\oplus T(std_{\GL_4}\times std_{\GL_2})).$$

\begin{figure}[h!]
\begin{tabular}{|  c |c|}
\hline
$(G,H,\rho_H)$ & $\hat{\rho}$  \\
\hline

$(\GL_6\times \GL_2,\GL_2\times S(\GL_4\times \GL_2),\wedge^2\otimes std_{\GL_2})$ & $\wedge^3\oplus T(std_{\SL_6}\otimes std_{\SL_2} )$ \\

\hline

$(\GL_4\times\GL_2,\GL_2\times \GL_2,0)$ & $std_{\SL_2}\otimes \wedge^2\oplus T(std_{\SL_4})$ \\

\hline

$(\GL_4\times \GSp_4, \GL_4\times \GSp_4, T(std_{\GL_4}\otimes std_{\GSp_4}))$ & $std_{\Sp_4}\otimes \wedge^2\oplus T(std_{\SL_4})$  \\

\hline

$(\GSpin_7\times \GSpin_6,\GSpin_6\times \GSpin_6,T(\HSpin_6\otimes \HSpin_6))$ & $std_{\Sp_6}\otimes std_{\Spin_6}\oplus T(\HSpin_6)$  \\

\hline

$\ast$ & $std_{\SL_2}\otimes \wedge^2\oplus T(std_{\SL_4}\otimes std_{\SL_2})$ \\

\hline

\end{tabular}
\captionof{table}{Dual quadruples of Table \ref{Table 12 red}}
\label{Table 12 red list}
\end{figure}

\subsection{The non-reductive case}\label{sec Table 12 non-red}
For (12.1), we first introduce a reductive quadruple which belongs to Table S of \cite{K}. Let $G=(\GL_2)^4$ and $H=S(\GL_2\times \GL_2\times \GL_2)$ where the embedding $H\rightarrow G$ is given by mapping the first two $\GL_2$-copies into the first two $\GL_2$-copy, and mapping the last $\GL_2$-copy diagonally into the third and fourth $\GL_2$-copy. Let $\rho_H=std_{\GL_2}\otimes std_{\GL_2}\otimes std_{\GL_2}\oplus T(std_{\GL_2,2})$ where $std_{\GL_2,i}$ represents the standard representation of the $i$-th $\GL_2$-copy and $\iota$ be trivial. This quadruple
\begin{equation}\label{GL(2) 4}
\Delta_0=(G,H,\rho_H,\iota)=((\GL_2)^4,S(\GL_2\times \GL_2\times \GL_2),std_{\GL_2}\otimes std_{\GL_2}\otimes std_{\GL_2}\oplus T(std_{\GL_2,2}),1)
\end{equation}
is almost the same as \eqref{GL(2) 5} except we replace the cusp form on one $\GL_2$-copy by theta series. It is obtained by combining Model (S.3) and (S.11) in Table S of \cite{K} with $n=4$ and $m=2$. We claim the dual quadruple is given by
$$\hat{\Delta}_0=(\hat{G},\widehat{G/Z_{\Delta}},\hat{\rho},1),\;\hat{\rho}=T(std_{\GL_2,1}\otimes std_{\GL_2,2})\oplus std_{\GL_2,1}\otimes std_{\GL_2,3}\otimes std_{\GL_2,4}.$$
We can use the same argument as in \eqref{GL(2) 5} to prove Theorem \ref{main thm 1} and Theorem \ref{main thm 2} for this case.

For (12.1), the associated quadruple $\Delta$ is
$$(\GSO_{12},S(\GL_2\times \GSO_4),T(std_{\GL_2}),\GL_4\times (\GL_1)^3).$$
The attached period integral is the same as model in \eqref{11.3} except we replace the cusp form on $\GL_2$ by theta series. In this case $\Delta_{red}$ is given by \eqref{GL(2) 4} and it is clear that Theorem \ref{main thm 3} holds in this case.

For (12.2), the associated quadruple $\Delta$ is
$$(\GSpin_{10}\times \GL_2, S(\GL_2\times \GSpin_6)\times \GL_2, T(\HSpin_6\otimes std_{\GL_2}),\GL_2\times (\GL_1)^{4}\times T_{\GL_2})$$
The nilpotent orbit $\iota$ induces a Fourier-Jacobi period on the unipotent radical of the parabolic subgroup $P=MU$ with $M=\GL_2\times \GSpin_6\times \GL_2$ whose stabilizer in $M$ is $H$. In this case $\Delta_{red}$ is given by \eqref{12.8}. It is clear that Theorem \ref{main thm 3} holds in this case.

For (12.3), the associated quadruple $\Delta$ is
$$(\GSO_{8}\times \GL_2,S(\GL_2\times \GSO_4),T(std_{\GL_2}),\GL_2\times (\GL_1)^3\times T_{\GL_2}).$$
The attached period integral is the same as the model \eqref{11.6} except we replace the cusp form on one $\GL_2$-copy by theta series. In this case $\Delta_{red}$ is given by \eqref{GL(2) 4} and it is clear that Theorem \ref{main thm 3} holds in this case.

For (12.6), we first introduce a reductive quadruple from Table S of \cite{K} (it is obtained by combining Model (S.10) and Model (S.3) with $n=4$)
\begin{equation}\label{GL(2) 3}
(G,H,\rho_H,\iota)=(\GL_2\times \GL_2\times \GL_2,\GL_2\times \GL_2,T(std_{\GL_2}\otimes std_{\GL_2}),1)
\end{equation}
where $H$ embeds into $G$ by mapping the first $\GL_2$-copy into the first $\GL_2$-copy and mapping the second $\GL_2$-copy diagonally into the second and third $\GL_2$-copy. We claim the dual quadruple is given by
$$(\hat{G},\widehat{G/Z_{\Delta}},\hat{\rho},1),\;\hat{\rho}=T(std_{\GL_2,1})\oplus std_{\GL_2,1}\otimes std_{\GL_2,2}\otimes std_{\GL_2,3}$$
where $std_{\GL_2,i}$ is the standard representation of the $i$-th $\GL_2$-copy. To justify the duality, we will prove Theorem \ref{main thm 1} and Theorem \ref{main thm 2} for this case.

We start with Theorem \ref{main thm 1}. By the theta correspondence for $\GL_2\times \GL_2$, the integral over the first $\GL_2$-copy of a cusp form in $\pi$ with the theta series gives a cusp form on $\GL_2$ (in the same space $\pi$, note though Theorem~\ref{theta GL} applied to the correspondence does introduce the central value of the standard L-function). Then the integral over the other $\GL_2$-copy is just the period integral for the trilinear $\GL_2$-model. As a result, Conjecture \ref{BSV conj}(1) and Theorem \ref{main thm 1} follow from the theta correspondence for $\GL_2\times \GL_2$ and the result in \cite{HK}. For the dual side, Conjecture \ref{BSV conj}(2) follows from the theta correspondence for $\GSp_2\times \GSO_4$ and the Rankin-Selberg integral of $\GL_2\times \GL_2$. This proves Theorem \ref{main thm 2} in this case. Later in Section 9, we will use a similar argument to prove Theorem \ref{main thm S} for some of the cases (more precisely for those cases containing model (S.10) of \cite{K}).

Now we can write down the associated quadruple $\Delta$ of (12.6). It is given by
$$(\GL_6,\GL_2\times \GL_2,0,\GL_2\times \GL_2\times \GL_1\times \GL_1).$$
The nilpotent orbit $\iota$ induces a Fourier-Jacobi period on the unipotent radical of the parabolic subgroup $P=MU$ with $M=\GL_2\times \GL_2\times \GL_2$ whose stabilizer in $M$ is $H$. In this case $\Delta_{red}$ is given by \eqref{GL(2) 3}. It is clear that Theorem \ref{main thm 3} holds in this case.

For (12.7) when $m>3$, the associated quadruple $\Delta$ is
$$(\GSpin_{2m+1}\times \GSpin_6,\GSpin_6\times \GSpin_6,T(\HSpin_6\otimes \HSpin_6),\GSpin_{2m-5}\times (\GL_1)^3\times (\GL_1)^4).$$
The nilpotent orbit $\iota$ induces a Bessel period on the unipotent radical of the parabolic subgroup $P=MU$ with $M=\GL_{1}^{m-3}\times \GSpin_7\times \GSpin_6$ whose stabilizer in $M$ is $H$. In this case $\Delta_{red}$ is given by \eqref{12.7}. It is clear that Theorem \ref{main thm 3} holds in this case. The unramified computation in \cite{II} and Theorem~\ref{theta GL unramified} applied to theta correspondence for $\GL_4\times \GL_4$ proves Theorem \ref{main thm 1} in this case. For the dual side, Conjecture \ref{BSV conj}(2) follows from the theta correspondence for $\GSp_{2m}\times \GSO_6$ and the Rankin-Selberg integral of $\GL_4\times \GL_4$. This proves Theorem \ref{main thm 2}.

By the discussion above, the strongly tempered quadruple associated to Table \ref{Table 12 non-red} is given as follows. Here for $\iota$, we only list the root type of the Levi subgroup $L$ of $G$ such that $\iota$ is principal in $L$ and
$$\ast=(\GSpin_{10}\times \GL_2, S(\GL_2\times \GSpin_6)\times \GL_2, T(\HSpin_6\otimes std_{\GL_2})).$$

\begin{figure}[h!]\leftskip-1cm
\begin{tabular}{| c | c | c |}
\hline
$(G,H,\rho_H)$ & $\iota$ & $\hat{\rho}$ \\
\hline

$(\GSO_{12},S(\GL_2\times \GSO_4),T(std_{\GL_2}))$ & $A_3$ & $\HSpin_{12}\oplus T(std_{\Spin_{12}})$   \\
\hline

$\ast$ & $A_1$ & $std_{\SL_2}\otimes std_{\Spin_{10}}\oplus T(std_{\Spin_{10}})$  \\

\hline

$(\GSO_{8}\times \GL_2,S(\GL_2\times \GSO_4),T(std_{\GL_2}))$ & $A_1$ & $std_{\SL_2}\otimes std_{\Spin_{8}}\oplus T(std_{\Spin_{8}})$  \\

\hline

$(\GL_6,\GL_2\times \GL_2,0)$ & $A_1\times A_1$ & $\wedge^3\oplus T(std_{\SL_6} )$  \\

\hline

$(\GSpin_{2m+1}\times \GSpin_6,\GSpin_6\times \GSpin_6,T(\HSpin_6\otimes \HSpin_6))$ & $B_{m-3}$ & $std_{\Sp_{2m}}\otimes std_{\SO_6}\oplus T(\HSpin_6)$ \\

\hline

\end{tabular}
\captionof{table}{Dual quadruples of Table \ref{Table 12 non-red}}
\label{Table 12 non-red list}
\end{figure}

\section{Models in Table 22}\label{sec Table 22}
In this section we will consider Table 22 of \cite{K}, this is for the case when $\hat{\rho}$ is the direct sum of four irreducible representations of $\hat{G}$ of the form $T(\rho_1)\oplus T(\rho_2)$. All the representations in Table 22 of \cite{K} are anomaly free, so we need to consider all of them. We still separate the cases based on whether
$\hat{\Fl}$ is abelian or not.

\begin{figure}[h!]
\begin{tabular}{| c | c | c |c|}
\hline
Number in \cite{K} & $(\hat{G},\hat{\rho})$ & $\hat{W}_V$ &  $\hat{\Fl}$ \\
\hline

(22.2), n=2m & $(\SL_n,T(\wedge^2)\oplus T(std_{\SL_n}))$ & $A_{m-1}\times A_{m-1}$ & 0 \\

\hline

(22.2), n=2m+1 & $(\SL_n,T(\wedge^2)\oplus T(std_{\SL_n}))$ & $A_{m}\times A_{m-1}$ & 0 \\

\hline

(22.3), m=n & $(\SL_m\times \SL_n,T(std_{\SL_m}\otimes std_{\SL_n})\oplus T(std_{\SL_n}))$ & $A_{n-1}\times A_{n-1}$ & 0 \\

\hline

(22.3), m=n+1 & $(\SL_m\times \SL_n,T(std_{\SL_m}\otimes std_{\SL_n})\oplus T(std_{\SL_n}))$ & $A_{n-1}\times A_{n-1}$ & 0 \\

\hline

(22.3), m=n-1 & $(\SL_m\times \SL_n,T(std_{\SL_m}\otimes std_{\SL_n})\oplus T(std_{\SL_n}))$ & $A_{m}\times A_{m-1}$ & 0 \\

\hline

(22.3), m=n-2 & $(\SL_m\times \SL_n,T(std_{\SL_m}\otimes std_{\SL_n})\oplus T(std_{\SL_n}))$ & $A_{m}\times A_{m-1}$ & 0 \\

\hline

(22.4), n=3 & $(\SL_3, T(std_{\SL_3})\oplus T(std_{\SL_3}))$ & $A_1$ & 0 \\

\hline

(22.5), m=2 & $(\Sp_{4}, T(std_{\Sp_{4}})\oplus T(std_{\Sp_{4}}))$ & $A_1\times A_1$ & 0 \\

\hline

\end{tabular}
\captionof{table}{Reductive models in Table 22 of \cite{K}}
\label{Table 22 red}
\end{figure}

\newpage

\begin{figure}[h!]
\begin{tabular}{| c | c | c |c|}
\hline
Number in \cite{K} & $(\hat{G},\hat{\rho})$ & $\hat{W}_V$ &  $\hat{\Fl}$ \\
\hline

(22.1) & $(\Spin_8, T(std_{\Spin_8})\oplus  T(\HSpin_8))$ & $A_1\times A_1\times A_1$ & $A_1$ \\

\hline

(22.3), $m>n+1$ & $(\SL_m\times \SL_n,T(std_{\SL_m}\otimes std_{\SL_n})\oplus T(std_{\SL_n}))$ & $A_{n-1}\times A_{n-1}$ & $A_{m-n+1}$ \\

\hline

(22.3), $m<n-2$ & $(\SL_m\times \SL_n,T(std_{\SL_m}\otimes std_{\SL_n})\oplus T(std_{\SL_n}))$ & $A_{m}\times A_{m-1}$ & $A_{n-m-2}$ \\

\hline

(22.4), $n>3$ & $(\SL_n, T(std_{\SL_n})\oplus T(std_{\SL_n}))$ & $A_1$ & $A_{n-3}$ \\

\hline

(22.5), $m>2$ & $(\Sp_{2m}, T(std_{\Sp_{2m}})\oplus T(std_{\Sp_{2m}}))$ & $A_1\times A_1$ & $C_{m-2}$ \\

\hline

\end{tabular}
\captionof{table}{Non-reductive models in Table 22 of \cite{K}}
\label{Table 22 non-red}
\end{figure}

\subsection{The reductive case}\label{sec Table 22 red}
For (22.2) with $n=2m$, the associated quadruple $\Delta$ is 
$$(\GL_{2m},\GL_m\times \GL_m,T(std_{\GL_m})).$$
The period integral in this case is exactly the Rankin-Selberg integral in \cite{BF90}. The result in loc. cit. proves Conjecture \ref{BSV conj}(1) and Theorem \ref{main thm 1} in this case. 

For (22.2) with $n=2m+1$, the associated quadruple $\Delta$ is 
$$(\GL_{2m+1},\GL_{m+1}\times \GL_m,T(std_{\GL_{m+1}})).$$
The period integral in this case is exactly the Rankin-Selberg integral in \cite{BF90}. The unramified computation in loc. cit. proves Conjecture \ref{BSV conj}(1) and Theorem \ref{main thm 1} in this case.

For (22.3) with $m=n$, the associated quadruple $\Delta$ is 
\begin{equation}\label{22.3 m=n}
(\GL_n\times \GL_n, \GL_n\times \GL_n, T(std_{\GL_n}\otimes std_{\GL_n}\oplus std_{\GL_n})).
\end{equation}
By the theta correspondence for $\GL_n\times \GL_n$, the integral over $\GL_n$ of a cusp form on $\GL_n$ with the theta series associated to $T(std_{\GL_n}\otimes std_{\GL_n})$ produces a cusp form on $\GL_n$. Then the integral over the other $\GL_n$-copy is just the Rankin-Selberg integral of $\GL_{n}\times \GL_n$. This quadruple is self-dual. The Rankin-Selberg integral of $\GL_n\times \GL_n$ and Theorems~\ref{theta GL} and \ref{theta GL unramified} applied to the theta correspondence for $\GL_n\times \GL_n$ proves Conjecture \ref{BSV conj}, Theorem \ref{main thm 1} and Theorem \ref{main thm 2}. Notice that the theta correspondence introduces an extra central value of the standard $L$-function in this case.

For (22.3) with $m=n+1$, the associated quadruple $\Delta$ is 
\begin{equation}\label{22.3 m=n+1}
(\GL_{n+1}\times \GL_n, \GL_n\times \GL_n, T(std_{\GL_n}\otimes std_{\GL_n})).
\end{equation}
By the theta correspondence for $\GL_n\times \GL_n$, the integral over $\GL_n$ of a cusp form on $\GL_n$ with the theta series associated to $\rho_H$ produces another cusp form on $\GL_n$. Then the integral over the other $\GL_n$-copy is just the Rankin-Selberg integral of $\GL_{n+1}\times \GL_n$. The Rankin-Selberg integral of $\GL_{n+1}\times \GL_n$ in \cite{JPSS} and Theorems~\ref{theta GL} and \ref{theta GL unramified} applied to the theta correspondence for $\GL_n\times \GL_n$ proves Conjecture \ref{BSV conj}(1) and Theorem \ref{main thm 1} in this case. Again notice that the theta correspondence introduces an extra central value of the standard $L$-function. For the dual side, Conjecture \ref{BSV conj}(2) follows from the theta correspondence of $\GL_{n+1}\times \GL_n$ with the Rankin-Selberg integral of $\GL_n\times \GL_n$. This proves Theorem \ref{main thm 2}.

For (22.3) with $m=n-1$, the associated quadruple $\Delta$ is 
\begin{equation}\label{22.3 m=n-1}
(\GL_n\times \GL_{n-1}, \GL_n\times \GL_{n-1}, T(std_{\GL_n}\otimes std_{\GL_{n-1}}\oplus std_{\GL_n})).
\end{equation}
By the theta correspondence for $\GL_n\times \GL_n$, the integral over $\GL_n$ of a cusp form on $\GL_n$ with the theta series associated to $\rho_H$ produces another cusp form on $\GL_n$. Then the integral over $\GL_{n-1}$ is just the Rankin-Selberg integral of $\GL_{n}\times \GL_{n-1}$. This quadruple is self-dual. The Rankin-Selberg integral of $\GL_n\times \GL_{n-1}$ and Theorems~\ref{theta GL} and \ref{theta GL unramified} applied to the theta correspondence for $\GL_n\times \GL_n$ proves Conjecture \ref{BSV conj}, Theorem \ref{main thm 1} and Theorem \ref{main thm 2} in this case. As before, the theta correspondence introduces an extra central value of the standard $L$-function.

For (22.3) with $m=n-2$, the associated quadruple $\Delta$ is 
\begin{equation}\label{22.3 m=n-2}
(\GL_n\times \GL_{n-2}, \GL_{n-1}\times \GL_{n-2}, T(std_{\GL_{n-1}}\otimes std_{\GL_{n-2}})).
\end{equation}
By the theta correspondence for $\GL_{n-1}\times \GL_{n-2}$, the integral over $\GL_{n-2}$ of a cusp form on $\GL_{n-2}$ with the theta series associated to $\rho_H$ produces an Eisenstein series on $\GL_{n-1}$ which is induced from the cuspidal automorphic representation on $\GL_{n-2}$ and the trivial character. Then the integral over $\GL_{n-1}$ is just the Rankin-Selberg integral of $\GL_{n}\times \GL_{n-1}$. The Rankin-Selberg integral of $\GL_{n-1}\times \GL_n$ in \cite{JPSS} and Theorems~\ref{theta GL} and \ref{theta GL unramified} applied to the theta correspondence for $\GL_{n-1}\times \GL_{n-2}$ proves Conjecture \ref{BSV conj}(1) and Theorem \ref{main thm 1} in this case. For the dual side, Conjecture \ref{BSV conj}(2) follows from the theta correspondence of $\GL_{n-1}\times \GL_n$ with the Rankin-Selberg integral of $\GL_{n-1}\times \GL_{n-2}$. This proves Theorem \ref{main thm 2}.

For (22.4) with $n=3$, the associated quadruple $\Delta$ is 
\begin{equation}\label{22.4}
(\GL_3,\GL_2\times \GL_1,T(std_{\GL_2})).
\end{equation}
The period integral is essentially the Rankin-Selberg integral of $\GL_{3}\times \GL_{2}$ except that we replace the cusp form on $\GL_2$ by theta series. The result in \cite{JPSS} proves Conjecture \ref{BSV conj}(1) and Theorem \ref{main thm 1} in this case.

For (22.5) with $m=2$, the associated quadruple $\Delta$ is 
\begin{equation}\label{22.5}
(\GSpin_5\times \GL_1,\GSpin_4\times \GL_1,T(\HSpin_{4}^{+}\oplus \HSpin_{4}^{-}\otimes std_{\GL_1})).
\end{equation}
The period integral is essentially the Gross-Prasad period for $\GSpin_5\times \GSpin_4$ except that we replace the cusp form on $\GSpin_4$ by theta series. The unramified computation in \cite{II} proves Theorem \ref{main thm 1} in this case.

By the discussion above, the strongly tempered quadruple associated to Table \ref{Table 12 red} is given as follows ($\iota$ is trivial for all these cases).

\begin{figure}[h!]
\begin{tabular}{| c | c |}
\hline
$(G,H,\rho_H)$ & $\hat{\rho}$  \\
\hline

$(\GL_{2m},\GL_m\times \GL_m,T(std_{\GL_m}))$ & $T(\wedge^2)\oplus T(std_{\GL_{2m}})$  \\

\hline

$(\GL_{2m+1},\GL_{m+1}\times \GL_m,T(std_{\GL_{m+1}}))$ & $T(\wedge^2)\oplus T(std_{\GL_{2m+1}})$ \\

\hline

$(\GL_n\times \GL_n, \GL_n\times \GL_n, T(std_{\GL_n}\otimes std_{\GL_n}\oplus std_{\GL_n}))$ & $T(std_{\GL_n}\otimes std_{\GL_n})\oplus T(std_{\GL_n})$\\

\hline

$(\GL_{n+1}\times \GL_n, \GL_n\times \GL_n, T(std_{\GL_n}\otimes std_{\GL_n}))$ & $T(std_{\GL_{n+1}}\otimes std_{\GL_n})\oplus T(std_{\GL_n})$ \\

\hline

$(\GL_n\times \GL_{n-1}, \GL_n\times \GL_{n-1}, T(std_{\GL_n}\otimes std_{\GL_{n-1}}\oplus std_{\GL_{n}}))$ & $T(std_{\GL_{n}}\otimes std_{\GL_{n-1}})\oplus T(std_{\GL_n})$ \\

\hline

$(\GL_n\times \GL_{n-2}, \GL_{n-1}\times \GL_{n-2}, T(std_{\GL_{n-1}}\otimes std_{\GL_{n-2}}))$ & $T(std_{\GL_{n}}\otimes std_{\GL_{n-2}})\oplus T(std_{\GL_n})$  \\

\hline

$(\GL_3,\GL_2\times \GL_1,T(std_{\GL_2}))$ & $T(std_{\SL_3})\oplus T(std_{\SL_3})$ \\

\hline

$(\GSpin_5\times \GL_1,\GSpin_4\times \GL_1,T(\HSpin_{4}^{+}\oplus \HSpin_{4}^{-}\otimes std_{\GL_1}))$ & $T(std_{\Sp_{4}})\oplus T(std_{\Sp_{4}})$  \\

\hline

\end{tabular}
\captionof{table}{Dual quadruples of Table \ref{Table 22 red}}
\label{Table 22 red list}
\end{figure}

\subsection{The non-reductive case}\label{sec Table 22 non-red}
For (22.1), we first introduce a reductive quadruple which belongs to Table S of \cite{K}. Let $G=(\GL_2)^3$, $H=S(\GL_2\times \GL_2\times \GL_2)$ and $\rho_H=std_{\GL_2}\otimes std_{\GL_2}\otimes std_{\GL_2}\oplus T(std_{\GL_2,2})\oplus T(std_{\GL_2,3})$ where $std_{\GL_2,i}$ represents the standard representation of the $i$-th $\GL_2$-copy and $\iota$ be trivial. This quadruple
\begin{equation}\label{GL(2) 31}
\Delta_0=(G,H,\rho_H,\iota)=((\GL_2)^3,S(\GL_2\times \GL_2\times \GL_2),std_{\GL_2}\otimes std_{\GL_2}\otimes std_{\GL_2}\oplus T(std_{\GL_2,2}\oplus T(std_{\GL_2,3}),1)
\end{equation}
is almost the same as \eqref{GL(2) 5} except we replace the cusp form on two $\GL_2$-copies by theta series. It is obtained by combining two copies of Model (S.11) in Table S of \cite{K} with $m=2$. We claim the dual quadruple is given by
$$\hat{\Delta}_0=(\hat{G},\widehat{G/Z_{\Delta}},\hat{\rho},1),\;\hat{\rho}=T(std_{\GL_2,1}\otimes std_{\GL_2,2})\oplus T(std_{\GL_2,1}\otimes std_{\GL_2,3}).$$
We can use the same argument as in \eqref{GL(2) 5} to prove Theorem \ref{main thm 1} and Theorem \ref{main thm 2} for this case.

For (22.1), the associated quadruple $\Delta$ is
$$(\GSO_{8},S(\GL_2\times \GSO_4),T(std_{\GL_2}\oplus std_{\GL_2}),\GL_2\times (\GL_1)^3).$$
The period integral is the same as \eqref{11.6} except we replace the cusp form on both $\GL_2$-copies by theta series. In this case $\Delta_{red}$ is given by \eqref{GL(2) 31} and it is clear that Theorem \ref{main thm 3} holds in this case.

For (22.3) when $m>n+1$, the associated quadruple $\Delta$ is
$$(\GL_{m}\times \GL_n,\GL_n\times \GL_n,T(std_{\GL_n}\otimes std_{\GL_n}),(\GL_1)^n\times \GL_{m-n}\times T_{\GL_n}).$$
When $n-m$ is odd (resp. even), the nilpotent orbit $\iota$ induces a Bessel period (resp. Fourier-Jacobi period) on the unipotent radical of the parabolic subgroup $P=MU$ with $M=\GL_{1}^{m-n-1}\times \GL_{n+1}\times \GL_n$ (resp. $M=\GL_{1}^{m-n}\times \GL_{n}\times \GL_n$) whose stabilizer in $M$ is $H$. In this case $\Delta_{red}$ is given by \eqref{22.3 m=n+1} (resp. \eqref{22.3 m=n}). It is clear that Theorem \ref{main thm 3} holds in this case. For the dual side, Conjecture \ref{BSV conj}(2) follows from Theorem~\ref{theta GL} applied to the theta correspondence of $\GL_{n}\times \GL_{m+1}$ and Gan-Gross-Prasad conjecture (Conjecture 9.11 of \cite{GGP2}) for non-tempered Arthur packet of the Rankin-Selberg integral of $\GL_{m+1}\times \GL_m$. This proves Theorem \ref{main thm 2}.

For (22.3) when $m<n-2$, the associated quadruple $\Delta$ is 
$$(\GL_{m}\times \GL_n,\GL_m\times \GL_{m+1},T(std_{\GL_m}\otimes std_{\GL_{m+1}}),T_{\GL_m}\times (\GL_1)^{m-1}\times \GL_{n-m-1}).$$
When $n-m-1$ is odd (resp. even), the nilpotent orbit $\iota$ induces a Bessel period (resp. Fourier-Jacobi period) on the unipotent radical of the parabolic subgroup $P=MU$ with $M=\GL_{1}^{n-m-2}\times \GL_{m+2}\times \GL_m$ (resp. $M=\GL_{1}^{n-m-1}\times \GL_{m+1}\times \GL_m$) whose stabilizer in $M$ is $H$. In this case $\Delta_{red}$ is given by \eqref{22.3 m=n-2} (resp. \eqref{22.3 m=n-1}). It is clear that Theorem \ref{main thm 3} holds in this case. For the dual side, Conjecture \ref{BSV conj}(2) follows from Theorem~\ref{theta GL} applied to the theta correspondence of $\GL_{n}\times \GL_{m+1}$ and the Rankin-Selberg integral of $\GL_{m+1}\times \GL_{m}$. This proves Theorem \ref{main thm 2}.

For (22.4) when $n>3$, we need to introduce another reductive quadruple from Table S of \cite{K} (it is obtained by combining two copies of Model (S.10))
\begin{equation}\label{GL(2) 1}
(G,H,\rho_H,\iota)=(\GL_2\times \GL_1,\GL_2\times \GL_1,T(std_{\GL_2}\oplus std_{\GL_2}\otimes std_{\GL_1}),1).
\end{equation}
We claim that the dual quadruple is given by 
$$(\hat{G},\hat{G},\hat{\rho},1),\;\hat{\rho}=T(std_{\GL_2}\oplus std_{\GL_2}\otimes std_{\GL_1}),$$
i.e., it is self-dual. We can use the same argument as in \eqref{GL(2) 3} to prove Theorem \ref{main thm 1} and Theorem \ref{main thm 2} for this case.

The associated quadruple $\Delta$ for (22.4) with $n>3$ is given by 
$$(\GL_n,\GL_2, T(std_{\GL_2}),\GL_{n-2}\times \GL_1\times \GL_1).$$
When $n-2$ is odd (resp. even), the nilpotent orbit $\iota$ induces a Bessel period (resp. Fourier-Jacobi period) on the unipotent radical of the parabolic subgroup $P=MU$ with $M=\GL_{1}^{n-3}\times \GL_3$ (resp. $M=\GL_{1}^{n-2}\times \GL_2$) whose stabilizer in $M$ is $H$. In this case $\Delta_{red}$ is given by \eqref{22.4} (resp. \eqref{GL(2) 1}). It is clear that Theorem \ref{main thm 3} holds in this case.

For (22.5) when $m>2$, the associated quadruple $\Delta$ is 
$$(\GSpin_{2m+1}\times \GL_1,\GSpin_4\times \GL_1,T(\HSpin_{4}^{+}\oplus \HSpin_{4}^{-}\otimes std_{\GL_1}),\GL_1\times \GL_1\times \GSpin_{2m-3}).$$
The nilpotent orbit $\iota$ induces a Bessel period on the unipotent radical of the parabolic subgroup $P=MU$ with $M=\GL_{1}^{m-2}\times \GSpin_5$ whose stabilizer in $M$ is $H$. In this case $\Delta_{red}$ is given by \eqref{22.5}. It is clear that Theorem \ref{main thm 3} holds in this case. The period integral is essentially the Gross-Prasad period for $\GSpin_{2m+1}\times \GSpin_4$ except that we replace the cusp form on $\GSpin_4$ by theta series. The unramified computation in \cite{II} proves Theorem \ref{main thm 1}.

By the discussion above, the strongly tempered quadruple associated to Table \ref{Table 22 non-red} is given as follows. Here for $\iota$, we only list the root type of the Levi subgroup $L$ of $G$ such that $\iota$ is principal in $L$ and
$$\ast=(\GSpin_{2m+1}\times \GL_1,\GSpin_4\times \GL_1,T(\HSpin_{4}^{+}\oplus \HSpin_{4}^{-}\otimes std_{\GL_1})).$$

\begin{figure}[h!]
\begin{tabular}{| c | c | c |}
\hline
$(G,H,\rho_H)$ & $\iota$ & $\hat{\rho}$ \\
\hline

$(\GSO_{8},S(\GL_2\times \GSO_4),T(std_{\GL_2}\oplus std_{\GL_2}))$ & $A_1$ & $T(std_{\Spin_8})\oplus  T(\HSpin_8)$ \\

\hline

$(\GL_{m}\times \GL_n,\GL_n\times \GL_n,T(std_{\GL_n}\otimes std_{\GL_n}))$ & $A_{m-n+1}$ & $T(std_{\SL_m}\otimes std_{\SL_n})\oplus T(std_{\SL_n})$ \\

\hline

$(\GL_{m}\times \GL_n,\GL_m\times \GL_{m+1},T(std_{\GL_m}\otimes std_{\GL_{m+1}}))$ & $A_{n-m-2}$ & $T(std_{\SL_m}\otimes std_{\SL_n})\oplus T(std_{\SL_n})$\\

\hline

$(\GL_n,\GL_2, T(std_{\GL_2})$ & $A_{n-3}$ & $T(std_{\SL_n})\oplus T(std_{\SL_n})$ \\

\hline

$\ast$ & $B_{m-2}$ & $T(std_{\Sp_{2m}})\oplus T(std_{\Sp_{2m}})$ \\

\hline

\end{tabular}
\captionof{table}{Dual quadruples of Table \ref{Table 22 non-red}}
\label{Table 22 non-red list}
\end{figure}

\newpage

\section{Summary}
We summarize our findings in this paper into the following 6 tables.

\begin{itemize}
\item Table \ref{red list} contains reductive strongly tempered quadruples for which we have provided evidence for Conjecture \ref{BSV conj}(1) and (2) (i.e., Theorem \ref{main thm 1} and \ref{main thm 2}).
\item Table \ref{red list extra} contains the remaining reductive strongly tempered quadruples. For all of them except $(\GL_6\times \GL_2,\GL_2\times S(\GL_4\times \GL_2),\wedge^2\otimes std_{\GL_2})$, we have provided evidence for Conjecture \ref{BSV conj}(1) (i.e. Theorem \ref{main thm 1}).
\item Table \ref{non-red list 1} contains non-reductive strongly tempered quadruples for which we have provided evidence for Conjecture \ref{BSV conj}(1) and (2) (i.e., Theorem \ref{main thm 1}, \ref{main thm 2} and \ref{main thm 3}).
\item Table \ref{non-red list 1 extra} contains non-reductive strongly tempered quadruples for which we have provided evidence only for Conjecture \ref{BSV conj}(1) (i.e., Theorem \ref{main thm 1} and \ref{main thm 3}).
\item Table \ref{non-red list 2} contains non-reductive strongly tempered quadruples for which we have provided evidence for Conjecture \ref{BSV conj}(1) by assuming Conjecture \ref{Whittaker induction} and we have provided evidence for Conjecture \ref{BSV conj} (2) (i.e. Theorem \ref{main thm 2} and \ref{main thm 3}).
\item Table \ref{non-red list 2 extra} contains the remaining non-reductive strongly tempered quadruples. For each of them, we have only provided evidence for Conjecture \ref{BSV conj}(1) by assuming Conjecture \ref{Whittaker induction} (i.e., Theorem \ref{main thm 3}).
\end{itemize}

For quadruples $(G,H,\rho_H,\iota)$ in Table \ref{red list} and \ref{red list extra}, the nilpotent orbit $\iota$ is trivial. For all the quadruples $\Delta=(G,H,\rho_H,,\iota)$ in Table \ref{red list}--\ref{non-red list 2 extra}, the dual quadruple is given by $(\hat{G},\widehat{G/Z_\Delta},\hat{\rho},1)$ where $\hat{\rho}$ is given in the tables and $Z_\Delta=Z_G\cap ker(\rho_H)$.

\newpage

\begin{figure}[h!]\leftskip-1cm
\begin{tabular}{| c | c  |c| c|}
\hline
\textnumero & (G, H, $\rho_H$)&  $\hat{\rho}$ \\
\hline
1  &($\SO_{2m+1}\times \SO_{2m}$, $\SO_{2m}$, 0)& $std_{\Sp_{2m}}\otimes std_{\SO_{2m}}$ \\
\hline

2  &($\SO_{2m+2}\times \SO_{2m+1}$, $\SO_{2m+1}$, 0)& $std_{\Sp_{2m}}\otimes std_{\SO_{2m+2}}$ \\
\hline

3  &($\GSp_6\times \GSpin_{7}$, $S(\GSp_6\times \GSpin_{7})$, $std_{\Sp_6}\otimes \Spin_7$) & $std_{\Sp_6}\otimes \Spin_7$ \\

\hline

4  &$(\GSp_6\times \GSpin_{9},S(\GSp_6\times \GSpin_{8}),std_{\Sp_6}\otimes \HSpin_8)$  & $std_{\Sp_8}\otimes \Spin_7$\\

\hline

5  &$(\GL_n\times \GL_n,\GL_n,T(std_{\GL_n}))$ & $T(std_{\GL_n}\otimes std_{\GL_n})$  \\
\hline

6  &$(\GL_{n+1}\times \GL_n,\GL_n,0)$ & $T(std_{\GL_{n+1}}\otimes std_{\GL_n})$  \\
\hline

7  &$(\GSp_4\times \GL_2,G(\SL_2\times \SL_2),T(std_{\GL_2,2}))$ & $T(Std_{\GSp_4}\otimes Std_{\GL_2})$ \\

\hline

8  &$(\GSp_4\times \GL_3,H=G,T(std_{\GSp_4}\otimes std_{\GL_3}))$ & $T(Std_{\GSp_4}\otimes Std_{\GL_3})$\\

\hline

9  &$(\GSp_4\times \GL_4, S(\GSp_4\times \GL_4),std_{\Sp_4}\otimes \wedge^2\oplus T(std_{\GL_4}))$ & $T(Std_{\GSp_4}\otimes Std_{\GL_4})$  \\

\hline

10  & $(\GSp_4\times \GL_5, S(\GSp_4\times \GL_4),std_{\Sp_4}\otimes \wedge^2)$ & $T(Std_{\GSp_4}\otimes Std_{\GL_5})$ \\

\hline

11  & $(\GSpin_7\times \GL_3,\GSpin_6\times \GL_3,T(\HSpin_6\otimes std_{\GL_3}))$ & $T(Std_{\GSp_6}\otimes Std_{\GL_3})$ \\

\hline

12 & $(\SO_{2m+1}\times \Sp_{2m},H=G, std_{\SO_{2m+1}}\otimes std_{\Sp_{2m}}\oplus std_{\Sp_{2m}})$ & $std_{\SO_{2m+1}}\otimes std_{\Sp_{2m}}\oplus std_{\Sp_{2m}}$  \\

\hline

13 & $(\SO_{2m+1}\times \Sp_{2m-2},\SO_{2m}\times \Sp_{2m-2}, std_{\SO_{2m}}\otimes std_{\Sp_{2m-2}})$ & $std_{\SO_{2m-1}}\otimes std_{\Sp_{2m}}\oplus std_{\Sp_{2m}}$  \\

\hline

14 & $(\GL_4\times \GSO_4,S(\GSp_4\times \GSO_4),std_{\SO_4}\times std_{\Sp_4})$ & $std_{\SL_2}\otimes std_{\SO_6}\oplus std_{\SO_6}\otimes std_{\SL_2}$  \\

\hline

15 & $(\GL_4\times\GL_2,\GL_2\times \GL_2,0)$ & $std_{\SL_2}\otimes \wedge^2\oplus T(std_{\SL_4})$ \\

\hline

16 & $(\GL_4\times \GSp_4, \GL_4\times \GSp_4, T(std_{\GL_4}\otimes std_{\GSp_4}))$ & $std_{\Sp_4}\otimes \wedge^2\oplus T(std_{\SL_4})$  \\

\hline

17 & $(\GSpin_7\times \GSpin_6,\GSpin_6\times \GSpin_6,T(\HSpin_6\otimes \HSpin_6))$ & $std_{\Sp_6}\otimes std_{\Spin_6}\oplus T(\HSpin_6)$  \\

\hline

18 & $(\GL_n\times \GL_n, \GL_n\times \GL_n, T(std_{\GL_n}\otimes std_{\GL_n}\oplus std_{\GL_n}))$ & $T(std_{\GL_n}\otimes std_{\GL_n})\oplus T(std_{\GL_n})$\\

\hline

19 &$(\GL_{n+1}\times \GL_n, \GL_n\times \GL_n, T(std_{\GL_n}\otimes std_{\GL_n}))$ & $T(std_{\GL_{n+1}}\otimes std_{\GL_n})\oplus T(std_{\GL_n})$ \\

\hline

20 &$(\GL_n\times \GL_{n-1}, \GL_n\times \GL_{n-1}, T(std_{\GL_n}\otimes std_{\GL_{n-1}}\oplus std_{\GL_{n}}))$ & $T(std_{\GL_{n-1}}\otimes std_{\GL_n})\oplus T(std_{\GL_n})$ \\

\hline

21 &$(\GL_n\times \GL_{n-2}, \GL_{n-1}\times \GL_{n-2}, T(std_{\GL_{n-1}}\otimes std_{\GL_{n-2}}))$ & $T(std_{\GL_{n}}\otimes std_{\GL_{n-2}})\oplus T(std_{\GL_n})$  \\

\hline

22 & $((\GL_2)^5,S(\GL_2\times \GL_2\times \GL_2),std_{\GL_2}\otimes std_{\GL_2}\otimes std_{\GL_2})$ & $\ast$\\

\hline

23 & $\sharp$ & $\ast\ast$\\

\hline

24 & $(\GL_2\times \GL_2\times \GL_2,\GL_2\times \GL_2,T(std_{\GL_2}\otimes std_{\GL_2}))$ & $\ast\ast\ast$ \\

\hline

25 & $\sharp\sharp$ & $\ast\ast\ast\ast$\\

\hline

26 & $(\GL_2\times \GL_1,\GL_2\times \GL_1,T(std_{\GL_2}\oplus std_{\GL_2}\otimes std_{\GL_1}))$ & $T(std_{\GL_2}\oplus std_{\GL_2}\otimes std_{\GL_1})$ \\
\hline

\end{tabular}
\captionof{table}{Reductive strongly tempered quadruples 1}
\label{red list}
\end{figure}

$$\sharp= ((\GL_2)^4,S(\GL_2\times \GL_2\times \GL_2),std_{\GL_2}\otimes std_{\GL_2}\otimes std_{\GL_2}\oplus T(std_{\GL_2,2})).$$

$$\sharp\sharp=((\GL_2)^3,S(\GL_2\times \GL_2\times \GL_2),std_{\GL_2}\otimes std_{\GL_2}\otimes std_{\GL_2}\oplus T(std_{\GL_2,2}\oplus T(std_{\GL_2,3})).$$

$$\ast=std_{\GL_2,1}\otimes std_{\GL_2,2}\otimes std_{\GL_2,3}\oplus std_{\GL_2,1}\otimes std_{\GL_2,4}\otimes std_{\GL_2,5}.$$

$$\ast\ast=T(std_{\GL_2,1}\otimes std_{\GL_2,2})\oplus std_{\GL_2,1}\otimes std_{\GL_2,3}\otimes std_{\GL_2,4}.$$

$$\ast\ast\ast=T(std_{\GL_2,1})\oplus std_{\GL_2,1}\otimes std_{\GL_2,2}\otimes std_{\GL_2,3}.$$

$$\ast\ast\ast\ast=T(std_{\GL_2,1}\otimes std_{\GL_2,2})\oplus T(std_{\GL_2,1}\otimes std_{\GL_2,3}).$$

\newpage

\begin{figure}[h!]\leftskip-1cm
\begin{tabular}{| c | c  |c| c|}
\hline
\textnumero & (G, H, $\rho_H$)&  $\hat{\rho}$ \\
\hline

1 &($\GSp_6\times \GSp_4$, $G(\Sp_4\times \Sp_2)$,0) & $\Spin_5\otimes \Spin_7$\\

\hline

2 & $(\GSp_6\times \GSO_4,S(\GSO_4\times G(\Sp_4\times \SL_2)), std_{\SO_4}\times std_{\Sp_4})$ & $std_{\SL_2}\otimes \Spin_7\oplus \Spin_7 \otimes std_{\SL_2}$    \\

\hline

3 & $\ast$ & $std_{\Sp_4}\otimes std_{\Spin_8}\oplus \HSpin_8\otimes std_{\SL_2}$   \\

\hline

4 & $(\GL_6\times \GL_2,\GL_2\times S(\GL_4\times \GL_2),\wedge^2\otimes std_{\GL_2})$ & $\wedge^3\oplus T(std_{\SL_6}\otimes std_{\SL_2} )$ \\

\hline

5 &$\ast \ast$ & $std_{\SL_2}\otimes \wedge^2\oplus T(std_{\SL_4}\otimes std_{\SL_2})$ \\

\hline

6 &$(\GL_{2m},\GL_m\times \GL_m,T(std_{\GL_m}))$ & $T(\wedge^2)\oplus T(std_{\GL_{2m}})$  \\

\hline

7 &$(\GL_{2m+1},\GL_{m+1}\times \GL_m,T(std_{\GL_{m+1}}))$ & $T(\wedge^2)\oplus T(std_{\GL_{2m+1}})$ \\

\hline

8 &$(\GL_3,\GL_2\times \GL_1,T(std_{\GL_2}))$ & $T(std_{\SL_3})\oplus T(std_{\SL_3})$ \\

\hline

9 &$(\GSpin_5\times \GL_1,\GSpin_4\times \GL_1,T(\HSpin_{4}^{+}\oplus \HSpin_{4}^{-}\otimes std_{\GL_1}))$ & $T(std_{\Sp_{4}})\oplus T(std_{\Sp_{4}})$  \\

\hline

\end{tabular}
\captionof{table}{Reductive strongly tempered quadruples 2}
\label{red list extra}
\end{figure}

$$\ast=(\GSp_4\times \GSpin_8\times \GL_2,S(\GSpin_8\times G(\Sp_4\times \SL_2)), std_{\Sp_4}\otimes std_{\Spin_8}\oplus \HSpin_8\otimes std_{\SL_2}).$$

$$\ast\ast=(\GL_2\times \GL_4\times \GL_2,S(\GL_2\times \GL_4)\times \GL_2, std_{\GL_2}\otimes \wedge^2\oplus T(std_{\GL_4}\times std_{\GL_2})).$$

\begin{figure}[h!]\leftskip-2cm
\begin{tabular}{| c| c | c | c |}
\hline
\textnumero & $(G,H,\rho_H)$ & $\iota$ &$\hat{\rho}$\\
\hline

 1& $(\SO_{2m+1}\times \SO_{2n},\SO_{2n},0)$ &$B_{m-n}$ & $std_{\Sp_{2m}}\otimes std_{\SO_{2n}}$  \\
\hline

 2&$(\SO_{2m+1}\times \SO_{2n},\SO_{2m+1},0)$ &$D_{n-m}$ & $std_{\Sp_{2m}}\otimes std_{\SO_{2n}}$ \\
\hline

 3& $(\GSpin_{2m+1}\times \GSp_6,S(\GSpin_8\times \GSp_6),std_{\Sp_6}\otimes \HSpin_8)$ &$B_{m-4}$ & $std_{\Sp_{2m}}\otimes \Spin_7$  \\

\hline

 4&$(\SO_{2m+1},\SO_2,0)$ & $B_{m-1}$ & $T(std_{\Sp_{2n}})$ \\
\hline

 5&$(\GSpin_{2m+1}\times \GL_2, G(\SL_2\times \SL_2), T(std_{\GL_2}))$ & $B_{m-2}$ & $T(Std_{\GSp_{2m}}\otimes Std_{\GL_2})$ \\

\hline

 6&$(\GSpin_{2m+1}\times \GL_3,\GSpin_6\times \GL_3,T(\HSpin_6\otimes std_{\GL_3}))$ & $B_{m-3}$ & $T(Std_{\Sp_{2m}}\otimes Std_{\SL_3})$ \\

\hline

 7&$(\SO_{2m+1}\times \Sp_{2n-2},\SO_{2n}\times \Sp_{2n-2},  std_{\SO_{2n}}\otimes std_{\Sp_{2n-2}})$ & $B_{m-n}$ & $std_{\SO_{2n-1}}\otimes std_{\Sp_{2m}}\oplus std_{\Sp_{2m}}$\\

\hline

8&$(\GSpin_{2k}\times \GSO_4,S(\GSp_4\times \GSO_4),std_{\SO_4}\times std_{\Sp_4})$ & $D_{k-2}$ & $std_{\SL_2}\otimes std_{\SO_{2k}}\oplus std_{\SO_{2k}}\otimes std_{\SL_2}$ \\

\hline

9& $(\GSpin_{2m+1}\times \GSpin_6,\GSpin_6\times \GSpin_6,T(\HSpin_6\otimes \HSpin_6))$ & $B_{m-3}$ & $std_{\Sp_{2m}}\otimes std_{\SO_6}\oplus T(\HSpin_6)$ \\

\hline

\end{tabular}
\captionof{table}{Non-reductive strongly tempered quadruples 1}
\label{non-red list 1}
\end{figure}

\newpage

\begin{figure}[h!]\leftskip-2cm
\begin{tabular}{| c| c | c | c |}
\hline
\textnumero & $(G,H,\rho_H)$ & $\iota$ &$\hat{\rho}$\\
\hline

1&$(\GSp_6\times \GL_2,\GL_2,0)$  & $A_2$& $std_{\GL_2}\otimes \Spin_7$ \\

\hline

2&$(\GSp_8\times \GL_2,G(\SL_2\times \SL_2),0)$ &$A_2$ & $std_{\GL_2}\otimes \Spin_9$  \\
\hline

3&$(\GSp_{10},\GL_2,0)$ &$A_4$& $\Spin_{11}$  \\
\hline

4&$(\GSO_{12},\GL_2,0)$ &$A_5$ & $\HSpin_{12}$ \\
\hline

5&$(\GL_6,\GL_2,0)$ & $A_2\times A_2$ & $\wedge^3$  \\
\hline

6&$(E_{7},\PGL_2,0)$ & $E_6$ & $std_{E_7}$   \\

\hline

7&$(\GL_{2m},\GL_m,T(std_{\GL_m}))$ & $(A_1)^m$ & $T(\wedge^2)$ \\
\hline

8&$(\GL_{2m+1},\GL_m,0)$ & $(A_1)^m$ & $T(\wedge^2)$ \\
\hline

9&$(\GSpin_{2k},\GSpin_3,T(\Spin_3))$ & $D_{k-1}$ & $T(std_{\SO_{2k}})$ \\

\hline

10&$(\GSp_6,\GL_2,T(std_{\GL_2}))$ & $A_2$ & $T(\Spin_7)$ \\

\hline

11&$(\GSp_8,G(\SL_2\times \SL_2),T(std_{\GL_2}))$ & $A_2$ & $T(\Spin_9)$ \\

\hline

12&$(GE_6, \GL_3,T(std_{\GL_3}))$ & $D_4$ & $T(std_{E_6})$ \\
\hline

13& $\ast$ & $B_{m-2}$ & $T(std_{\Sp_{2m}})\oplus T(std_{\Sp_{2m}})$ \\

\hline

\end{tabular}
\captionof{table}{Non-reductive strongly tempered quadruples 2}
\label{non-red list 1 extra}
\end{figure}

$$\ast=(\GSpin_{2m+1}\times \GL_1,\GSpin_4\times \GL_1,T(\HSpin_{4}^{+}\oplus \HSpin_{4}^{-}\otimes std_{\GL_1})).$$

\begin{figure}[h!]\leftskip-2cm
\begin{tabular}{|c| c | c | c |}
\hline
\textnumero & $(G,H,\rho_H)$ & $\iota$ &$\hat{\rho}$\\
\hline

 1& $(\GL_m\times \GL_n,\GL_n,0)$ & $A_{m-n-1}$ & $T(std_{\GL_m}\otimes std_{\GL_n})$ \\

\hline

 2& $(\GSp_4\times \GL_n, S(\GSp_4\times \GL_4),std_{\Sp_4}\otimes \wedge^2)$ & $A_{n-5}$ & $T(Std_{\Sp_4}\otimes Std_{\SL_m})$ \\

\hline

 3& $(\SO_{2m+1} \times \Sp_{2k},\SO_{2m+1}\times \Sp_{2m},std_{\SO_{2m+1}}\otimes std_{\Sp_{2m}})$ & $C_{k-m}$ & $std_{\SO_{2k+1}}\otimes std_{\Sp_{2m}}\oplus std_{\Sp_{2m}}$\\

\hline

 4 & $(\GL_{m}\times \GL_n,\GL_n\times \GL_n,T(std_{\GL_n}\otimes std_{\GL_n}))$ & $A_{m-n+1}$ & $T(std_{\SL_m}\otimes std_{\SL_n})\oplus T(std_{\SL_n})$ \\

\hline

 5 & $(\GL_{m}\times \GL_n,\GL_m\times \GL_{m+1},T(std_{\GL_m}\otimes std_{\GL_{m+1}}))$ & $A_{n-m-2}$ & $T(std_{\SL_m}\otimes std_{\SL_n})\oplus T(std_{\SL_n})$\\

\hline

\end{tabular}
\captionof{table}{Non-reductive strongly tempered quadruples 3}
\label{non-red list 2}
\end{figure}

\newpage

\begin{figure}[h!]\leftskip-2cm
\begin{tabular}{|c| c | c | c |}
\hline
\textnumero & $(G,H,\rho_H)$ & $\iota$ &$\hat{\rho}$\\
\hline

1& $(\GSp_{12},\GSp_4,0)$ &$A_2\times A_2$ & $\Spin_{13}$   \\
\hline

2& $(\PGSO_{10},\GL_2,0)$ & $A_3$ & $T(\HSpin_{10})$ \\

\hline

3& $(\GSO_{12},S(\GSp_4\times \GSO_4),0)$ & $A_1\times A_1$ & $\HSpin_{12}^{+}\oplus \HSpin_{12}^{-}$  \\

\hline

4& $(\GSO_{12}\times \PGL_2,S(\GL_2\times \GSO_4),0)$ & $A_3$ &  $std_{\SL_2}\otimes std_{\Spin_{12}}\oplus \HSpin_{12}$ \\

\hline

5& $\ast$ & $A_1$ & $std_{\Sp_4}\otimes std_{\Spin_{12}}\oplus \HSpin_{12}$ \\

\hline

6& $(\GSO_{8}\times \GSO_4,S(\GL_2\times \GSO_4),0)$ & $A_1$ & $std_{\SL_2}\otimes std_{\Spin_8}\oplus \HSpin_8 \otimes std_{\SL_2}$ \\

\hline

7& $(\GSpin_7, S(\GL_2\times \GL_2), std_{\GL_2})$ & $A_1$ & $\wedge^3\oplus std_{\Sp_6}$ \\

\hline

8& $(\GSO_{12},S(\GL_2\times \GSO_4),T(std_{\GL_2}))$ & $A_3$ & $\HSpin_{12}\oplus T(std_{\Spin_{12}})$   \\
\hline

9& $\ast\ast$ & $A_1$ & $std_{\SL_2}\otimes std_{\Spin_{10}}\oplus T(std_{\Spin_{10}})$  \\

\hline

10& $(\GSO_{8}\times \GL_2,S(\GL_2\times \GSO_4),T(std_{\GL_2}))$ & $A_1$ & $std_{\SL_2}\otimes std_{\Spin_{8}}\oplus T(std_{\Spin_{8}})$  \\

\hline

11& $(\GL_6,\GL_2\times \GL_2,0)$ & $A_1\times A_1$ & $\wedge^3\oplus T(std_{\SL_6} )$  \\

\hline

12 & $(\GSO_{8},S(\GL_2\times \GSO_4),T(std_{\GL_2}\oplus std_{\GL_2}))$ & $A_1$ & $T(std_{\Spin_8})\oplus  T(\HSpin_8)$ \\

\hline

13 & $(\GL_n,\GL_2, T(std_{\GL_2})$ & $A_{n-3}$ & $T(std_{\SL_n})\oplus T(std_{\SL_n})$ \\

\hline

\end{tabular}
\captionof{table}{Non-reductive strongly tempered quadruples 4}
\label{non-red list 2 extra}
\end{figure}

$$\ast=(\GSpin_4\times \GSpin_{12}, S(\GSpin_8\times G(\Sp_4\times \SL_2)), std_{\Sp_4}\otimes std_{\Spin_8}).$$

$$\ast\ast=(\GSpin_{10}\times \GL_2, S(\GL_2\times \GSpin_6)\times \GL_2, T(\HSpin_6\otimes std_{\GL_2})).$$

\newpage

\section{Table S of \cite{K}}\label{sec table S}
In this section we will discuss Table S of \cite{K}. This table contains several models from Table 1, 2, 11, 12, 22 of \cite{K} with some $A_1$ components (i.e. $\hat{G}=\hat{G}_1\times \hat{G_2}$ where $\hat{G}_1$ is of Type $A_1$). Then one can obtain more multiplicity free representations by gluing those representations together via the $A_1$-components. Two of the models in Table S contain two $A_1$-components (i.e. (S.1) and (S.2)) and hence it can be used to create infinite many multiplicity free representations. We refer the reader to \cite{K} for the details of the gluing process. In this section we will discuss how to write down the dual of those models by defining a gluing process for the dual of models in Table S. 

We first recall Table S from \cite{K}. Note that underlined section of the $\SL_2$ part is where we can glue the representations. Model (S.1) and (S.2) contains two underlined $\SL_2$ and we can use it glue representations with any arbitrary length.

\begin{figure}[h!]
	\begin{tabular}{|c| c | c| }
		\hline
		\textnumero in \cite{K} & $\hat{G}$ & $\hat{\rho}$\\
		\hline
		
(S.1) & $\underline{\SL}_2\times \Sp_{2m}\times \underline{\SL}_2$   & $std_{\SL_2}\otimes std_{\Sp_{2m}}\otimes std_{\SL_2}$ \\
\hline

(S.2) &  $\underline{\SL}_2\times \Spin_{8}\times \underline{\SL}_2$  & $std_{\SL_2}\otimes std_{\Spin_{8}}\oplus \HSpin_8\otimes  std_{\SL_2}$ \\
\hline

(S.3) & $\SO_n\times \underline{\SL}_2$  & $std_{\SO_n}\otimes std_{\SL_2}$ \\
\hline

(S.4) & $\Spin_{12}\times \underline{\SL}_2$  & $\HSpin_{12}\oplus std_{\Spin_{12}}\otimes std_{\SL_2}$ \\
\hline

(S.5) & $\Spin_{9}\times \underline{\SL}_2$  & $\Spin_{9}\otimes std_{\SL_2}$ \\
\hline

(S.6) & $\Spin_{8}\times \underline{\SL}_2$  & $T(std_{\Spin_8})\oplus \HSpin_{8}\otimes std_{\SL_2}$ \\
\hline

(S.7) &   $\Spin_{7}\times \underline{\SL}_2$  & $\Spin_{7}\otimes std_{\SL_2}$  \\
\hline

(S.8) & $\SL_2\times \Spin_{7}\times \underline{\SL}_2$  & $std_{\SL_2}\otimes std_{\Spin_{7}}\oplus \Spin_7\otimes  std_{\SL_2}$ \\
\hline

(S.9) & $\underline{\SL}_2$  & $std_{\SL_2}$ \\
\hline

(S.10) & $\underline{\SL}_2$  & $T(std_{\SL_2})$ \\
\hline

(S.11) & $\SL_m\times \underline{\SL}_2$  & $T(std_{\SL_m}\otimes std_{\SL_2})$ \\
\hline

(S.12) & $\SL_4\times \underline{\SL}_2$  & $T(std_{\SL_4})\oplus \wedge^2\otimes std_{\SL_2}$  \\
\hline

(S.13) & $\Sp_{2m}\times \underline{\SL}_2$  & $std_{\Sp_{2m}}\otimes \Sym^2$ \\
\hline

(S.14) & $\Sp_{2m}\times \underline{\SL}_2$  & $T(std_{\Sp_{2m}}\otimes std_{\SL_2})$ \\
\hline

(S.15) & $\Spin_{5}\times \underline{\SL}_2$  & $\Spin_5\oplus std_{\Spin_5}\otimes std_{\SL_2}$ \\
\hline

(S.16) & $G_2\times \underline{\SL}_2$  & $std_{G_2}\otimes std_{\SL_2}$ \\
\hline

\end{tabular}
\captionof{table}{Table S of \cite{K}}
\label{Knop Table S}
\end{figure}

It is clear that if one glue some anomaly-free representations with some non anomaly-free representations in Table S, one will get a non anomaly-free representation. Hence we can consider them separately. We first consider the anomaly-free representations in Table S, this corresponds to Model (S.1), (S.2), (S.3) when $n$ is even, (S.4)-(S.7), (S.10)-(S.12) and (S.14). For each of them we have already write down its dual quadruple in the previous sections. We just need to describe how to glue the dual together \footnote{(S.1) is Model 1 of Table \ref{non-red list 1} with n=4, (S.2) is Model 6 of Table \ref{non-red list 2 extra}, (S.3) when $n$ is even is Model 2 of Table \ref{non-red list 1} with $m=1$, (S.4) is Model 4 of Table \ref{non-red list 2 extra}, (S.5) is Model 2 of Table \ref{non-red list 1 extra}, (S.6) is Model 10 of Table \ref{non-red list 2 extra}, (S.7) is Model 1 of Table \ref{non-red list 1 extra}, (S.10) is Model 6 of Table \ref{red list} with $n=1$, (S.11) is Model 1 of Table \ref{non-red list 2} with $n=2$, (S.12) is Model 15 of Table \ref{red list}, and (S.14) is Model 5 of Table \ref{non-red list 1}.}. 

Let $\hat{\Delta}=(\hat{G},\hat{G},\hat{\rho},1)$ be one of such model and let $\Delta=(G,H,\rho_H,\iota)$ be its dual. If the model is not (S.1), (S.2) or (S.10), then $\hat{G}=\hat{G}_1\times \hat{G}_2$ and $G=G_1\times G_2$ with $\hat{G_1}, G_1$ being of Type $A_1$.  Moreover, by our description of $\Delta$ in the previous section, we know that the projection map $H\rightarrow G_1$ is surjective and we can write the group $HG_1$ (i.e. the group generated by $H$ and $G_1$) as $G_1\times H_1 \times H_2$ with $H_2=H\cap G_2$ and $H_1\simeq G_1$ is in the centralizer of $H_2$ in $G_1H\cap G_2$. Moreover  the image of the diagonal embedding from $H_1$ into $G_1\times H_1$ belongs to $H$. In particular, the representation $\rho_H$ induces a representation (still denoted by $\rho_{H}$) on $H_1\times H_2$ (on $H_2$-part this is given by restriction and on $H_1$ part it is given by restriction and the diagonal embedding from $H_1$ into $G_1\times H_1$). Finally, the nilpotent orbit $\iota$ is the product of some nilpotent oribt of $G_2$ with the trivial nilpotent orbit of $G_1$. For example, consider Model (S.3) when $n=4$, the dual is the trilinear $\GL_2$ model
$$(G,H,\rho_H,\iota)=(\PGL_{2}^3,\PGL_2,0,1)$$
and in this case 
$$G_1=\{(1,1,h)|\; h\in \PGL_2\},\; H=\{(h,h,h)|\; h\in \PGL_2\},\; H_1=\{(h,h,1)|\; h\in \PGL_2\},\; H_2=\{1\}.$$

If the Model is (S.1) or (S.2), then $\hat{G}=\hat{G}_{11}\times \hat{G}_{12}\times \hat{G}_2$ and $G=G_{11}\times G_{12}\times G_2$ with $\hat{G}_{11},\hat{G}_{12}, G_{11},G_{12}$ being of Type $A_1$.  Moreover, by our description of $\Delta$ in the previous section, we know that the projection map $H\rightarrow G_{11}\times G_{12}$ is surjective and we can write the group $HG_{11}G_{12}$ as $G_{11}\times G_{12}\times H_{11}\times H_{12}\times H_2$ with $H_2\subset H$, $H_{1i}\simeq G_{1i}$,  and the image of the diagonal embedding from $H_{1i}$ into $G_{1i}\times H_{1i}$ belongs to $H$ for $i=1,2$. Also the representation $\rho_H$ would be $0$ in this case. Finally, the nilpotent orbit $\iota$ is the product of some nilpotent oribt of $G_2$ with the trivial nilpotent orbit of $G_{11}\times G_{12}$.

Lastly, if the model is $(S.10)$, then $\hat{\Delta}=(\SL_2,\SL_2,T(std),1)$ and $\Delta=(\PGL_2,\GL_1,0,1)$. In particular $\hat{G}=\hat{G_1}$ and $G=G_1$ are of Type $A_1$.

Now we can describe the gluing process on the dual side. Suppose we are gluing two representations $(\hat{G},\hat{G},\hat{\rho})$ and $(\hat{G}',\hat{G}',\hat{\rho}')$. In particular we can write 
$$\hat{G}=\hat{G}_1\times \hat{G_2},\; \hat{G}'=\hat{G}_1'\times \hat{G}_2'$$
and we are gluing $\hat{G}_1$ with $\hat{G}_1'$. The goal is to write down the dual of 
$$\hat{\Delta}_{glue}=(\hat{G}_2\times \hat{G_1}\times \hat{G}_2',\hat{G}_2\times \hat{G_1}\times \hat{G}_2',\hat{\rho}\oplus \hat{\rho}',1).$$
Here we consider $\hat\rho$ (or $\hat{\rho}'$) as a representation of
	$\hat{G}_2\times \hat{G_1}\times \hat{G}_2'$ where the $\hat{G}_2'$ (or $\hat{G}_2$) component acts trivially.
Let $\Delta=(G,H,\rho_H,\iota)$ and $\Delta'=(G',H',\rho_H',\iota')$ be the dual of $\hat{\Delta}$ and $\hat{\Delta}'$.

There are two cases. First we consider the case when both representations are not (S.10). In this case, by our discussion above, we have the decomposition \footnote{If we are in the case of (S.1) or (S.2), then $G_1$ would be one of the $G_{1i}$ and $H_1$ would be the corresponding $H_{1i}$}
$$G_1H=G_1\times H_1\times H_2,\;G_1'H'=G_1'\times H_1'\times H_2'.$$
Then the dual would be given by 
$$\Delta_{glue}=(G_2\times G_1\times G_2', H_2\times H_1\times G_1\times H_1'\times H_2', \rho_H\oplus \rho_H'\oplus \rho',\iota\times \iota')$$
where $\rho'$ is the tensor product representation of $H_1\times G_1\times H_1'$. Note that when the model is not (S.1), (S.2) or (S.10), we have explained how to view $\rho_H$ (resp. $\rho_H'$) as a representation of $H_1\times H_2$ (resp. $H_1'\times H_2'$). In the cases of $(S.1)$ or $(S.2)$ the representation $\rho_H$ (resp. $\rho_H'$) is just 0. Also note that $\iota$ (resp. $\iota'$) is the product of some nilpotent orbit of $G_2$ (resp. $G_2'$) with the trivial nilpotent orbit of $G_1=G_1'$ and hence we can view $\iota\times \iota'$ as a nilpotent orbit of $G_1\times G_1\times G_2'$. Model \ref{GL(2) 5} and \ref{GL(2) 4} are examples of this case.

Now let's prove Theorem \ref{main thm S} in this case. The idea is similar to the proof of Conjecture \ref{BSV conj} for the model \ref{GL(2) 5}. Roughly speaking, we will show that the period integral associated to $\Delta_{glue}$ (resp. $\hat{\Delta}_{glue}$) is a combination of the period integrals associated to $\Delta,\Delta',\Delta_1$ (resp. $\hat{\Delta},\hat{\Delta}', \hat{\Delta}_1$) \footnote{$(\Delta_1,\hat{\Delta}_1)$ is just Model 2 in Table \ref{red list} when $m=1$, the period integral associated to $\Delta_1$ is the theta correspondence of $\Sp_2\times \SO_4$ while the period integral associated to $\hat{\Delta}_1$ is the triple product integral}
where
$$\Delta_1=(\SL_{2}^{3},\SL_{2}^{3},std\otimes std\otimes std,1), \hat{\Delta}_1=(\PGL_{2}^3,\PGL_2,0,1).$$
In particular, Conjecture \ref{BSV conj} for $(\Delta_{glue},\hat{\Delta}_{glue})$ will follow from Conjecture \ref{BSV conj} for $(\Delta,\hat{\Delta})$, $(\Delta',\hat{\Delta}')$ and $(\Delta_1,\hat{\Delta}_1)$. As Conjecture \ref{BSV conj} is know for $(\Delta_1,\hat{\Delta}_1)$ (by the Rallis inner product formula and the work of Harris-Kudla \cite{HK} for the triple product period), we know that Conjecture \ref{BSV conj} for $(\Delta_{glue},\hat{\Delta}_{glue})$ will follow from Conjecture \ref{BSV conj} for $(\Delta,\hat{\Delta})$ and $(\Delta',\hat{\Delta}')$. This proves Theorem \ref{main thm S}.

It remains to explain why the the period integral associated to $\Delta_{glue}$ (resp. $\hat{\Delta}_{glue}$) is a combination of the period integrals associated to $\Delta,\Delta',\Delta_1$ (resp. $\hat{\Delta},\hat{\Delta}', \hat{\Delta}_1$). We start with the period integral associated to $\Delta_{glue}$. In these case we start with an automorphic form $\phi=\phi_{G_2}\phi_{G_1}\phi_{G_2'}$ on $G_2\times G_1\times G_2'$. We first integrate over $G_1$ (note that the projection of nilpotent orbit $\iota\times \iota'$ to $G_1$ is the trivial orbit, so the unipotent integral associated to $\iota\times \iota'$ commutes with the integral over $G_1$). Since the symplectic representation $\rho_H$ (resp. $\rho_{H}'$) is on the group $H_1\times H_2$ (resp. $H_1'\times H_2'$), the integral over $G_1$ is just the integral of $\phi_{G_1}$ with theta function associated to $\rho'$ (recall that $\rho'$ is the tensor product representation of $H_1\times G_1\times H_1'$). By the theta correspondence of $\Sp_2\times \SO_4$, the integral
$$\int_{G_1(k)\back G_1(\BA)} \phi_{G_1}(g_1)\Theta_{\rho'}(h_1,g_1,h_1')dg_1$$
gives an automorphic form $\phi_{H_1\times H_1'}(h_1,h_1')$ in the irreducible space $\pi\otimes\pi$ on $H_1\times H_1'$ (assuming $\phi\in \pi$ of $G_1\simeq H_1\simeq H_2$).
We may as well assume $ \phi_{H_1\times H_1'}(h_1,h_1')$ has the form $\phi_{H_1}(h_1)\phi_{H_1'}(h_1')$. Note that by Rallis inner product formula $\|\phi_{G_1}\|``=" \|\phi_{H_1}\|\|\phi_{H_1'}\|$.
 Then the remaining integrals (i.e. the unipotent integral associated to $\iota\times \iota'$ and integral over $H_2\times H_1\times H_1'\times H_2'$) become the product of the period integrals of the automorphic forms $\phi_{G_2}\times \phi_{H_1}$ and $\phi_{H_1'}\times \phi_{G_2'}$ associated to the quadruples $(G_2\times H_1,H_2\times H_1,\rho_H,\iota)$ and $(H_1'\times G_2',H_1'\times H_2',\rho_H',\iota')$ respectively \footnote{note here for the embedding of $H_1\times H_2$ (respectively $H_1'\times H_2'$), the component $H_1$ (resp. $H_1'$) diagonally embeds into $G_2\times H_1$ (resp. $H_1'\times G_2'$)}. But these two quadruples are just $\Delta$ and $\Delta'$ via the isomorphism $H_1\simeq G_1\simeq G_1'\simeq H_1'$. As a result, Conjecture \ref{BSV conj}(1) for $\Delta_{glue}$ would follow from Conjecture \ref{BSV conj}(1) for $\Delta$ and $\Delta'$.

For the other direction, the period integral associated to $\hat{\Delta}_{glue}$ is given by

\begin{equation}\label{equation s.1}
\int_{\hat{G}_1(k)\back \hat{G}_1(\BA)}\int_{\hat{G}_2(k)\back \hat{G}_2(\BA)}\int_{\hat{G}_2'(k)\back \hat{G}_2'(\BA)} \phi_{\hat{G}_1}(g_1) \phi_{\hat{G}_2}(g_2)\phi_{\hat{G}_2'}(g_2') \Theta_{\hat{\rho}}(g_1,g_2)\Theta_{\hat{\rho}'}(g_1,g_2')dg_2'dg_2dg_1.
\end{equation}
By Conjecture \ref{BSV conj}(2) for $\Delta$ and $\Delta'$, we know that
\begin{itemize}
\item the integral
$$\int_{\hat{G}_2(k)\back \hat{G}_2(\BA)}  \phi_{\hat{G}_2}(g_2) \Theta_{\hat{\rho}}(g_1,g_2)dg_2$$
is non-vanishing only if the Arthur parameter of $\phi_{\hat{G}_2}$ factors through $H_1\times H_2\rightarrow G_2$, i.e. it is the lifting of an Arthur packet $\Pi_{\hat{H}_1}\otimes \Pi_{\hat{H}_2}$ of $\hat{H}_1(\BA)\times \hat{H}_2(\BA)$. Moreover, if this is the case and the packet $\Pi_{\hat{H}_1}\otimes \Pi_{\hat{H}_2}$ is tempered, then the automorphic function
$$\phi_1(g_1):= \int_{\hat{G}_2(k)\back \hat{G}_2(\BA)}  \phi_{\hat{G}_2}(g_2) \Theta_{\hat{\rho}}(g_1,g_2)dg_2,\; g_1\in \hat{G}_1(\BA)$$
belongs to the packet $\Pi_{\hat{H}_1}$ (as $H_1\simeq G_1$ we can view $\Pi_{\hat{H}_1}$ as a packet for $\hat{G}_1$);
\item the integral
$$\int_{\hat{G}_2'(k)\back \hat{G}_2'(\BA)}  \phi_{\hat{G}_2'}(g_2') \Theta_{\hat{\rho}'}(g_1,g_2')dg_2'$$
is non-vanishing only if the Arthur parameter of $\phi_{\hat{G}_2'}$ factors through $H_1'\times H_2'\rightarrow G_2'$, i.e. it is the lifting of an Arthur packet $\Pi_{\hat{H}_1'}\otimes \Pi_{\hat{H}_2'}$ of $\hat{H}_1'(\BA)\times \hat{H}_2'(\BA)$. Moreover, if this is the case and the packet $\Pi_{\hat{H}_1'}\otimes \Pi_{\hat{H}_2'}$ is tempered, then the automorphic function
$$\phi_1'(g_1):= \int_{\hat{G}_2'(k)\back \hat{G}_2'(\BA)}  \phi_{\hat{G}_2'}(g_2') \Theta_{\hat{\rho}'}(g_1,g_2')dg_2',\;g_1\in \hat{G}_1(\BA)$$
belongs to the packet $\Pi_{\hat{H}_1'}$ (as $H_1'\simeq G_1$ we can view $\Pi_{\hat{H}_1'}$ as a packet for $\hat{G}_1$).
\end{itemize}

By the above two facts, the integral \eqref{equation s.1} becomes
$$\int_{\hat{G}_1(k)\back \hat{G}_1(\BA)} \phi_{\hat{G}_1}(g_1)\phi_1(g_1)\phi_1'(g_1)dg_1$$
which is exactly a triple product integral on $\hat H_1\times \hat G_1\times \hat H_1'$. In particular, Conjecture \ref{BSV conj}(2) for $\Delta_{glue}$ follows from the work of Harris-Kudla \cite{HK} for the triple product period and Conjecture \ref{BSV conj}(2) for $\Delta$ and $\Delta'$. This finishes the proof of Theorem \ref{main thm S} for this case.

Next we consider the case when at least one of the representation is (S.10). We may assume that $(\hat{G}',\hat{G}',\hat{\rho}')$ is $(S.10)$. Then we have the decomposition 
$$G_1H=G_1\times H_1\times H_2,\;  G'\simeq G_1.$$
The dual of
$$\hat{\Delta}_{glue}=(\hat{G}_2\times \hat{G}_1,\hat{G}_2\times \hat{G}_1,\hat{\rho}\oplus T(std_{G_1}),1).$$
would be given by
$$\Delta_{glue}=(G_2\times G_1,H_2\times H_1\times G_1,\rho_H\oplus T(\rho'), \iota)$$
where $\rho'$  is the tensor product representation of $H_1\times G_1$.  Model \ref{GL(2) 3} and \ref{GL(2) 1} are examples of this case.

Now let's prove Theorem \ref{main thm S} in this case. The idea is similar to the proof of Conjecture \ref{BSV conj} for the model \ref{GL(2) 3}. Roughly speaking, we will show that the period integral associated to $\Delta_{glue}$ (resp. $\hat{\Delta}_{glue}$) is a combination of the period integrals associated to $\Delta,\Delta_2$ (resp. $\hat{\Delta}, \hat{\Delta}_2$) where \footnote{$(\Delta_2,\hat{\Delta}_2)$ is just Model 5 in Table \ref{red list} when $n=2$, the period integral associated to $\Delta_2$ is the theta correspondence of $\GL_2\times \GL_2$ while the period integral associated to $\hat{\Delta}_2$ is the Rankin-Selberg integral of $\GL_2\times \GL_2$}
$$\Delta_2=(\GL_2\times \GL_2,\GL_2\times \GL_2,T(std\otimes std),1),\;\hat{\Delta}_2=(\GL_2\times \GL_2,\GL_2,T(std),1).$$

In particular, Conjecture \ref{BSV conj} for $(\Delta_{glue},\hat{\Delta}_{glue})$ will follow from Conjecture \ref{BSV conj} for $(\Delta,\hat{\Delta})$, and $(\Delta_2,\hat{\Delta}_2)$. As Conjecture \ref{BSV conj} is know for $(\Delta_2,\hat{\Delta}_2)$ (by the theory of Rankin-Selberg integral and the Rallis inner product formula), we know that Conjecture \ref{BSV conj} for $(\Delta_{glue},\hat{\Delta}_{glue})$ will follow from Conjecture \ref{BSV conj} for $(\Delta,\hat{\Delta})$. This proves Theorem \ref{main thm S}.

It remains to explain why the period integral associated to $\Delta_{glue}$ (resp. $\hat{\Delta}_{glue}$) is a combination of the period integrals associated to $\Delta,\Delta_2$ (resp. $\hat{\Delta}, \hat{\Delta}_2$). The argument is very similar to the previous case as well as the case of the Model \ref{GL(2) 3}, we will skip it here.

This completes the description of the dual BZSV quadruples  associated to representations glued from anomaly-free representations in Table S of \cite{K}, as well as the proof of Theorem~\ref{main thm S} for those cases.

It remains to consider the non anomaly-free representations in Table S of \cite{K}, which are  (S.3) when $n$ is odd, (S.8), (S.9), (S.13), (S.15) and (S.16). It is easy to see that if we glue the model (S.8), (S.13) or (S.15) with another model, then the representation we get is not anomaly-free. Hence we just need to consider  (S.3) when $n$ is odd, (S.9) and (S.16). There are 6 different cases.

If we glue (S.3) when $n$ is odd with (S.9), we get the model (11.11) of Table 11 with $m=1$, this has already been considered in Section \ref{sec Table 11}. If we glue (S.9) with itself, the representation we get is just $T(std)$ of $\SL_2$ which is model (S.10). For the remaining four cases, the generic stabilizer of the representation is  not connected \footnote{while in this paper we suggested the form of dual quadruples for some representations with non-connected stabilizer,  it is not clear to us what the dual should be for the cases here}.


\begin{thebibliography}{99}

\bibitem{BSV}
D. Ben-Zvi, Y. Sakellaridis and A. Venkatesh,
{\it Relative Langlands duality}. preprint

\bibitem{Beu2}
R. Beuzart-Plessis,
{\it A local trace formula for the Gan-Gross-Prasad conjecture for unitary groups: the archimedean case.} Ast\'erisque no. 418 (2020).

\bibitem{BFGT}
A. Braverman, M. Finkelberg, V. Ginzburg, and R. Travkin 
{\it Mirabolic Satake equivalence and supergroups.} 
Compos. Math.157(2021), no.8, 1724–1765.

\bibitem{BFT}
A. Braverman, M. Finkelberg, and R. Travkin
{\it Gaiotto conjecture for ${\rm Rep_q(GL(N-1 | N))}$.} 
arXiv:2107.02653, 2021.




\bibitem{BF90}
D. Bump and S. Friedberg,
\newblock The exterior square automorphic {$L$}-functions on {${\rm GL}(n)$}.
\newblock In {\em Festschrift in honor of {I}. {I}. {P}iatetski-{S}hapiro on
  the occasion of his sixtieth birthday, {P}art {II} ({R}amat {A}viv, 1989)},
  volume~3 of {\em Israel Math. Conf. Proc.}, pages 47--65. Weizmann,
  Jerusalem, 1990.

\bibitem{BG92a}
D. Bump and D. Ginzburg,
{\em Spin {$L$}-functions on symplectic groups.}
Internat. Math. Res. Notices, (8):153--160, 1992.

\bibitem{BG00}
D. Bump and D. Ginzburg,
{\it Spin {$L$}-functions on {${\rm GSp}_8$} and {${\rm GSp}_{10}$}.}
Trans. Amer. Math. Soc., 352(2):875--899, 2000.


\bibitem{BG92-ann}
D. Bump and D. Ginzburg,
{\it Symmetric square {$L$}-functions on {${\rm GL}(r)$}.}
Ann. of Math. (2), 136(1):137--205, 1992.

\bibitem{D}
Sanath Devalapurkar, {\it ku-theoretic spectral decompositions for spheres and projective spaces}. arXiv:2402.03995

\bibitem{FTY}
M. Finkelberg, R. Travkin, and R. Yang
{\it Gaiotto Conjecture for ${\rm Rep_q(F(4))}$.}
arXiv:2406.05521, 2024.



\bibitem{FW}
Tony Feng and Jonathan Wang. {\it Geometric Langlands duality for periods}, 2024. arXiv:2402.00180

\bibitem{GGP}
W. Gan, B. Gross and D. Prasad, 
{\it Symplectic local root numbers, central critical L values, and restriction problems in the representation theory of classical groups.} 
Sur les conjectures de Gross et Prasad. I. Ast\'erisque No. 346 (2012), 1–109. ISBN: 978-2-85629-348-5

\bibitem{GGP2}
W. Gan, B. Gross and D. Prasad, 
{\it Branching laws for classical groups: the non-tempered case.} Compositio Mathematica. 2020;156(11):2298-2367

\bibitem{G2}
D. Ginzburg, {\it On standard L-functions for $E_6$ and $E_7$.} J. Reine Angew. Math., 465:101-131, 1995.


\bibitem{Gin93}
D. Ginzburg,
{\it On the standard {$L$}-function for {$G_2$}.}
 Duke Math. J., 69(2):315--333, 1993.

\bibitem{G95}
D. Ginzburg,
{\it On spin {$L$}-functions for orthogonal groups.}
Duke Math. J., 77(3):753--798, 1995.

\bibitem{GH04}
D. Ginzburg and J. Hundley,
{\it Multivariable {R}ankin-{S}elberg integrals for orthogonal groups.}
 Int. Math. Res. Not., (58):3097--3119, 2004.

\bibitem{GJ}
W. Gan and B. Jun,
{\it Generalized Whittaker models as instances of relative Langlands duality}. arXiv:2309.08874

\bibitem{GQT}
W. T. Gan, Y. Qiu, and S. Takeda,
{\it The regularized Siegel-Weil formula (the second term identity) and the Rallis inner product formula}
Invent. Math. 198 (2014), no. 3, 739–831.


\bibitem{GJ1}
S. Gelbart and H. Jacquet, 
{\it A relation between automorphic representations of GL(2)
 and GL(3),} Ann. Sci. Ecole Normale Sup., 4e serie, 11 (1978), 471-552.


\bibitem{GJR01}
D. Ginzburg, D. Jiang, and S. Rallis,
{\it Nonvanishing of the central critical value of the third symmetric
power L-functions.} Forum Math., 13(1):109-132, 2001


\bibitem{Gross}
B. Gross,
{\it On the motive of a reductive group,} Invent.Math.130(2), 287–313.


\bibitem{GRS}
D. Ginzburg, S. Rallis, D. Soudry,
{\it L-functions for symplectic groups}. Bulletin de la S. M. F., tome 126, no 2 (1998), p. 181-244

\bibitem{HK}
M. Harris and  S. Kudla,
{\it On a conjecture of Jacquet.} Contributions to automorphic forms, geometry, and number theory. Johns Hopkins Univ. Press,
2004, 355–371


\bibitem{II}
A. Ichino and T. Ikeda,
{\it On the periods of automorphic forms on special orthogonal groups and the Gross--Prasad conjecture.} Geometric and Functional Analysis 19 (2010), no. 5, 1378-1425.

\bibitem{JPSS}
H. Jacquet, I. I. Piatetskii-Shapiro and J. A. Shalika,
{\it Rankin-Selberg Convolutions}, American Journal of Mathematics
Vol. 105, No. 2 (Apr., 1983), pp. 367-464.

\bibitem{JS90}
H. Jacquet and J. Shalika.
\newblock Exterior square {$L$}-functions.
\newblock In {\em Automorphic forms, {S}himura varieties, and {$L$}-functions,
  {V}ol.\ {II} ({A}nn {A}rbor, {MI}, 1988)}, volume~11 of {\em Perspect.
  Math.}, pages 143--226. Academic Press, Boston, MA, 1990.

\bibitem{K}
F. Knop,
{\it Classification of multiplicity free symplectic representations}. Journal of Algebra Volume 301, Issue 2, 531-553.

\bibitem{Knop}
F. Knop and B. Schalke, 
{\it The dual group of a spherical variety}. Trans. Moscow Math. Soc. 2017, 187-216.

\bibitem{KR}
S. Kudla and S. Rallis,
{\it A regularized Siegel–Weil formula: the first term identity.} 
Ann. Math. 140, 1--80 (1994)

\bibitem{LM}
E. Lapid and Z. Mao,
{\it A conjecture on Whittaker–Fourier coefficients of cusp forms}, Journal of Number Theory 146, 448-505.

\bibitem{L}
J.-S. Li,
{\it Nonvanishing theorems for the cohomology of certain arithmetic quotients}
J. Reine Angew. Math. 428 (1992), 177–217.

\bibitem{Lo}
I. Losev,
{\it Coisotropic representations of reductive groups.} Trudy Mosc. Mat. Ob-va,
66(2005), p. 157-181 (in Russian). English translation in: Trans. Moscow Math. Soc.
(2005), 143-168.

\bibitem{MWZ}
Z. Mao, C. Wan and L. Zhang,
{\it BZSV duality for some strongly tempered spherical varieties}. arXiv:2310.17837

\bibitem{PPS89}
S.~J. Patterson and I.~I. Piatetski-Shapiro,
{\it The symmetric-square {$L$}-function attached to a cuspidal
  automorphic representation of {${\rm GL}_3$}.}
Math. Ann., 283(4):551--572, 1989.

\bibitem{PS}
A. Pollack and S. Shah,
{\it Multivariate Rankin-Selberg integrals on $\GL_4$ and $\GU_{2,2}$}. Canadian Mathematical Bulletin, Vol. 61 (4), 2018, 822-835


\bibitem{S}
Y. Sakellaridis,
{\it Functorial transfer between relative trace formulas in rank one}. Duke Math. J., 170(2):279-364, 2021.

\bibitem{S1}
Y. Sakellaridis,
{\it Spherical functions on spherical varieties.} Amer. J. Math., 135(5):1291-1381, 2013. 



\bibitem{SV}
Y. Sakellaridis and A. Venkatesh, \emph{Periods and harmonic analysis
  on spherical varieties}, Ast\'{e}risque (2017), no.~396, viii+360.

\bibitem{Tk14}
S. Takeda,
{\it The twisted symmetric square {$L$}-function of {${\rm GL}(r)$}.}
 Duke Math. J., 163(1):175--266, 2014.
 
 \bibitem{TY}
 R. Travkin and R. Yang
{\it  Twisted Gaiotto equivalence for ${\rm GL(M | N)}$.}
arXiv:2306.09556, 2023

\bibitem{Wal2}
J.-L. Waldspurger,
{\it Une formule int\'egrale reli\'ee \`a la conjecture locale de Gross-Prasad, 2e partie : extension aux repr\'esentations temp\'er\'ees.}  in "Sur les conjectures de Gross et Prasad. I" Ast\'erisque No. 346 (2012), 171-312


\bibitem{Wan2}
C. Wan,
{\it Multiplicity One Theorem for the Ginzburg-Rallis Model: the tempered case.} Trans. Amer. Math. Soc. 371 (2019), 7949-7994.



\bibitem{WZ}
C. Wan and L. Zhang,
 {\it Periods of automorphic forms associated to strongly tempered spherical varieties}. Accepted by Memoir of AMS. arxiv 2102.03695, 109 pages.
 
 \bibitem{Y}
 S. Yamana,
{\it $L$-functions and theta correspondence for classical groups.}
 Invent. Math.  196 (2014), 651--732.




\end{thebibliography}
\end{document}